\documentclass[11pt, reqno]{amsart}

\usepackage{amscd,amsmath}
\usepackage{mathrsfs}
\usepackage{amsfonts}
\usepackage{amssymb}
\usepackage{enumerate}
\usepackage{setspace}
\usepackage{color}
\usepackage{esint}
\usepackage{dsfont}
\usepackage{booktabs}
\usepackage{multirow}
\usepackage[table,xcdraw]{xcolor}

\usepackage{tabu}
\usepackage{tikz} 
\usepackage{setspace}
\usepackage[margin=2.2cm]{geometry}

\usepackage[english]{babel}
\usepackage{graphicx}

\definecolor{luh-dark-blue}{rgb}{0.0, 0.313, 0.608}

\usepackage[colorlinks=true, citecolor=blue]{hyperref}

\newtheorem{theorem}{Theorem}[section]
\newtheorem{lemma}{Lemma}[section]

\newtheorem{corollary}{Corollary}[section]
\newtheorem{definition}{Definition}[section]

\newtheorem{remark}{Remark}[section]
\newtheorem{proposition}{Proposition}[section]

\newcommand{\eqn}{\begin{eqnarray}}
\newcommand{\een}{\end{eqnarray}}

\numberwithin{equation}{section}

\DeclareMathOperator{\dv}{div}

\begin{document}

\title[Oldroyd type models in hybrid  Besov spaces]{Global existence and long time behavior of solutions to some Oldroyd type models in hybrid  Besov spaces}

\author{Hantaek Bae}
\address{Department of Mathematical Sciences, Ulsan National Institute of Science and Technology, Republic of Korea}
\email{hantaek@unist.ac.kr}

\author{Jaeyong Shin}
\address{Department of Mathematical Sciences, Ulsan National Institute of Science and Technology, Republic of Korea}
\email{sinjaey@unist.ac.kr}

\date{\today}

\subjclass[2010]{35Q35; 76D03}

\keywords{Viscoelastic fluids; Hybrid Besov spaces; Existence; Decay rates}

\begin{abstract}
In this paper,  we deal with some Oldroyd type models, which describe incompressible viscoelastic fluids. There are 3 parameters in these models: the viscous coefficient of fluid $\nu_{1}$, the viscous coefficient of the elastic part of the stress tensor $\nu_{2}$, and the damping coefficient of the elastic part of the stress tensor $\alpha$. In this paper, we assume at least one of the parameters is zero: $(\nu_{1}, \nu_{2},\alpha)=(+,0,+), (+,0,0), (0,+,+), (0,+,0), (0,0,+)$ and prove the global existence of unique solutions to all these 5 cases in the framework of hybrid Besov spaces. We also derive decay rates of the solutions except for the case $(\nu_{1}, \nu_{2},\alpha)=(+,0,0)$. To the best of our knowledge, decay rates in this paper are the first results in this framework, and can improve some previous works.

\end{abstract}

\maketitle

\vspace{-5ex}

\section{Introduction}
In this paper, we investigate some incompressible viscoelastic fluids. We start from the equations of homogeneous incompressible fluids:
\[
\begin{split}
&u_{t}+u\cdot\nabla u=\dv \sigma, \\
&  \dv u=0,
\end{split}
\]
where $u$ is the velocity of the fluid, and $\sigma$ is a symmetric total Cauchy stress tensor. We now decompose $\sigma=\widetilde{\tau}-pI_{d}$, where $\widetilde{\tau}$ is symmetric extra stress tensor and $p$ is the pressure of the fluid.  Then, we obtain the following system of the incompressible fluid:
\eqn \label{incompressible fluid eq}
\begin{split}
& u_{t}+u\cdot\nabla u+\nabla p=\dv \widetilde{\tau},\\
& \dv u=0.
\end{split}
\een

To complete the system (\ref{incompressible fluid eq}), we need to add the constitutive law of $\widetilde{\tau}$. In a Newtonian fluid, the constitutive law is given by 
\eqn \label{Newtonian tau}
\widetilde{\tau}=2\nu_{1}D(u), \quad D(u)=\frac{\nabla u+\nabla u^{T}}{2}
\een
with the viscosity coefficient of the fluid $\nu_{1}$. Then, (\ref{incompressible fluid eq}) becomes the  incompressible Navier-Stokes equations. On the other hand, in a viscoelastic fluid, $\widetilde{\tau}$ depends not only on the current motion of the fluid, but also on the history of the motion \cite[Chapter 2]{Renardy}. One of the linear models that relate the history of the motion of the fluid to $\widetilde{\tau}$ is the Jeffreys model:
\eqn \label{Jeffreys model}
\widetilde{\tau}+\lambda_{1}\frac{\partial \widetilde{\tau}}{\partial t}=2\eta_{0}\left(D(u)+\lambda_{2}\frac{\partial D(u)}{\partial t}\right),
\een
where $\lambda_{1}$ and $\lambda_{2}$ are the relaxation time and retardation time, respectively, satisfying $0\leq\lambda_{2}<\lambda_{1}$. However, (\ref{Jeffreys model}) cannot explain several rheological phenomena of viscoelastic fluids. Moreover,  (\ref{Jeffreys model}) does not satisfy the principle of material frame-indifference, or the objectivity. We will present the nonlinear and objective constitutive models after introducing the meaning of the objectivity.

In classical mechanics, an observer, or a frame, is defined as a rigid body with a clock. A change of frame is the transformation from the pair $(x,t)$ in one frame to the pair $(x^{*},t^{*})$ in a different frame, where $x$ and $x^{*}$ are the position vectors and we may take $t=t^{*}$. Since the two frames are rigid bodies, the change of frame is given by 
\eqn \label{change of frame}
x^{*}=c(t)+Q(t)(x-x_{b}),
\een
where $c(t)$ is the relative displacement of the base point $x_{b}$, $Q(t)$ is a time-dependent orthogonal tensor. In other words, the frames move relative to each other, and can also rotate or reflect. A quantity is called \emph{objective} if it is frame-independent.

\begin{definition} \label{objective quantities} \cite{Lai} \upshape
Under the change of frame given by (\ref{change of frame}),
\begin{enumerate}
\item a vector $v$ is objective if $v^{*}=Q(t)v$,
\item a tensor $\tau$ is objective if  $\tau^{*}=Q(t)\tau Q^{T}(t)$.
\end{enumerate}
\end{definition}

According to the principle of material indifference of continuum mechanics, the Cauchy stress tensor $\sigma=\widetilde{\tau}-pI_{d}$ is objective, and so is $\widetilde{\tau}$. We next see whether $\nabla u$ is objective or not. Since 
\[
u^{*}=\frac{dx^{*}}{dt}=\dot{c}+\dot{Q}(x-x_{b})+Q\frac{dx}{dt}=\dot{c}+\dot{Q}(x-x_{b})+Qu, \quad \frac{dx^{*}}{dx}=Q(t), 
\]
we derive $\nabla^{*}u^{*}=Q\nabla u Q^{T}+\dot{Q}Q^{T}$ which tells that $\nabla u$ is also not an objective tensor. We now turn to the symmetric and skew-symmetric parts of $\nabla u$:
\[
D(u)=\frac{\nabla u+\nabla u^{T}}{2},\quad\Omega(u)=\frac{\nabla u-\nabla u^{T}}{2}.
\] 
Using $\nabla^{*}u^{*}=Q\nabla u Q^{T}+\dot{Q}Q^{T}$ and orthogonality of $Q(t)$, we can find that
\eqn \label{starred deformation and vorticity tensors}
D^{*}=QDQ^{T},\quad \Omega^{*}=Q\Omega Q^{T}+\dot{Q}Q^{T}.
\een
Thus, $D(u)$ is objective, but $\Omega(u)$ is not objective.

In order to construct the proper constitutive equations, we need objective derivatives for an objective tensor $\tau$. At first, let us check the material derivative 
\[
\frac{D}{Dt}=\frac{\partial}{\partial t}+u^{j}\frac{\partial}{\partial x_{j}}.
\]
Noting that $\displaystyle \frac{D^{*}}{Dt}=\frac{D}{Dt}$ when $t=t^{\ast}$ and $\tau^{*}=Q(t)\tau Q^{T}(t)$, we see that 
\eqn \label{starred material derivative}
\frac{D}{Dt}\tau^{*}=\left(\frac{DQ}{Dt}\right)\tau Q^{T}+Q\left(\frac{D\tau}{Dt}\right)Q^{T}+Q\tau\left(\frac{DQ^{T}}{Dt}\right)=\dot{Q}\tau Q^{T}+Q\frac{D\tau}{Dt}Q^{T}+Q\tau\dot{Q}^{T}
\een
and so the material derivative  is not objective. We now define the Oldroyd derivative as
\[
\frac{\mathcal{D}\tau}{\mathcal{D}t}=\frac{D\tau}{Dt}+\tau\Omega-\Omega\tau
\]
which is objective by (\ref{starred deformation and vorticity tensors}) and (\ref{starred material derivative}). Similarly, the following derivatives are  objective:
\[
\begin{split}
&\frac{\mathcal{D}_{b}\tau}{\mathcal{D}t}=\frac{D\tau}{Dt}+\mathcal{Q}_{b}(\tau,\nabla u),\\& \mathcal{Q}_{b}(\tau,\nabla u)=\tau\Omega-\Omega\tau-b(D(u)\tau+\tau D(u)), \quad b\in [-1,1].
\end{split}
\]

\subsection{Constitutive equations}\label{sec:1.2}
We here provide some nonlinear and objective constitutive models. The most well-known model is  described by Oldroyd \cite{Oldroyd 1} with
\eqn \label{Oldroyd model}
\widetilde{\tau}+\lambda_{1}\frac{\mathcal{D}_{b} \widetilde{\tau}}{\mathcal{D} t}=2\eta_{0}\left(D(u)+\lambda_{2}\frac{\mathcal{D}_{b} D(u)}{\mathcal{D} t}\right)
\een
which can be obtained by replacing the partial derivative in time in (\ref{Jeffreys model}) with the above objective derivative. (\ref{Oldroyd model}) is called Oldroyd-A model when $b=-1$ and Oldroyd-B model when $b=1$. Oldroyd-B model is preferred than Oldroyd-A model since Oldroyd-B model exhibits the rod climbing effect as rheological phenomena, but the other does the effect in reverse. 

We now reorganize (\ref{Oldroyd model}) into manageable forms by decomposing the extra stress tensor as 
\eqn \label{tensor decompostion}
\widetilde{\tau}=\tau_{N}+\tau, \quad \tau_{N}=2\eta_{0}\frac{\lambda_{2}}{\lambda_{1}}D(u),
\een
where $\tau_{N}$ and $\tau$ correspond to the Newtonian part and the elastic part, respectively. Then, (\ref{Oldroyd model}) is reduced to 
\eqn \label{Oldroyd model2}
\alpha\tau+\frac{\mathcal{D}_{b}\tau}{\mathcal{D}t}=\mu D(u)
\een
allowing to choose proper non-negative constants $\alpha$ and $\mu$. We also introduce the Johnson-Segalman model, the Phan-Thien-Tanner (PTT) model, the Giesekus model, and the generalized Oldroyd model. To do so, we first give
\[
\begin{split}
& \overset{\triangledown}{\tau}=\tau_{t}+u\cdot \nabla \tau-(\nabla u)^{T} \cdot \tau-\tau \cdot \nabla u\\
& \overset{\vartriangle}{\tau}=\tau_{t}+u\cdot \nabla \tau+\tau \cdot (\nabla u)^{T} +\nabla u \cdot \tau.
\end{split}
\]

\begin{enumerate}[]
\item \textbullet \ \ Johnson-Segalman model \cite{Johnson}:
\[
\alpha\tau+\left(\frac{1+a}{2}\overset{\triangledown}{\tau}+\frac{1-a}{2} \overset{\vartriangle}{\tau}\right)=\mu D(u), \quad  a\in [-1,1]
\]
\item \textbullet \ \ Phan-Thien-Tanner (PTT) model \cite{PTT}:
\[
f(\tau)\tau + \frac{\mathcal{D}_{b}\tau}{\mathcal{D}t}=\mu D(u), \quad f(\tau)=\left\{ \begin{array}{lll}
\alpha e^{\text{tr}\tau} & \ \text{(exponential)}\\
\alpha+ \beta\text{tr}\tau & \ \text{(linear)},
\end{array} \right.
\]
where $b=1$ and $\text{tr}\tau$ is the trace of $\tau$. In this paper, we take the linear $f(\tau)$ in (PTT) model. 
\item \textbullet \ \ Giesekus model \cite{Giesekus}:
\[
\alpha\tau+\overset{\triangledown}{\tau}+\beta \tau^{2}=\mu D(u).
\]
\item \textbullet \ \ Generalized Oldroyd model \cite{Oldroyd 2}:
\[
\begin{split}
& \alpha\tau+\frac{\mathcal{D}_{b}\tau}{\mathcal{D}t}+\widetilde{\mathcal{Q}}(\tau,\nabla u)=\mu D(u),\\
& \widetilde{\mathcal{Q}}(\tau,\nabla u)=c_{1}\text{tr}\tau D(u)+c_{2}(\tau:D(u))I+c_{3}(D(u))^{2}+c_{4}(D(u):D(u))I.
\end{split}
\]
\end{enumerate}

Combined with (\ref{incompressible fluid eq}) and (\ref{tensor decompostion}), the models described above can be written as 
\begin{subequations} \label{our model}
\begin{align} 
& u_{t}+u\cdot\nabla u -\nu_{1}\Delta u+\nabla p=\dv \tau, \label{our model a}\\
& \tau_{t}+u\cdot\nabla \tau+\alpha\tau-\nu_{2}\Delta\tau+\mathcal{Q}(\tau, \nabla u)= \mu D(u), \label{our model b}\\
& \dv u=0,
\end{align}
\end{subequations}
where $\nu_{1},\nu_{2}, \alpha$ are non-negative constants, and we set $\mu=1$ from now on for simplicity. $\mathcal{Q}(\tau, \nabla u)$, which is a quadratic function of $\tau$ and $\nabla u$, represents the nonlinear term of each model described above.


\subsection{Our results}

To determine function spaces for initial data and solutions, we may use a scaling-invariant property.  In fact, (\ref{our model}) is not scaling-invariant, but if we neglect $\dv \tau$ and $\alpha\tau$, and restrict $\mathcal{Q}(\tau,\nabla u)$ to $\tau\nabla u$, (\ref{our model}) is invariant under the transformation 
\[
\left(u(t,x),\tau(t,x)\right)\longmapsto  \left(\lambda u(\lambda^{2}t,\lambda x),\tau(\lambda^{2}t,\lambda x)\right)
\]
when $\nu_{1}>0$ and $\nu_{2}>0$. Accordingly we can take the scaling invariant property for initial data
\eqn \label{critical transformation data}
\left(u_{0}(x),\tau_{0}(x)\right)\longmapsto \left(\lambda u_{0}(\lambda x),\tau_{0}(\lambda x)\right).
\een
So, we may take $(u_{0},\tau_{0})$ belonging to $\dot{B}^{\frac{d}{p}-1}_{p,1}\times \dot{B}^{\frac{d}{p}}_{p,1}$ which is invariant under (\ref{critical transformation data}). But, in this paper, we take the model (\ref{our model}) with at least one of $(\nu_{1}, \nu_{2},\alpha)$ being  zero:
\begin{enumerate}[]
\item $\vartriangleright$  Case I: $(\nu_{1}, \nu_{2}, \alpha)=(+,0,+)$;
\item $\vartriangleright$ Case II: $(\nu_{1}, \nu_{2}, \alpha)=(+,0,0)$;
\item $\vartriangleright$ Case III: $(\nu_{1}, \nu_{2}, \alpha)=(0,+,+)$;
\item $\vartriangleright$ Case IV:  $(\nu_{1}, \nu_{2}, \alpha)=(0,+,0)$;
\item $\vartriangleright$ Case V: $(\nu_{1}, \nu_{2}, \alpha)=(0,0,+)$,
\end{enumerate}
where  $+$ means a positive constant. Therefore, we need a space stronger  than $\dot{B}^{\frac{d}{2}-1}_{2,1}\times\dot{B}^{\frac{d}{2}}_{2,1}$. Motivated by the compressible Navier-Stokes equations \cite{Danchin 1, Danchin 2}, we take initial data in hybrid Besov space  whose norm is given by 
\[
\left\|u\right\|_{\dot{B}^{s,t}}=\sum_{j\leq j_{0}}2^{js}\left\|\Delta_{j}u\right\|_{L^{2}}+\sum_{j\geq j_{0}+1}2^{jt}\left\|\Delta_{j}u\right\|_{L^{2}}.
\]

The goal of this paper is twofold:
\begin{enumerate}[]
\item (1) The existence of global-in-time solutions of (\ref{our model}) in the framework of hybrid Besov spaces: there exists a constant $\epsilon>0$ such that if $(u_{0},\tau_{0})\in \mathcal{E}_{0}$ with $\left\|(u_{0},\tau_{0})\right\|_{\mathcal{E}_{0}}\leq \epsilon$, then there is a unique solution of (\ref{our model}) such that $\left\|(u,\tau)\right\|_{\mathcal{E}_{T}}\lesssim \left\|(u_{0},\tau_{0})\right\|_{\mathcal{E}_{0}}$ for all $T>0$. The spaces $\mathcal{E}_{0}$ and $\mathcal{E}_{T}$ are given by the following table where we classify the possible form of $\mathcal{Q}(\tau,\nabla u)$ in each case (Since we treat $\mathcal{Q}(\tau,\nabla u)$ as a perturbation of the linear part of (\ref{our model}), the exact form of $\mathcal{Q}(\tau,\nabla u)$ is not expressed, but is written as in Table.) This result covers all known results in Besov spaces and energy spaces. Moreover, we provide new existence results with $\mathcal{Q}(\tau,\nabla u)$ having additional terms. 
\item (2) The decay rates of $(u,\tau)$ except for Case $\text{II}$: 
\eqn \label{decay rate}
\left\|u(t)\right\|_{\dot{B}^{l(u)+s_{0},h(u)}}+\left\|\tau(t)\right\|_{\dot{B}^{l(\tau)+s_{0},h(\tau)}}\leq C(\left\|(u_{0},\tau_{0})\right\|_{\mathcal{E}_{0}}, s_{0})(1+t)^{-\frac{s_{0}}{2}}, \quad\forall s_{0}\in(0,2].
\een
Here, $l(u)$ and $h(u)$ denote low  and high frequency regularity of $u$, respectively, that guarantee the global existence of unique solutions in hybrid Besov spaces. Similarly, $l(\tau)$ and $h(\tau)$ denote low and high frequency regularity of $\tau$. Our approach to derive (\ref{decay rate}) is straightforward, but this type of decay results are totally new. Furthermore, as an application of (\ref{decay rate}), we can obtain further decay rates by imposing additional initial data condition ($L^{2}$ or $L^{1}$ type), which cover some recent works related to decay rates for Oldroyd type models \cite{Hieber, Huang, Wang}. We emphasize that (\ref{decay rate}) indicates that decay rates depend only on regularities in low frequency parts, and even in deriving further decay rates, low frequency parts play a crucial role.
\end{enumerate}

\begin{table}[h!]
\centering
\resizebox{\columnwidth}{!}{%
\begin{tabular}{@{}|c|ccc|c|c|c|c@{}}
\toprule
Case & $\nu_{1}$ & $\nu_{2}$ & $\alpha$ & $\mathcal{E}_{0}$ & $\mathcal{E}_{T}$ & $\mathcal{Q}(\tau, \nabla u)\simeq$ \\
\midrule
I& + & 0 & + & $\dot{B}^{\frac{d}{2}-1,\frac{d}{2}+1}\times\dot{B}^{\frac{d}{2},\frac{d}{2}+1}$ & $L^{\infty}_{T}\dot{B}^{\frac{d}{2}-1,\frac{d}{2}+1}\cap L^{1}_{T}\dot{B}^{\frac{d}{2}+1,\frac{d}{2}+2}\times L^{\infty}_{T}\dot{B}^{\frac{d}{2},\frac{d}{2}+1}\cap L^{1}_{T}\dot{B}^{\frac{d}{2},\frac{d}{2}+1}$ & $\tau^{2}+\tau\nabla u+ (\nabla u)^{2}$  \\
\hline
II& + & 0 & 0 & $\dot{B}^{\frac{d}{2}-1,\frac{d}{2}+1}\times\dot{B}^{\frac{d}{2}-1,\frac{d}{2}+1}$ & $L^{\infty}_{T}\dot{B}^{\frac{d}{2}-1,\frac{d}{2}+1}\cap L^{1}_{T}\dot{B}^{\frac{d}{2}+1,\frac{d}{2}+2}\times L^{\infty}_{T}\dot{B}^{\frac{d}{2}-1,\frac{d}{2}+1}$  &  $\tau\nabla u+ (\nabla u)^{2}$  \\
 \hline 
III & 0 & + & + & $\dot{B}^{\frac{d}{2}-1,\frac{d}{2}+1}\times\dot{B}^{\frac{d}{2}}_{2,1}$ & $L^{\infty}_{T}\dot{B}^{\frac{d}{2}-1,\frac{d}{2}+1}\cap L^{1}_{T}\dot{B}^{\frac{d}{2}+1}_{2,1}\times L^{\infty}_{T}\dot{B}^{\frac{d}{2}}_{2,1}\cap L^{1}_{T}\dot{B}^{\frac{d}{2},\frac{d}{2}+2}$  & $\tau^{2}+\tau\nabla u+ (\nabla u)^{2}$ \\
 \hline
IV & 0 & + & 0 & $\dot{B}^{\frac{d}{2}-1,\frac{d}{2}+1}\times\dot{B}^{\frac{d}{2}-1,\frac{d}{2}}$ & $L^{\infty}_{T}\dot{B}^{\frac{d}{2}-1,\frac{d}{2}+1}\cap L^{1}_{T}\dot{B}^{\frac{d}{2}+1}_{2,1}\times L^{\infty}_{T}\dot{B}^{\frac{d}{2}-1,\frac{d}{2}}\cap L^{1}_{T}\dot{B}^{\frac{d}{2}+1,\frac{d}{2}+2}$  & $\tau\nabla u+ (\nabla u)^{2}$ \\
 \hline
V & 0 & 0 & + & $\dot{B}^{\frac{d}{2}-1,\frac{d}{2}+1}\times\dot{B}^{\frac{d}{2},\frac{d}{2}+1}$ & $L^{\infty}_{T}\dot{B}^{\frac{d}{2}-1,\frac{d}{2}+1}\cap L^{1}_{T}\dot{B}^{\frac{d}{2}+1}_{2,1}\times L^{\infty}_{T}\dot{B}^{\frac{d}{2},\frac{d}{2}+1}\cap L^{1}_{T}\dot{B}^{\frac{d}{2},\frac{d}{2}+1}$  & $\tau^{2}$ \\ \bottomrule
\end{tabular}
}
\label{table}
\end{table}

There has been a lot of work on (\ref{our model}) and related models: see \cite{Bejaoui, Chemin, Chen-Hao, Chen-Liu, Chen, Chen 1, Chen 2, Chupin, Constantin 1, Constantin 2, Constantin 3, Constantin 4, Constantin 5, Constantin 6, De Anna, Elgindi 1, Elgindi 2, Fang 1, Fang 2, Fernandez, Guillope, Hieber, Hu, Huang, Lei 1, Lei, Lions, Masmoudi, Renardy 1, Renardy 2, Renardy, Wan, Wang, Wu, Ye 1, Ye 2, Zhu, Zi, Zi-Fang-Zhang} and references therein. Here we compare some results with ours in terms of the set of parameters $(\nu_{1}, \nu_{2}, \alpha)$, the function spaces and decay rates.
\begin{enumerate}[]
\item {\bf Case I ($(\nu_{1}, \nu_{2},\alpha)=(+,0,+)$):} This is the case which have been studied widely, especially with $\mathcal{Q}(\tau, \nabla u)=\mathcal{Q}_{b}(\tau,\nabla u)$: see \cite{Chemin, Chen 2, Fang 1, Fang 2, Fernandez, Guillope, Hieber, Hu, Huang, Lei, Lions, Wan, Zi-Fang-Zhang}. In particular, global well-posedness in critical Besov spaces has been established in \cite{Chemin} and improved in \cite{Zi-Fang-Zhang}. We note that $\mathcal{Q}_{b}(\tau,\nabla u) \simeq \tau\nabla u$. Compared to these results, we find a proper hybrid Besov space to show the global well-posedness with $\mathcal{Q}(\tau,\nabla u)\simeq \tau^{2}+\tau \nabla u+(\nabla u)^{2}$ (Theorem \ref{case1 theorem}). Furthermore, we obtain (\ref{decay rate}) type of decay rates for both scenarios (Theorem \ref{case1 decay} and Corollary \ref{case1 decay2}). Also, to compare previous works related to decay rates with ours, we assume additional initial data condition and obtain further decay rates, which improve the results in \cite{Hieber, Huang}: see Corollary \ref{case1 further decay} and Remark \ref{case1 further decay remark} for more details. 
\item {\bf Case II ($(\nu_{1}, \nu_{2},\alpha)=(+,0,0)$):} This case with $\mathcal{Q}(\tau, \nabla u)=\mathcal{Q}_{b}(\tau,\nabla u)$ has been studied in \cite{Chen-Hao, Zhu}, and the global well-posedness in hybrid Besov spaces is dealt with in \cite{Chen-Hao}. As Case I, we find a proper hybrid Besov space that guarantee the global well-posedness with $\mathcal{Q}(\tau,\nabla u)\simeq \tau \nabla u+(\nabla u)^{2}$ (Theorem \ref{case2 theorem}).
\item {\bf Case III ($(\nu_{1}, \nu_{2},\alpha)=(0,+,+)$):} The global well-posedness with small initial data in $H^{s}$ for $s>\frac{d}{2}+1$, is considered in \cite{Elgindi 1, Elgindi 2}. By contrast, we prove such results in hybrid Besov spaces, which is  larger than $H^{s}$ (Theorem \ref{case3 theorem}). The decay rates of solutions with initial data in $L^{1}\cap H^{s}$ are shown in \cite{Huang}, and improved by ours (Theorem \ref{case3 decay} and Corollary \ref{case3 further decay}).
\item {\bf Case IV ($(\nu_{1}, \nu_{2},\alpha)=(0,+,0)$):} This case with $\mathcal{Q}(\tau,\nabla u)=\mathcal{Q}_{b}(\tau,\nabla u)$ has been studied in \cite{Chen-Liu, Constantin 6, Wang, Wu}, and the global well-posedness in hybrid Besov spaces is dealt with in \cite{Chen-Liu, Wu}. We can check similar results with $\mathcal{Q}(\tau,\nabla u)\simeq \tau\nabla u+(\nabla u)^{2}$ (Theorem \ref{case4 theorem}). We also derive decay rates of solutions which improve \cite{Wang}: See Theorem \ref{case4 decay}, Corollary \ref{case4 further decay} and Remark \ref{case4 further decay remark}.
\item {\bf Case V ($(\nu_{1}, \nu_{2},\alpha)=(0,0,+)$):}
This case is dealt with in \cite{Zi}  in 3D with  $(u_{0}, \tau_{0}) \in \dot{B}^{\frac{1}{2},\frac{5}{2}}\times\dot{B}^{\frac{1}{2},\frac{5}{2}}$ and $\mathcal{Q}(\tau, \nabla u)=0$. In this paper, we improve this result with $(u_{0}, \tau_{0}) \in \dot{B}^{\frac{d}{2}-1,\frac{d}{2}+1}\times\dot{B}^{\frac{d}{2},\frac{d}{2}+1}$ with $d=2,3$, and $\mathcal{Q}(\tau, \nabla u)\simeq\tau^{2}$ (Theorem \ref{case5 theorem}). We also derive decay rates of solutions (Theorem \ref{case5 decay} and Corollary \ref{case5 further decay}).
\end{enumerate}

\subsection{Detailed plan of the proofs of our results}
To prove our results, we may use a fixed point argument. But, in principle, the calculations used to derive a priori estimates are easily applied to the fixed point argument. So, we only provide a priori estimates for the existence part and  then show the uniqueness of solutions and derive decay rates.  The proofs of our results look similar in each case, but we will  provide as much of the proof process as possible, rather than referring to some bounds obtained before,  because the function spaces are slightly different in each case.  

\begin{enumerate}[]
\item \textbullet  \ In each section, we begin with a  linearized model of (\ref{our model}) with a divergence-free vector field $v$,
\begin{subequations} \label{Linearized Equations}
\begin{align} 
& u_{t}+\mathbb{P}(v\cdot\nabla u) -\nu_{1}\Delta u-\mathbb{P}\dv \tau=F, \label{Linearized Equations a}\\
& \tau_{t}+v\cdot\nabla \tau+\alpha\tau-\nu_{2}\Delta\tau-D(u)=G, \label{Linearized Equations b} \\
& (\mathbb{P}\dv \tau)_{t} +\mathbb{P}\dv (v\cdot \nabla \tau)+\alpha\mathbb{P} \dv \tau-\nu_{2}\Delta\mathbb{P} \dv \tau- \frac{1}{2}\Delta u=\mathbb{P}\dv G, \label{Linearized Equations c}
\end{align}
\end{subequations}
where $\mathbb{P}$ is the Leray projection operator. (\ref{Linearized Equations c}) is derived by taking $\mathbb{P}\dv$ to  (\ref{Linearized Equations b}) because  (\ref{Linearized Equations a}) contains $\mathbb{P}\dv \tau$. We first bound solutions of (\ref{Linearized Equations}) in Lemma \ref{case1 lemma}, Lemma \ref{case2 lemma}, Lemma \ref{case3 lemma}, Lemma \ref{case4 lemma}, and Lemma \ref{case5 lemma}. These lemmas, with $v=u$, $F=0$, and $G=-\mathcal{Q}(\tau,\nabla u)$, are used to derive the bound 
$\left\|(u,\tau)\right\|_{\mathcal{E}_{T}}\lesssim \left\|(u_{0},\tau_{0})\right\|_{\mathcal{E}_{0}}$ for all $T>0$ when $\left\|(u_{0},\tau_{0})\right\|_{\mathcal{E}_{0}}$ is sufficiently small. 

We can also derive decay rates of solutions (\ref{decay rate}) by using the bounds obtained in the proofs of Lemma \ref{case1 lemma}, Lemma \ref{case3 lemma}, Lemma \ref{case4 lemma}, and Lemma \ref{case5 lemma}. However, we are not able to derive the decay rates of solutions in Case II: see Remark \ref{case2 decay remark} for more details. Furthermore, we can derive the faster decay rates by intersecting $\mathcal{E}_{0}$ with $L^{2}$ or $L^{1}$ type. Simply speaking, decay rates with $L^{2}$ intersection come from applying interpolation with $L^{2}$ when deriving (\ref{decay rate}). In three dimension, owing to obtain decay rates with $L^{1}$ intersection, we use a variant of Schonbek's Fourier splitting method \cite{Schonbek} in our setting. To this end, we consider $\dot{B}^{-\frac{3}{2}}_{2,\infty}$ since it is larger than $L^{1}$ and more natural in this framework: see Corollary \ref{case1 further decay}, Corollary \ref{case3 further decay}, Corollary \ref{case4 further decay} and Corollary \ref{case5 further decay} for more details.

\item \textbullet  \ To prove the uniqueness of solutions, we assume that there are two solutions $(u_{1}, \tau_{1})$ and $(u_{2}, \tau_{2})$ and deal with the following equations for $(u,\tau)=(u_{1}-u_{2},\tau_{1}-\tau_{2})$:
\[
\begin{split}
& u_{t}+\mathbb{P}(u_{2}\cdot\nabla u)-\nu_{1}\Delta u-\mathbb{P}\dv \tau=-\mathbb{P}(u\cdot\nabla u_{1}),\\
& \tau_{t}+u_{2}\cdot\nabla\tau+\alpha\tau-\nu_{2}\Delta \tau-D(u)=-u\cdot\nabla\tau_{1}-\left[\mathcal{Q}(\tau_{1},\nabla u_{1})-\mathcal{Q}(\tau_{2},\nabla u_{2})\right], \\
&\mathcal{Q}(\tau_{1},\nabla u_{1})-\mathcal{Q}(\tau_{2},\nabla u_{2}) \simeq a(\tau_{1}+\tau_{2})\tau+b(\tau\nabla u_{1}+\tau_{2}\nabla  u)+c(\nabla u_{1}+\nabla u_{2})\nabla u
\end{split}
\]
for given $\mathcal{Q}(\tau,\nabla u)\simeq a\tau^{2}+b\tau\nabla u+c(\nabla u)^{2}$. Then $(u,\tau)$ satisfy (\ref{Linearized Equations}) with $v=u_{2}$, $F\mapsto \delta F=-\mathbb{P}(u\cdot\nabla u_{1})$, and $G\mapsto \delta G=-u\cdot\nabla\tau_{1}-\left[\mathcal{Q}(\tau_{1},\nabla u_{1})-\mathcal{Q}(\tau_{2},\nabla u_{2})\right]$. We then modify Lemma \ref{case1 lemma}, Lemma \ref{case2 lemma}, Lemma \ref{case3 lemma}, Lemma \ref{case4 lemma}, and Lemma \ref{case5 lemma} to bound $(u, \tau)$ in $\mathcal{F}_{T}$ (which will be selected  in each section), which is slightly larger than $\mathcal{E}_{T}$, as follows
\[
\left\|(u, \tau)\right\|_{\mathcal{F}_{T}}\lesssim \Big(\left\|(u_{1},\tau_{1})\right\|_{\mathcal{E}_{T}}+\left\|(u_{2},\tau_{2})\right\|_{\mathcal{E}_{T}}\Big)\left\|(u, \tau)\right\|_{\mathcal{F}_{T}}.
\]
Since $\left\|(u_{1},\tau_{1})\right\|_{\mathcal{E}_{T}}+\left\|(u_{2},\tau_{2})\right\|_{\mathcal{E}_{T}}$ is sufficiently small, we can confirm the uniqueness of solutions.
\end{enumerate}

After providing some auxiliary materials in Section  \ref{preliminaries}, we will give our results from Section \ref{case1}.

\section{Preliminaries} \label{preliminaries}

\subsection{Notation.}
\begin{enumerate}[]
\item (1) All  constants will be denoted by $C$ and we follow the convention that such constants can vary from expression to expression and even between two occurrences within the same expression.  $A \lesssim B $ means there is a constant $C$ such that $A \leq C B $. $A\simeq B$ means $A\lesssim B$ and $B\lesssim A$. 
\item (2) $L^{p}_{T}X=L^{p}([0,T];X)$.
\item (3) $(f, g)$ denotes the $L^{2}$ inner product of $f$ and $g$. This notation is also used to bound two quantities in the same space.
\end{enumerate}

\subsection{Littlewood-Paley decomposition} \cite[Chapter 2]{Bahouri}
Let $\mathcal{C}=\left\{\xi\in \mathbb{R}^{d}: \frac{3}{4}\leq |\xi|\leq \frac{8}{3}\right\}$ and take a smooth radial function $\phi$ with values in $[0,1]$ supported on  $\mathcal{C}$  and satisfies
\[
\sum^{\infty}_{j=-\infty}\phi\left(2^{-j}\xi\right)=1 \ \  \forall \ \xi \in \mathbb{R}^{d}\setminus\{0\}.
\]
Using this, we define dyadic blocks:
\[
\Delta_{j}f=2^{jd} \int_{\mathbb{R}^{d}} h\left(2^{j}y\right)f(x-y)dy, \quad \widehat{h}(\xi)=\phi(\xi).
\] 
Then, the homogeneous Littlewood-Paley decomposition is given by
\[
f=\sum_{j\in \mathbb{Z}} \Delta_{j}f \ \  \text{in} \ \  \mathcal{S}^{'}_{h},
\]
where $\mathcal{S}^{'}_{h}$ is the space of tempered distributions $u\in \mathcal{S}^{'}$ such that $\displaystyle \lim_{j\rightarrow -\infty}S_{j}u=0$ in $\mathcal{S}^{'}$. 


We define the homogeneous Besov spaces: for $1\leq p,q\leq \infty$ and $s\in\mathbb{R}$,
\[
\dot{B}^{s}_{p,q}=\left\{f\in \mathcal{S}^{'}_{h};\   \|f\|_{\dot{B}^{s}_{p,q}}=\left\|2^{js}\left\|\Delta_{j}f\right\|_{L^{p}}\right\|_{l^{q}(\mathbb{Z})}<\infty \right\}.
\]
We also define the hybrid Besov spaces with $p=2$, $q=1$: for $s,t\in\mathbb{R}$
\[
\left\|f\right\|_{\dot{B}^{s,t}}=\sum_{j\leq j_{0}}2^{js}\left\|\Delta_{j}f\right\|_{L^{2}}+\sum_{j\geq j_{0}+1}2^{jt}\left\|\Delta_{j}f\right\|_{L^{2}}
\]
for some $j_{0}$ and $\dot{B}^{s,t}=\left\{f\in \mathcal{S}^{'}_{h};\   \|f\|_{\dot{B}^{s,t}}<\infty \right\}$. We denote 
\eqn \label{low and high norm}
\begin{split}
&\left\|f\right\|_{\dot{B}^{s}}:=\left\|f\right\|_{\dot{B}^{s,s}}=\left\|f\right\|_{\dot{B}^{s}_{2,1}},\\
&\left\|f\right\|^{l}_{\dot{B}^{s}}:=\sum_{j\leq j_{0}}2^{js}\left\|\Delta_{j}f\right\|_{L^{2}},\quad \left\|f\right\|^{h}_{\dot{B}^{s}}:=\sum_{j\geq j_{0}+1}2^{js}\left\|\Delta_{j}f\right\|_{L^{2}}.
\end{split}
\een
From the definition, we easily check that
\eqn \label{hybrid Besov embedding}
\dot{B}^{s,t}=\dot{B}^{s}_{2,1}\cap\dot{B}^{t}_{2,1}\ \ \text{if}\ \ s\leq t,\quad \dot{B}^{s_{1},t_{1}}\hookrightarrow \dot{B}^{s_{2},t_{2}}\ \ \text{if}\ \ s_{1}\leq s_{2},\ \ t_{1}\geq t_{2}. 
\een 
We note that the definition of $\dot{B}^{s,t}$ does not depend on the choice of $j_{0}$. This is the reason why we can choose $j_{0}$ not at the beginning but during the proofs of our results as in \cite{Danchin 1, Danchin 2}. We note that homogeneous Besov spaces satisfy the following interpolation relationship: for $\theta\in[0,1]$, 
\eqn \label{interpolation}
\left\|f\right\|_{\dot{B}^{\theta s_{1}+(1-\theta)s_{2}}_{p,q}}\leq \left\|f\right\|^{\theta}_{\dot{B}^{s_{1}}_{p,q}}\left\|f\right\|^{1-\theta}_{\dot{B}^{s_{2}}_{p,q}},
\een
and for $0<s<s'$,
\eqn \label{interpolation 2}
\left\|f\right\|_{\dot{B}^{s}_{p,1}}\leq C\left\|f\right\|^{1-\frac{s}{s'}}_{\dot{B}^{0}_{p,\infty}}\left\|f\right\|^{\frac{s}{s'}}_{\dot{B}^{s'}_{p,\infty}}\leq C\left\|f\right\|^{1-\frac{s}{s'}}_{L^{p}}\left\|f\right\|^{\frac{s}{s'}}_{\dot{B}^{s'}_{p,1}}.
\een
We also note that (\ref{interpolation}) holds for $\left\|\cdot\right\|^{l}_{\dot{B}^{s}_{p,q}}$ and $\left\|\cdot\right\|^{h}_{\dot{B}^{s}_{p,q}}$.

\subsection{Product estimates}
For two distributions $f$ and $g$, we introduce Bony's decomposition:
\[
fg=\mathcal{T}_{f}g+\mathcal{T}_{g}f+\mathcal{R}(f,g), \quad \mathcal{T}_{f}g=\sum_{j\in\mathbb{Z}}S_{j-1}f\Delta_{j}g,\quad \mathcal{R}(f,g)=\sum_{|j-j'|\leq 1}\Delta_{j}f\Delta_{j'}g.
\]

\begin{proposition} \cite[Appendix]{Danchin 2}\label{paraproduct prop} \upshape
For all $s_{1},s_{2},t_{1},t_{2}$ such that $s_{1},s_{2}\leq \frac{d}{2}$, we have
\[
\left\|\mathcal{T}_{f}g\right\|_{\dot{B}^{s_{1}+t_{1}-\frac{d}{2},s_{2}+t_{2}-\frac{d}{2}}}\lesssim \left\|f\right\|_{\dot{B}^{s_{1},s_{2}}}\left\|g\right\|_{\dot{B}^{t_{1},t_{2}}}.
\]
If $\min{(s_{1}+t_{1},s_{2}+t_{2})}>0$, then
\[
\left\|\mathcal{R}(f,g)\right\|_{\dot{B}^{s_{1}+t_{1}-\frac{d}{2},s_{2}+t_{2}-\frac{d}{2}}}\lesssim \left\|f\right\|_{\dot{B}^{s_{1},s_{2}}}\left\|g\right\|_{\dot{B}^{t_{1},t_{2}}}.
\]
If $f\in L^{\infty}$,
\[
\left\|\mathcal{T}_{f}g\right\|_{\dot{B}^{t_{1},t_{2}}}\lesssim \left\|f\right\|_{L^{\infty}}\left\|g\right\|_{\dot{B}^{t_{1},t_{2}}},
\]
and if $\min{(t_{1},t_{2})}>0$, then
\[
\left\|\mathcal{R}(f,g)\right\|_{\dot{B}^{t_{1},t_{2}}}\lesssim \left\|f\right\|_{L^{\infty}}\left\|g\right\|_{\dot{B}^{t_{1},t_{2}}}.
\]
\end{proposition}

By Proposition \ref{paraproduct prop} and $\dot{B}^{\frac{d}{2}}_{2,1}\hookrightarrow L^{\infty}$, we deduce the following inequalities. For $-\frac{d}{2}<s_{1},s_{2}\leq\frac{d}{2}$ and $t_{1},t_{2}>0$,
\eqn \label{product estimate}
\begin{split}
&\left\|fg\right\|_{\dot{B}^{s_{1},s_{2}}}\lesssim \left\|f\right\|_{\dot{B}^{s_{1},s_{2}}}\left\|g\right\|_{\dot{B}^{\frac{d}{2}}},\\
& \left\|fg\right\|_{\dot{B}^{t_{1},t_{2}}}\lesssim \left\|f\right\|_{\dot{B}^{t_{1},t_{2}}}\left\|g\right\|_{\dot{B}^{\frac{d}{2}}}+\left\|f\right\|_{\dot{B}^{\frac{d}{2}}}\left\|g\right\|_{\dot{B}^{t_{1},t_{2}}}.
\end{split}
\een

\vspace{1ex}

We next deal with the convection terms. We note that $\dv v=0$ enables to extend the range of indices in the following proposition.

\begin{proposition}\cite[Appendix]{Danchin 2} \label{convection prop} \upshape
Let $v$ be a divergence-free vector field. Let $-\frac{d}{2}-1<s_{1},s_{2},t_{1},t_{2}\leq \frac{d}{2}+1$ and $c_{j}\in l^{1}(\mathbb{Z})$. Then, the following bounds hold: for $k\in\mathbb{Z}$
\eqn \label{convection term estimate}
|(\Delta_{j}\Lambda^{k}(v\cdot\nabla u),\Delta_{j}\Lambda^{k} u)| \lesssim c_{j}2^{-j(\varphi^{s_{1},s_{2}}(j)-k)}\left\|v\right\|_{\dot{B}^{\frac{d}{2}+1}}\left\|u\right\|_{\dot{B}^{s_{1},s_{2}}}\left\|\Delta_{j}\Lambda^{k} u\right\|_{L^{2}}, 
\een
and
\eqn \label{convection term estimate 2}
\begin{split}
& |(\Delta_{j}(v\cdot\nabla u),\Delta_{j}w)+(\Delta_{j}(v\cdot\nabla w),\Delta_{j}u)| \\
&\lesssim c_{j}\left\|v\right\|_{\dot{B}^{\frac{d}{2}+1}}\left(2^{-j\varphi^{s_{1},s_{2}}(j)}\left\|w\right\|_{\dot{B}^{s_{1},s_{2}}}\left\|\Delta_{j}u\right\|_{L^{2}}+2^{-j\varphi^{t_{1},t_{2}}(j)}\left\|u\right\|_{\dot{B}^{t_{1},t_{2}}}\left\|\Delta_{j}w\right\|_{L^{2}}\right),
\end{split}
\een
where the function $\varphi^{s_{1},s_{2}}(j)$ is defined by 
\[
\varphi^{s_{1},s_{2}}(j)=
	\begin{cases}
		s_{1},\quad \text{if}\ \ j\leq j_{0}, \\
		s_{2},\quad \text{if}\ \ j\geq j_{0}+1.
	\end{cases}
\]
\end{proposition}

We also recall the following identity introduced in \cite{Zhu}.

\begin{proposition} \label{decompose prop} \upshape
For any smooth tensor $\tau$ and a vector $u$, 
\[
\begin{split}
&\mathbb{P}\dv(u\cdot\nabla\tau)=\mathbb{P}(u\cdot\nabla\mathbb{P}\dv\tau)+\mathbb{B}(\nabla u,\nabla\tau),\\
&\mathbb{B}(\nabla u, \nabla\tau):=\mathbb{P}(\nabla u\cdot\nabla\tau)-\mathbb{P}(\nabla u\cdot\nabla\Delta^{-1}\nabla\cdot\dv \tau).
\end{split}
\]
\end{proposition}
As in \cite{Chen-Hao}, for a divergence-free vector field $u$ 
\eqn \label{decompose divergence free}
\begin{split}
&(\nabla u\cdot\nabla\tau)^{i}=\sum_{j,k}\partial_{j}u^{k}\partial_{k}\tau^{i,j}=\sum_{j,k}\partial_{k}(\partial_{j}u^{k}\tau^{i,j}), \\
&(\nabla u\cdot\nabla\Delta^{-1}\nabla\cdot\dv \tau)^{i}=\sum_{k}\partial_{i}u^{k}\partial_{k}\Delta^{-1}\nabla\cdot\dv\tau=\sum_{k}\partial_{k}(\partial_{i}u^{k}\Delta^{-1}\nabla\cdot\dv\tau).
\end{split}
\een
We finally fix one constant $C_{0}$ which is defined as
\eqn \label{C0}
\left\|\mathbb{P}\dv\tau\right\|_{L^{2}}\leq C_{0}\left\|\Lambda\tau\right\|_{L^{2}}
\een
although we can see that $C_{0}=d$. Here, $\Lambda$ is the fractional Laplacian: $\Lambda=\sqrt{-\Delta}$.

\section{\bf Case I: $(\nu_{1}, \nu_{2},\alpha)=(+,0,+)$} \label{case1} 
In this section, we consider (\ref{our model}) with $(\nu_{1}, \nu_{2},\alpha)=(+,0,+)$:
\begin{subequations} \label{our model Case I}
\begin{align} 
& u_{t}+u\cdot\nabla u -\nu_{1}\Delta u+\nabla p=\dv \tau, \label{our model Case I a}\\
& \tau_{t}+u\cdot\nabla \tau +\alpha\tau+\mathcal{Q}(\tau, \nabla u)=D(u), \label{our model Case I b}\\
& \dv u=0
\end{align}
\end{subequations}
with $\mathcal{Q}(\tau,\nabla u)\simeq\tau^{2}+\tau\nabla u+(\nabla u)^{2}$. We first define the spaces $\mathcal{E}_{0}$ and $\mathcal{E}_{T}$ for initial data and solutions:
\[
\begin{split}
&\|(u_{0}, \tau_{0})\|_{\mathcal{E}_{0}}=\left\|u_{0}\right\|_{\dot{B}^{\frac{d}{2}-1,\frac{d}{2}+1}}+\left\|\tau_{0}\right\|_{\dot{B}^{\frac{d}{2},\frac{d}{2}+1}}, \quad \|(u, \tau)\|_{\mathcal{E}_{T}}=\|(u, \tau)\|_{\mathcal{L}_{T}}+\|(u, \tau)\|_{\mathcal{H}_{T}},\\
& \|(u, \tau)\|_{\mathcal{L}_{T}}=\left\|u\right\|_{L^{\infty}_{T}\dot{B}^{\frac{d}{2}-1,\frac{d}{2}+1}}+\left\|\tau\right\|_{L^{\infty}_{T}\dot{B}^{\frac{d}{2},\frac{d}{2}+1}},\quad \|(u, \tau)\|_{\mathcal{H}_{T}}=\left\|u\right\|_{L^{1}_{T}\dot{B}^{\frac{d}{2}+1,\frac{d}{2}+2}}+\left\|\tau\right\|_{L^{1}_{T}\dot{B}^{\frac{d}{2},\frac{d}{2}+1}}.
\end{split}
\]

\begin{theorem} \label{case1 theorem}\upshape
There exists $\epsilon>0$ such that if $(u_{0},\tau_{0})\in \mathcal{E}_{0}$ with $\left\|(u_{0},\tau_{0})\right\|_{\mathcal{E}_{0}}\leq \epsilon$, then there exists a unique solution $(u,\tau)$ of (\ref{our model Case I}) in $\mathcal{E}_{T}$  such that $\left\|(u,\tau)\right\|_{\mathcal{E}_{T}}\lesssim \left\|(u_{0},\tau_{0})\right\|_{\mathcal{E}_{0}}$ for all $T>0$.
\end{theorem}

\begin{theorem} \label{case1 decay}\upshape
The solution of Theorem \ref{case1 theorem} has the following decay rates:
\[
\left\|u(t)\right\|_{\dot{B}^{\frac{d}{2}-1+s_{0},\frac{d}{2}+1}}+\left\|\tau(t)\right\|_{\dot{B}^{\frac{d}{2}+s_{0},\frac{d}{2}+1}}\leq C(\left\|(u_{0},\tau_{0})\right\|_{\mathcal{E}_{0}}, s_{0})(1+t)^{-\frac{s_{0}}{2}}, \quad s_{0}\in(0,2].
\]
\end{theorem}

Moreover, using Theorem \ref{case1 decay} we can derive the further decay rates of the solution of Theorem \ref{case1 theorem} with additional initial data condition.

\begin{corollary} \label{case1 further decay}
Let $(u,\tau)\in\mathcal{E}_{T}$ be the solution of Theorem \ref{case1 theorem} with $(u_{0},\tau_{0})\in\mathcal{E}_{0}$.
\item{\textbullet} Assume that $(u_{0},\tau_{0})\in L^{2}$. Then, $(u,\tau)$ is in $L^{\infty}_{T}L^{2}\cap L^{2}_{T}\dot{H}^{1}\times L^{\infty}_{T}L^{2}\cap L^{2}_{T}L^{2}$ and satisfies
\begin{align}
&\left\|u(t)\right\|_{\dot{B}^{s}}\lesssim (1+t)^{-\frac{s}{2}},\quad s\in(0,\frac{d}{2}+1], \label{case1 further decay of u}\\
\left\|\tau(t)\right\|_{L^{2}}\lesssim (1+t)^{-\frac{1}{2}},\qquad &\left\|\tau(t)\right\|_{\dot{B}^{s}}\lesssim 
\begin{cases}
(1+t)^{-\frac{1}{2}-\frac{s}{2}},\quad s\in (0,\frac{d}{2}],\\ 
(1+t)^{-\frac{1}{2}-\frac{d}{4}},\quad s\in (\frac{d}{2},\frac{d}{2}+1]. \label{case1 further decay of tau}
\end{cases}
\end{align}
\item{\textbullet} Let $d=3$ and assume that $(u_{0},\tau_{0})\in \dot{B}^{-\frac{3}{2}}_{2,\infty}\times\dot{B}^{-\frac{1}{2}}_{2,\infty}$. Then, $(u,\tau)$ has the further decay rates:
\begin{align}
&\left\|u(t)\right\|_{\dot{B}^{s}}\lesssim (1+t)^{-\frac{3}{4}-\frac{s}{2}},\quad s\in[0,\frac{5}{2}], \label{case1 further decay2 of u}\\
&\left\|\tau(t)\right\|_{\dot{B}^{s}}\lesssim
\begin{cases}
(1+t)^{-\frac{5}{4}-\frac{s}{2}},\quad s\in[0,\frac{3}{2}], \\
(1+t)^{-2},\quad s\in (\frac{3}{2},\frac{5}{2}].
\end{cases} \label{case1 further decay2 of tau}
\end{align}
\end{corollary}

\begin{remark}\upshape \label{case1 further decay remark}
Here, we compare the results of recent work \cite{Huang} and Corollary \ref{case1 further decay}. One of the main results in \cite{Huang} is following: in three dimension, for $(u_{0},\tau_{0})\in L^{1}\cap H^{3}$ with $\left\|(u_{0},\tau_{0})\right\|_{H^{3}}\leq \epsilon$,
\[
\begin{split}
&\left\|u(t)\right\|_{\dot{H}^{k}}\lesssim (1+t)^{-\frac{3}{4}-\frac{k}{2}},\quad k=0,1,2, \\
&\left\|\tau(t)\right\|_{\dot{H}^{k}}\lesssim (1+t)^{-\frac{5}{4}-\frac{k}{2}},\quad k=0,1.
\end{split}
\]
First, we can derive (\ref{case1 further decay of u}) and (\ref{case1 further decay of tau}) without $L^{1}$ assumption while $L^{2}\cap\mathcal{E}_{0}$ is larger than $H^{3}$. We also note that $L^{1}\times L^{\frac{3}{2}}$ is continuously embedded in $\dot{B}^{-\frac{3}{2}}_{2,\infty}\times\dot{B}^{-\frac{1}{2}}_{2,\infty}$, which is the initial data space of $(u_{0},\tau_{0})$ assumed to derive (\ref{case1 further decay2 of u}) and (\ref{case1 further decay2 of tau}), and $L^{\frac{3}{2}}$ is weaker low regularity condition than $L^{1}$ in terms of $\tau_{0}$. Furthermore, we improve \cite{Huang} in the sense of the embedding $\dot{B}^{k}\hookrightarrow\dot{H}^{k}$.
\end{remark}

The first step to prove Theorem \ref{case1 theorem} and Theorem \ref{case1 decay} is to use the linearized system of (\ref{our model Case I}) with  a divergence-free vector field $v$:
\eqn \label{case1 linear eq}
\begin{split} 
& u_{t}+\mathbb{P}(v\cdot\nabla u)-\nu_{1}\Delta u-\mathbb{P}\dv\tau=F, \\
& \tau_{t}+v\cdot\nabla\tau+\alpha\tau-D(u)=G,\\
&(\mathbb{P}\dv\tau)_{t}+\mathbb{P}\dv(v\cdot\nabla\tau)+\alpha\mathbb{P}\dv\tau-\frac{1}{2}\Delta u=\mathbb{P}\dv G.
\end{split}
\een

\begin{lemma} \label{case1 lemma} \upshape
Let $-\frac{d}{2}-1<s\leq \frac{d}{2}+1$,  $u_{0}\in \dot{B}^{s,\frac{d}{2}+1}$, $\tau_{0}\in \dot{B}^{s+1,\frac{d}{2}+1}$, and $v\in L^{1}_{T}\dot{B}^{\frac{d}{2}+1}$. Then, a solution $(u,\tau)$ of (\ref{case1 linear eq}) satisfies  the following bound
\begin{align*}
&\left\|u\right\|_{L^{\infty}_{T}\dot{B}^{s,\frac{d}{2}+1}}+\left\|\tau\right\|_{L^{\infty}_{T}\dot{B}^{s+1,\frac{d}{2}+1}}+\left\|u\right\|_{L^{1}_{T}\dot{B}^{s+2,\frac{d}{2}+2}}+\left\|\tau\right\|_{L^{1}_{T}\dot{B}^{s+1,\frac{d}{2}+1}} \\
& \lesssim \left\|u_{0}\right\|_{\dot{B}^{s,\frac{d}{2}+1}}+\left\|\tau_{0}\right\|_{\dot{B}^{s+1,\frac{d}{2}+1}}+\left\|F\right\|_{L^{1}_{T}\dot{B}^{s,\frac{d}{2}+1}}+\left\|G\right\|_{L^{1}_{T}\dot{B}^{s+1,\frac{d}{2}+1}} +\left\|\mathbb{B}(\nabla v,\nabla\tau)\right\|^{l}_{L^{1}_{T}\dot{B}^{s}} \\
& +\Big(\left\|u\right\|_{L^{\infty}_{T}\dot{B}^{s,\frac{d}{2}+1}}+\left\|\tau\right\|_{L^{\infty}_{T}\dot{B}^{\min{(s,\frac{d}{2})}+1,\frac{d}{2}+1}}\Big)\left\|v\right\|_{L^{1}_{T}\dot{B}^{\frac{d}{2}+1}}.
\end{align*}
\end{lemma}

Let $s=\frac{d}{2}-1$, $v=u$, $F=0$ and $-G=\mathcal{Q}(\tau,\nabla u)\simeq \tau^{2}+\tau\nabla u+(\nabla u)^{2}$. Since 
\[
\left\|\mathcal{Q}(\tau,\nabla u)\right\|_{L^{1}_{T}\dot{B}^{\frac{d}{2}, \frac{d}{2}+1}}+\left\|\mathbb{B}(\nabla u,\nabla\tau)\right\|^{l}_{L^{1}_{T}\dot{B}^{\frac{d}{2}-1}}\lesssim \|(u, \tau)\|_{\mathcal{L}_{T}}  \|(u, \tau)\|_{\mathcal{H}_{T}},
\]
a direct application of Lemma \ref{case1 lemma} yields
\[
\left\|(u,\tau)\right\|_{\mathcal{L}_{T}}+\left\|(u,\tau)\right\|_{\mathcal{H}_{T}}\lesssim \left\|(u_{0},\tau_{0})\right\|_{\mathcal{E}_{0}}+\|(u, \tau)\|_{\mathcal{L}_{T}}  \|(u, \tau)\|_{\mathcal{H}_{T}}.
\]
So, we have the desired result, $\left\|(u,\tau)\right\|_{\mathcal{E}_{T}}\lesssim \left\|(u_{0},\tau_{0})\right\|_{\mathcal{E}_{0}}$ for all $T>0$, when $\left\|(u_{0},\tau_{0})\right\|_{\mathcal{E}_{0}}$ is sufficiently small. By combining with the uniqueness part in Section \ref{case1 unique}, we can prove Theorem \ref{case1 theorem}. (The same argument will be used to the remaining cases.) The proof of Lemma \ref{case1 lemma} is also used to derive decay rates in Theorem \ref{case1 decay} and Corollary \ref{case1 further decay}.

\begin{remark}\upshape
Using the same argument  in Lemma \ref{case1 lemma}, we   bound  (\ref{case1 linear eq}) in high frequency part as 
\eqn \label{case1 high estimate2}
\begin{split}
\left\|u\right\|^{h}_{L^{\infty}_{T}\dot{B}^{\frac{d}{2}-1}}&+\left\|\tau\right\|^{h}_{L^{\infty}_{T}\dot{B}^{\frac{d}{2}}}+\left\|u\right\|^{h}_{L^{1}_{T}\dot{B}^{\frac{d}{2}+1}}+\left\|\tau\right\|^{h}_{L^{1}_{T}\dot{B}^{\frac{d}{2}}}  \lesssim \left\|u_{0}\right\|^{h}_{\dot{B}^{\frac{d}{2}-1}}+\left\|\tau_{0}\right\|^{h}_{\dot{B}^{\frac{d}{2}}}+\left\|F\right\|^{h}_{L^{1}_{T}\dot{B}^{\frac{d}{2}-1}}\\
&+\left\|G\right\|^{h}_{L^{1}_{T}\dot{B}^{\frac{d}{2}}}+\left\|\mathbb{B}(\nabla v,\nabla \tau)\right\|^{h}_{L^{1}_{T}\dot{B}^{\frac{d}{2}-1}}+\Big(\left\|u\right\|_{L^{\infty}_{T}\dot{B}^{s,\frac{d}{2}-1}}+\left\|\tau\right\|_{L^{\infty}_{T}\dot{B}^{\frac{d}{2}+1,\frac{d}{2}}}\Big)\left\|v\right\|_{L^{1}_{T}\dot{B}^{\frac{d}{2}+1}}.
\end{split}
\een
Combining (\ref{case1 high estimate2}) and (\ref{case1 low estimate}) with $s=\frac{d}{2}-1$ gives us
\eqn \label{case1 uniqueness estimate}
\begin{split}
&\left\|u\right\|_{L^{\infty}_{T}\dot{B}^{\frac{d}{2}-1}}+\left\|\tau\right\|_{L^{\infty}_{T}\dot{B}^{\frac{d}{2}}}+\left\|u\right\|_{L^{1}_{T}\dot{B}^{\frac{d}{2}+1}}+\left\|\tau\right\|_{L^{1}_{T}\dot{B}^{\frac{d}{2}}} \\
& \lesssim \left\|u_{0}\right\|_{\dot{B}^{\frac{d}{2}-1}}+\left\|\tau_{0}\right\|_{\dot{B}^{\frac{d}{2}}}+\left\|F\right\|_{L^{1}_{T}\dot{B}^{\frac{d}{2}-1}} +\left\|G\right\|_{L^{1}_{T}\dot{B}^{\frac{d}{2}}} +\Big(\left\|u\right\|_{L^{\infty}_{T}\dot{B}^{\frac{d}{2}-1}}+\left\|\tau\right\|_{L^{\infty}_{T}\dot{B}^{\frac{d}{2}}}\Big)\left\|v\right\|_{L^{1}_{T}\dot{B}^{\frac{d}{2}+1}}.
\end{split}
\een
Let $v=u$, $F=0$ and $-G=\mathcal{Q}(\tau,\nabla u)\simeq \tau^{2}+\tau\nabla u$. We also set the spaces of initial data and solutions:  
\eqn \label{case1 uniqueness space}
\mathcal{F}_{0}=\dot{B}^{\frac{d}{2}-1}\times\dot{B}^{\frac{d}{2}},\qquad \mathcal{F}_{T}=L^{\infty}_{T}\dot{B}^{\frac{d}{2}-1}\cap L^{1}_{T}\dot{B}^{\frac{d}{2}+1}\times L^{\infty}_{T}\dot{B}^{\frac{d}{2}}\cap L^{1}_{T}\dot{B}^{\frac{d}{2}}.
\een
Then, (\ref{case1 uniqueness estimate}) yields $\left\|(u,\tau)\right\|_{\mathcal{F}_{T}} \lesssim \left\|(u_{0},\tau_{0})\right\|_{\mathcal{F}_{0}}$ for all $T>0$ when $\left\|(u_{0},\tau_{0})\right\|_{\mathcal{F}_{0}}$ is sufficiently small. This corresponds to the results in \cite{Chemin, Zi-Fang-Zhang}. 
\end{remark}

\begin{theorem} \label{case1 theorem2}\upshape \cite{Chemin, Zi-Fang-Zhang}
 Let $\mathcal{Q}(\tau,\nabla u)\simeq \tau^{2}+\tau\nabla u$ and define $\mathcal{F}_{0}$ and $\mathcal{F}_{T}$ as in (\ref{case1 uniqueness space}). There exists $\epsilon>0$ such that if $(u_{0},\tau_{0})\in \mathcal{F}_{0}$ with $\left\|(u_{0},\tau_{0})\right\|_{\mathcal{F}_{0}}\leq \epsilon$, then there exists a unique solution $(u,\tau)$ of (\ref{our model Case I}) in $\mathcal{F}_{T}$  such that $\left\|(u,\tau)\right\|_{\mathcal{F}_{T}}\lesssim \left\|(u_{0},\tau_{0})\right\|_{\mathcal{F}_{0}}$ for all $T>0$.
\end{theorem}

Furthermore, we can  derive the decay rates of the solutions in $\mathcal{F}_{T}$ which are not in \cite{Chemin, Zi-Fang-Zhang}. 

\begin{corollary} \label{case1 decay2}\upshape
The solution of Theorem \ref{case1 theorem2} has the following decay rates: for all $t>0$ and $s_{0}\in(0,2]$
\[
\left\|u(t)\right\|_{\dot{B}^{\frac{d}{2}-1+s_{0},\frac{d}{2}-1}}+\left\|\tau(t)\right\|_{\dot{B}^{\frac{d}{2}+s_{0},\frac{d}{2}}}\leq C(\left\|(u_{0},\tau_{0})\right\|_{\mathcal{F}_{0}}, s_{0})(1+t)^{-\frac{s_{0}}{2}}.
\]
\end{corollary}

\subsection{Proof of Lemma \ref{case1 lemma}}

\subsubsection{\bf Low frequency part}
First of all, it follows from (\ref{case1 linear eq}) that 
\eqn \label{case1 eq1}
\frac{1}{2}\frac{d}{dt}\left\|\Delta_{j}u\right\|^{2}_{L^{2}}+\nu_{1}\left\|\Delta_{j}\Lambda u\right\|^{2}_{L^{2}}-(\Delta_{j}u,\Delta_{j}\dv \tau)=(\Delta_{j}F,\Delta_{j}u)-(\Delta_{j}(v\cdot\nabla u),\Delta_{j}u),
\een 
\eqn \label{case1 eq2}
\frac{1}{2}\frac{d}{dt}\left\|\Delta_{j}\Lambda\tau\right\|^{2}_{L^{2}}+\alpha\left\|\Delta_{j}\Lambda\tau\right\|^{2}_{L^{2}}+(\Delta_{j}\Lambda u,\Delta_{j}\Lambda\dv\tau)=(\Delta_{j}\Lambda G,\Delta_{j}\Lambda\tau)-(\Delta_{j}\Lambda(v\cdot\nabla\tau),\Delta_{j}\Lambda\tau), 
\een 
and 
\eqn \label{case1 eq3}
\begin{split}
&\frac{d}{dt}(\Delta_{j}u,\Delta_{j}\mathbb{P}\dv\tau)-\left\|\Delta_{j}\mathbb{P}\dv\tau\right\|^{2}_{L^{2}}+\frac{1}{2}\left\|\Delta_{j}\Lambda u\right\|^{2}_{L^{2}}+\alpha(\Delta_{j}u,\Delta_{j}\dv\tau)+\nu_{1}(\Delta_{j}\Lambda^{2}u,\Delta_{j}\dv\tau) \\
&=(\Delta_{j}F,\Delta_{j}\mathbb{P}\dv\tau)+(\Delta_{j}\mathbb{P}\dv G,\Delta_{j}u)-(\Delta_{j}(v\cdot\nabla u),\Delta_{j}\mathbb{P}\dv\tau)-(\Delta_{j}\mathbb{P}\dv(v\cdot\nabla\tau),\Delta_{j}u).
\end{split}
\een
We now proceed to add $\alpha$(\ref{case1 eq1}), $K_{1}$(\ref{case1 eq2}), and (\ref{case1 eq3}) to deduce that 
\[
\frac{1}{2}\frac{d}{dt}f^{2}_{j}+h^{2}_{j} =\widetilde{F}^{1}_{j}-\widetilde{F}^{2}_{j}-\widetilde{F}^{3}_{j}-\widetilde{F}^{4}_{j} \qquad\text{for}\ j \leq j_{0},
\]
where $K_{1}>0$ is determined below,  and
\begin{align*}
f^{2}_{j}&=\alpha\left\|\Delta_{j}u\right\|^{2}_{L^{2}}+K_{1}\left\|\Delta_{j}\Lambda\tau\right\|^{2}_{L^{2}}+2(\Delta_{j}u,\Delta_{j}\mathbb{P}\dv\tau), \\
h^{2}_{j}&=\Big(\alpha\nu_{1}+\frac{1}{2}\Big)\left\|\Delta_{j}\Lambda u\right\|^{2}_{L^{2}}+\alpha K_{1}\left\|\Delta_{j}\Lambda\tau\right\|^{2}_{L^{2}}-\left\|\Delta_{j}\mathbb{P}\dv\tau\right\|^{2}_{L^{2}} +\left(K_{1}+\nu_{1}\right)(\Delta_{j}\Lambda^{2}u,\Delta_{j}\dv\tau),\\
\widetilde{F}^{1}_{j}&=\alpha(\Delta_{j}F,\Delta_{j}u)+K_{1}(\Delta_{j}\Lambda G,\Delta_{j}\Lambda\tau)+(\Delta_{j}F,\Delta_{j}\mathbb{P}\dv\tau)+(\Delta_{j}\mathbb{P}\dv G,\Delta_{j}u), \\
\widetilde{F}^{2}_{j}&=\alpha(\Delta_{j}(v\cdot\nabla u),\Delta_{j}u)+(\Delta_{j}(v\cdot\nabla u),\Delta_{j}\mathbb{P}\dv\tau)+(\Delta_{j}(v\cdot\nabla\mathbb{P}\dv\tau),\Delta_{j}u), \\
\widetilde{F}^{3}_{j}&=K_{1}(\Delta_{j}\Lambda(v\cdot\nabla\tau),\Delta_{j}\Lambda\tau), \\
\widetilde{F}^{4}_{j}&=(\Delta_{j}\mathbb{B}(\nabla v,\nabla\tau),\Delta_{j}u),
\end{align*}
where we apply Proposition \ref{decompose prop} to $\mathbb{P}\dv(v\cdot\nabla\tau)$. By choosing $K_{1}=\frac{2C^{2}_{0}}{\alpha}$ where $C_{0}$ is defined in (\ref{C0}), and by taking a sufficiently small $j_{0}$, we find   
\[
f^{2}_{j}\simeq \left\|\Delta_{j}u\right\|^{2}_{L^{2}}+\left\|\Delta_{j}\Lambda\tau\right\|^{2}_{L^{2}}, \quad h^{2}_{j}\simeq \left\|\Delta_{j}\Lambda u\right\|^{2}_{L^{2}}+\left\|\Delta_{j}\Lambda\tau\right\|^{2}_{L^{2}} \gtrsim 2^{2j}f^{2}_{j}.
\]
We now bound $\widetilde{F}_{j}$ terms using (\ref{convection term estimate}) and (\ref{convection term estimate 2})
\begin{align*}
&\left|\widetilde{F}^{1}_{j}\right|+\left|\widetilde{F}^{4}_{j}\right| \lesssim c_{j}2^{-sj}\left(\left\|F\right\|^{l}_{\dot{B}^{s}}+\left\|G\right\|^{l}_{\dot{B}^{s+1}}+\left\|\mathbb{B}(\nabla v,\nabla \tau)\right\|^{l}_{\dot{B}^{s}}\right)f_{j}, \\
&\left|\widetilde{F}^{2}_{j}\right| \lesssim c_{j}2^{-sj}\left\|v\right\|_{\dot{B}^{\frac{d}{2}+1}}\left(\left\|u\right\|_{\dot{B}^{s,\frac{d}{2}-1}}+\left\|\tau\right\|_{\dot{B}^{s+1,\frac{d}{2}}}\right) f_{j}, \\
&\left|\widetilde{F}^{3}_{j}\right| \lesssim c_{j}2^{-\min{(s,\frac{d}{2})}j}\left\|v\right\|_{\dot{B}^{\frac{d}{2}+1}}\left\|\tau\right\|_{\dot{B}^{\min{(s,\frac{d}{2})}+1,\frac{d}{2}}} \left\|\Delta_{j}\Lambda\tau\right\|_{L^{2}} \lesssim c_{j}2^{-sj}\left\|v\right\|_{\dot{B}^{\frac{d}{2}+1}}\left\|\tau\right\|_{\dot{B}^{\min{(s,\frac{d}{2})}+1,\frac{d}{2}}} f_{j}
\end{align*}
whenever $j\leq j_{0}$ and $-\frac{d}{2}-1<s\leq \frac{d}{2}+1$. From these bounds, we obtain  
\[
\begin{split}
\left\|u\right\|^{l}_{L^{\infty}_{T}\dot{B}^{s}}&+\left\|\tau\right\|^{l}_{L^{\infty}_{T}\dot{B}^{s+1}}+\left\|u\right\|^{l}_{L^{1}_{T}\dot{B}^{s+2}}+\left\|\tau\right\|^{l}_{L^{1}_{T}\dot{B}^{s+3}}  \lesssim \left\|u_{0}\right\|^{l}_{\dot{B}^{s}}+\left\|\tau_{0}\right\|^{l}_{\dot{B}^{s+1}}+\left\|F\right\|^{l}_{L^{1}_{T}\dot{B}^{s}}\\
&+\left\|G\right\|^{l}_{L^{1}_{T}\dot{B}^{s+1}}+\left\|\mathbb{B}(\nabla v,\nabla \tau)\right\|^{l}_{L^{1}_{T}\dot{B}^{s}}+\Big(\left\|u\right\|_{L^{\infty}_{T}\dot{B}^{s,\frac{d}{2}-1}}+\left\|\tau\right\|_{L^{\infty}_{T}\dot{B}^{\min{(s,\frac{d}{2})}+1,\frac{d}{2}}}\Big)\left\|v\right\|_{L^{1}_{T}\dot{B}^{\frac{d}{2}+1}}.
\end{split}
\]
From (\ref{case1 eq2}), we also have
\[
\left\|\tau\right\|^{l}_{L^{\infty}_{T}\dot{B}^{s+1}}+\left\|\tau\right\|^{l}_{L^{1}_{T}\dot{B}^{s+1}} \lesssim \left\|\tau_{0}\right\|^{l}_{\dot{B}^{s+1}}+\left\|u\right\|^{l}_{L^{1}_{T}\dot{B}^{s+2}}+\left\|G\right\|^{l}_{L^{1}_{T}\dot{B}^{s+1}}+\left\|\tau\right\|_{L^{\infty}_{T}\dot{B}^{\min{(s,\frac{d}{2})}+1,\frac{d}{2}}}\left\|v\right\|_{L^{1}_{T}\dot{B}^{\frac{d}{2}+1}}.
\]
Combining these two inequalities and using (\ref{hybrid Besov embedding}), we arrive at  
\eqn \label{case1 low estimate}
\begin{split}
&\left\|u\right\|^{l}_{L^{\infty}_{T}\dot{B}^{s}}+\left\|\tau\right\|^{l}_{L^{\infty}_{T}\dot{B}^{s+1}}+\left\|u\right\|^{l}_{L^{1}_{T}\dot{B}^{s+2}}+\left\|\tau\right\|^{l}_{L^{1}_{T}\dot{B}^{s+1}}  \lesssim \left\|u_{0}\right\|^{l}_{\dot{B}^{s}}+\left\|\tau_{0}\right\|^{l}_{\dot{B}^{s+1}}+\left\|F\right\|^{l}_{L^{1}_{T}\dot{B}^{s}}\\
&+\left\|G\right\|^{l}_{L^{1}_{T}\dot{B}^{s+1}}+\left\|\mathbb{B}(\nabla v,\nabla \tau)\right\|^{l}_{L^{1}_{T}\dot{B}^{s}}+\Big(\left\|u\right\|_{L^{\infty}_{T}\dot{B}^{s,\frac{d}{2}-1}}+\left\|\tau\right\|_{L^{\infty}_{T}\dot{B}^{\min{(s,\frac{d}{2})}+1,\frac{d}{2}}}\Big)\left\|v\right\|_{L^{1}_{T}\dot{B}^{\frac{d}{2}+1}}.
\end{split}
\een

\subsubsection{\bf High frequency part}
In going to bound the high frequency part of $(u, \tau)$, we first see 
\eqn \label{case1 eq4}
\frac{1}{2}\frac{d}{dt}\left\|\Delta_{j}\tau\right\|^{2}_{L^{2}}+\alpha\left\|\Delta_{j}\tau\right\|^{2}_{L^{2}}+(\Delta_{j}u,\Delta_{j}\dv \tau)=(\Delta_{j}G,\Delta_{j}\tau)-(\Delta_{j}(v\cdot\nabla\tau),\Delta_{j}\tau)
\een
which has one regularity less than (\ref{case1 eq2}). By  $(\ref{case1 eq1})+(\ref{case1 eq4})$, we derive the following for $ j\geq j_{0}+1$
\eqn \label{case1 high eq}
\frac{1}{2}\frac{d}{dt}f^{2}_{j}+h^{2}_{j}=\widetilde{G}^{1}_{j}-\widetilde{G}^{2}_{j},
\een
where $j_{0}$ is fixed above and 
\[
\begin{split}
f^{2}_{j} &=\left\|\Delta_{j}u\right\|^{2}_{L^{2}}+\left\|\Delta_{j}\tau\right\|^{2}_{L^{2}}, \quad h^{2}_{j} =\nu_{1}\left\|\Delta_{j}\Lambda u\right\|^{2}_{L^{2}}+\alpha\left\|\Delta_{j}\tau\right\|^{2}_{L^{2}}\gtrsim f^{2}_{j},\\
\widetilde{G}^{1}_{j} &=(\Delta_{j}F,\Delta_{j}u)+(\Delta_{j}G,\Delta_{j}\tau), \quad \widetilde{G}^{2}_{j}=(\Delta_{j}(v\cdot\nabla u),\Delta_{j}u)+(\Delta_{j}(v\cdot\nabla\tau),\Delta_{j}\tau).
\end{split}
\]
Using (\ref{convection term estimate}), we bound the right-hand side of (\ref{case1 high eq}) as follows:
\[
\left|\widetilde{G}^{1}_{j}\right|\lesssim c_{j}2^{-\left(\frac{d}{2}+1\right)j}\left\|(F,G)\right\|^{h}_{\dot{B}^{\frac{d}{2}+1}} f_{j},\quad \left|\widetilde{G}^{2}_{j}\right|\lesssim c_{j}2^{-\left(\frac{d}{2}+1\right)j}\left\|v\right\|_{\dot{B}^{\frac{d}{2}+1}}\Big(\left\|u\right\|_{\dot{B}^{s,\frac{d}{2}+1}}+\left\|\tau\right\|_{\dot{B}^{\frac{d}{2}+1}}\Big)f_{j}
\]
when $j\ge j_{0}+1$. So, we obtain  
\eqn \label{eq:3.12}
\begin{split}
& \left\|u\right\|^{h}_{L^{\infty}_{T}\dot{B}^{\frac{d}{2}+1}}+\left\|\tau\right\|^{h}_{L^{\infty}_{T}\dot{B}^{\frac{d}{2}+1}}+\left\|u\right\|^{h}_{L^{1}_{T}\dot{B}^{\frac{d}{2}+1}}+\left\|\tau\right\|^{h}_{L^{1}_{T}\dot{B}^{\frac{d}{2}+1}} \\
& \lesssim \left\|u_{0}\right\|^{h}_{\dot{B}^{\frac{d}{2}+1}}+\left\|\tau_{0}\right\|^{h}_{\dot{B}^{\frac{d}{2}+1}}+\left\|(F,G)\right\|^{h}_{L^{1}_{T}\dot{B}^{\frac{d}{2}+1}} +\Big(\left\|u\right\|_{L^{\infty}_{T}\dot{B}^{s,\frac{d}{2}+1}}+\left\|\tau\right\|_{L^{\infty}_{T}\dot{B}^{\frac{d}{2}+1}}\Big)\left\|v\right\|_{L^{1}_{T}\dot{B}^{\frac{d}{2}+1}}.
\end{split}
\een 
From (\ref{case1 eq1}), we also have 
\eqn \label{eq:3.13}
\left\|u\right\|^{h}_{L^{\infty}_{T}\dot{B}^{\frac{d}{2}}}+\left\|u\right\|^{h}_{L^{1}_{T}\dot{B}^{\frac{d}{2}+2}}\lesssim \left\|u_{0}\right\|^{h}_{\dot{B}^{\frac{d}{2}}}+\left\|\tau\right\|^{h}_{L^{1}_{T}\dot{B}^{\frac{d}{2}+1}}+\left\|F\right\|^{h}_{L^{1}_{T}\dot{B}^{\frac{d}{2}}}+\left\|u\right\|_{L^{\infty}_{T}\dot{B}^{s,\frac{d}{2}}}\left\|v\right\|_{L^{1}_{T}\dot{B}^{\frac{d}{2}+1}}.
\een
The two bounds (\ref{eq:3.12}) and (\ref{eq:3.13}), and (\ref{hybrid Besov embedding}) give  
\eqn \label{case1 high estimate}
\begin{split}
& \left\|u\right\|^{h}_{L^{\infty}_{T}\dot{B}^{\frac{d}{2}+1}}+\left\|\tau\right\|^{h}_{L^{\infty}_{T}\dot{B}^{\frac{d}{2}+1}}+\left\|u\right\|^{h}_{L^{1}_{T}\dot{B}^{\frac{d}{2}+2}}+\left\|\tau\right\|^{h}_{L^{1}_{T}\dot{B}^{\frac{d}{2}+1}} \lesssim \left\|u_{0}\right\|^{h}_{\dot{B}^{\frac{d}{2}+1}}+\left\|\tau_{0}\right\|^{h}_{\dot{B}^{\frac{d}{2}+1}}\\
& +\left\|(F,G)\right\|^{h}_{L^{1}_{T}\dot{B}^{\frac{d}{2}+1}}+\Big(\left\|u\right\|_{L^{\infty}_{T}\dot{B}^{s,\frac{d}{2}+1}}+\left\|\tau\right\|_{L^{\infty}_{T}\dot{B}^{\frac{d}{2}+1}}\Big)\left\|v\right\|_{L^{1}_{T}\dot{B}^{\frac{d}{2}+1}}.
\end{split}
\een
By (\ref{case1 low estimate}) and (\ref{case1 high estimate}), we complete the proof of Lemma \ref{case1 lemma}.

\subsection{Uniqueness} \label{case1 unique}
Suppose there are two solutions $(u_{1},\tau_{1})$ and $(u_{2},\tau_{2})$  of (\ref{our model Case I}) with the same initial data. Let $(u, \tau)=(u_{1}-u_{2},\tau_{1}-\tau_{2})$. Then, $(u, \tau)$ satisfies $\dv u_{2}=0$ and 
\eqn \label{case1 difference eq}
\begin{split}
& u_{t}+\mathbb{P}(u_{2}\cdot\nabla u)-\nu_{1}\Delta u-\mathbb{P}\dv \tau=-\mathbb{P}(u\cdot\nabla u_{1})=\delta F,\\
& \tau_{t}+u_{2}\cdot\nabla\tau+\alpha\tau-D(u)=-u\cdot\nabla\tau_{1}-\left[\mathcal{Q}(\tau_{1},\nabla u_{1})-\mathcal{Q}(\tau_{2},\nabla u_{2})\right]=\delta G.
\end{split}
\een
We now show the uniqueness  in $\mathcal{F}_{T}$ defined in (\ref{case1 uniqueness space}). By applying (\ref{case1 uniqueness estimate}) to (\ref{case1 difference eq}), we have
\[
\left\|(u, \tau)\right\|_{\mathcal{F}_{T}} \lesssim \left\|\delta F\right\|_{L^{1}_{T}\dot{B}^{\frac{d}{2}-1}}+\left\|\delta G\right\|_{L^{1}_{T}\dot{B}^{\frac{d}{2}}} +\left\|u_{2}\right\|_{L^{1}_{T}\dot{B}^{\frac{d}{2}+1}}\left\|(u, \tau)\right\|_{\mathcal{F}_{T}}.
\]
We now bound the nonlinear terms. We first deal with $u\cdot\nabla u_{1}$ and $u\cdot\nabla \tau_{1}$ by using (\ref{interpolation}) and (\ref{product estimate}): 
\[
\begin{split}
\left\|u\cdot\nabla u_{1}\right\|_{L^{1}_{T}\dot{B}^{\frac{d}{2}-1}}+ \left\|u\cdot\nabla \tau_{1}\right\|_{L^{1}_{T}\dot{B}^{\frac{d}{2}}} &\lesssim \left\|u_{1}\right\|_{L^{1}_{T}\dot{B}^{\frac{d}{2}+1}}\|u\|_{L^{\infty}_{T}\dot{B}^{\frac{d}{2}-1}}+\left\|\tau_{1}\right\|_{L^{2}_{T}\dot{B}^{\frac{d}{2}+1}}\|u\|_{L^{2}_{T}\dot{B}^{\frac{d}{2}}} \\
& \lesssim \left\|(u_{1},\tau_{1})\right\|_{\mathcal{E}_{T}}\left\|(u,\tau)\right\|_{\mathcal{F}_{T}}.
\end{split}
\]
Since $\mathcal{Q}(\tau_{1},\nabla u_{1})-\mathcal{Q}(\tau_{2},\nabla u_{2}) \simeq (\tau_{1}+\tau_{2})\tau+\tau\nabla u_{1}+\tau_{2}\nabla  u+(\nabla u_{1}+\nabla u_{2})\nabla u$,   
\[
\left\|\mathcal{Q}(\tau_{1},\nabla u_{1})-\mathcal{Q}(\tau_{2},\nabla u_{2})\right\|_{L^{1}_{T}\dot{B}^{\frac{d}{2}}}\lesssim \Big(\left\|(u_{1},\tau_{1})\right\|_{\mathcal{E}_{T}}+\left\|(u_{2},\tau_{2})\right\|_{\mathcal{E}_{T}}\Big)\left\|(u, \tau)\right\|_{\mathcal{F}_{T}}.
\]
By all these bounds together, we deduce that 
\[
\left\|(u, \tau)\right\|_{\mathcal{F}_{T}}\lesssim \Big(\left\|(u_{1},\tau_{1})\right\|_{\mathcal{E}_{T}}+\left\|(u_{2},\tau_{2})\right\|_{\mathcal{E}_{T}}\Big)\left\|(u, \tau)\right\|_{\mathcal{F}_{T}}.
\]
Since $\left\|(u_{1},\tau_{1})\right\|_{\mathcal{E}_{T}}+\left\|(u_{2},\tau_{2})\right\|_{\mathcal{E}_{T}}\lesssim\epsilon$, $(u,\tau)=0$ in $\mathcal{F}_{T}$, we  complete the uniqueness part.

\subsection{Proof of Theorem \ref{case1 decay}}
Let $v=u$, $F=0$, and $-G=\mathcal{Q}(\tau,\nabla u)\simeq \tau^{2}+\tau\nabla u+(\nabla u)^{2}$. If we do not take time integrations when bounding solutions in the proof of Lemma \ref{case1 lemma}, we have  
\eqn \label{case1 ineq for decay}
\begin{split}
& \frac{d}{dt}\left(\left\|u\right\|_{\dot{B}^{s,\frac{d}{2}+1}}+\left\|\tau\right\|_{\dot{B}^{s+1,\frac{d}{2}+1}}\right)+\left\|u\right\|_{\dot{B}^{s+2,\frac{d}{2}+2}}+\left\|\tau\right\|_{\dot{B}^{s+1,\frac{d}{2}+1}} \\
&\lesssim \left\|\mathcal{Q}(\tau,\nabla u)\right\|_{\dot{B}^{s+1,\frac{d}{2}+1}}+\left\|\mathbb{B}(\nabla u,\nabla \tau)\right\|^{l}_{\dot{B}^{s}}+\left\|u\right\|_{\dot{B}^{\frac{d}{2}+1}}\left(\left\|u\right\|_{\dot{B}^{s,\frac{d}{2}+1}}+\left\|\tau\right\|_{\dot{B}^{\min{(s,\frac{d}{2})}+1,\frac{d}{2}+1}}\right).
\end{split}
\een
From (\ref{hybrid Besov embedding}) and (\ref{product estimate}),  we bound two nonlinear terms
\[
\begin{split}
\left\|\mathcal{Q}(\tau,\nabla u)\right\|_{\dot{B}^{s+1,\frac{d}{2}+1}}&\lesssim \Big(\left\|u\right\|_{\dot{B}^{\frac{d}{2}+1}}+\left\|\tau\right\|_{\dot{B}^{\frac{d}{2}}}\Big)\left(\left\|u\right\|_{\dot{B}^{s+2,\frac{d}{2}+2}}+\left\|\tau\right\|_{\dot{B}^{s+1,\frac{d}{2}+1}}\right),\\
\left\|\mathbb{B}(\nabla u,\nabla \tau)\right\|^{l}_{\dot{B}^{s}}&\lesssim \left\|\mathbb{B}(\nabla u,\nabla \tau)\right\|_{\dot{B}^{\min{(s,\frac{d}{2})},\frac{d}{2}}}\lesssim \left\|u\right\|_{\dot{B}^{\frac{d}{2}+1}}\left\|\tau\right\|_{\dot{B}^{\min{(s,\frac{d}{2})}+1,\frac{d}{2}+1}}
\end{split}
\]
for $s>-\frac{d}{2}$. Thus, using (\ref{hybrid Besov embedding}) we can rearrange (\ref{case1 ineq for decay}) to get 
\eqn \label{eq:3.15}
\begin{split}
& \frac{d}{dt}\left(\left\|u\right\|_{\dot{B}^{s,\frac{d}{2}+1}}+\left\|\tau\right\|_{\dot{B}^{s+1,\frac{d}{2}+1}}\right)+\left\|u\right\|_{\dot{B}^{s+2,\frac{d}{2}+2}}+\left\|\tau\right\|_{\dot{B}^{s+1,\frac{d}{2}+1}} \\
&\lesssim \Big(\left\|u\right\|_{\dot{B}^{\frac{d}{2}+1}}+\left\|\tau\right\|_{\dot{B}^{\frac{d}{2}}}\Big)\left(\left\|u\right\|_{\dot{B}^{s+2,\frac{d}{2}+2}}+\left\|\tau\right\|_{\dot{B}^{s+1,\frac{d}{2}+1}}\right)+\left\|u\right\|_{\dot{B}^{\frac{d}{2}+1}}\left\|u\right\|^{l}_{\dot{B}^{s}}+\left\|u\right\|_{\dot{B}^{\frac{d}{2}+1}}\left\|\tau\right\|^{l}_{\dot{B}^{\frac{d}{2}+1}}.
\end{split}
\een
Since $\left\|(u,\tau)\right\|_{\mathcal{E}_{T}}\lesssim\epsilon$, the first term of the right-hand side of (\ref{eq:3.15}) can be absorbed into the left-hand side of (\ref{eq:3.15}). We now bound the remaining terms when $\frac{d}{2}-1\leq s\leq \frac{d}{2}+1$. First of all, when $s\geq \frac{d}{2}-1$,
\eqn \label{bound for decay 1}
\begin{split}
\left\|u\right\|_{\dot{B}^{\frac{d}{2}+1}}\left\|u\right\|^{l}_{\dot{B}^{s}} &=\left\|u\right\|^{l}_{\dot{B}^{\frac{d}{2}+1}}\left\|u\right\|^{l}_{\dot{B}^{s}}+\left\|u\right\|^{h}_{\dot{B}^{\frac{d}{2}+1}}\left\|u\right\|^{l}_{\dot{B}^{s}} \\
& \lesssim \left\|u\right\|^{l}_{\dot{B}^{\frac{d}{2}-1}}\left\|u\right\|^{l}_{\dot{B}^{s+2}} +\left\|u\right\|^{h}_{\dot{B}^{\frac{d}{2}+1}}\left\|u\right\|^{l}_{\dot{B}^{s}} \lesssim \left\|u\right\|^{l}_{\dot{B}^{\frac{d}{2}-1}}\left\|u\right\|_{\dot{B}^{s+2,\frac{d}{2}+1}}. 
\end{split}
\een
When $\frac{d}{2}-1\leq s\leq \frac{d}{2}$, we have
\[
\left\|u\right\|_{\dot{B}^{\frac{d}{2}+1}}\left\|\tau\right\|^{l}_{\dot{B}^{\frac{d}{2}+1}}=\left\|u\right\|^{l}_{\dot{B}^{\frac{d}{2}+1}}\left\|\tau\right\|^{l}_{\dot{B}^{\frac{d}{2}+1}}+\left\|u\right\|^{h}_{\dot{B}^{\frac{d}{2}+1}}\left\|\tau\right\|^{l}_{\dot{B}^{\frac{d}{2}+1}}\lesssim \left\|u\right\|^{l}_{\dot{B}^{\frac{d}{2}-1}}\left\|\tau\right\|^{l}_{\dot{B}^{s+1}}+\left\|\tau\right\|^{l}_{\dot{B}^{\frac{d}{2}}}\left\|u\right\|^{h}_{\dot{B}^{\frac{d}{2}+1}}.
\]
When $\frac{d}{2}<s\leq\frac{d}{2}+1$, we deal with  $\left\|u\right\|_{\dot{B}^{\frac{d}{2}+1}}\left\|\tau\right\|^{l}_{\dot{B}^{\frac{d}{2}+1}}$ as follows:
\begin{align*}
\left\|u\right\|_{\dot{B}^{\frac{d}{2}+1}}\left\|\tau\right\|^{l}_{\dot{B}^{\frac{d}{2}+1}} &\lesssim \Big(\left\|u\right\|^{l}_{\dot{B}^{\frac{d}{2}-1}}\Big)^{1-\theta}\left(\left\|u\right\|^{l}_{\dot{B}^{s+2}}\right)^{\theta}\Big(\left\|\tau\right\|^{l}_{\dot{B}^{\frac{d}{2}}}\Big)^{\beta}\left(\left\|\tau\right\|^{l}_{\dot{B}^{s+1}}\right)^{1-\beta} +\left\|u\right\|^{h}_{\dot{B}^{\frac{d}{2}+1}}\left\|\tau\right\|^{l}_{\dot{B}^{\frac{d}{2}+1}}\\
&\lesssim \Big(\left\|\tau\right\|^{l}_{\dot{B}^{\frac{d}{2}}}\left\|u\right\|^{l}_{\dot{B}^{s+2}}\Big)^{\beta}\Big(\left\|u\right\|^{l}_{\dot{B}^{\frac{d}{2}-1}}\left\|\tau\right\|^{l}_{\dot{B}^{s+1}}\Big)^{1-\beta} +\left\|u\right\|^{h}_{\dot{B}^{\frac{d}{2}+1}}\left\|\tau\right\|^{l}_{\dot{B}^{\frac{d}{2}+1}}\\
&\lesssim \left\|\tau\right\|^{l}_{\dot{B}^{\frac{d}{2}}}\left\|u\right\|^{l}_{\dot{B}^{s+2}} +\left\|u\right\|^{l}_{\dot{B}^{\frac{d}{2}-1}}\left\|\tau\right\|^{l}_{\dot{B}^{s+1}}+\left\|\tau\right\|^{l}_{\dot{B}^{\frac{d}{2}}}\left\|u\right\|^{h}_{\dot{B}^{\frac{d}{2}+1}} \\
&\lesssim \Big(\left\|u\right\|^{l}_{\dot{B}^{\frac{d}{2}-1}}+\left\|\tau\right\|^{l}_{\dot{B}^{\frac{d}{2}}}\Big)\left(\left\|u\right\|_{\dot{B}^{s+2,\frac{d}{2}+1}}+\left\|\tau\right\|^{l}_{\dot{B}^{s+1}}\right),
\end{align*}
where $\theta=\frac{2}{s+3-d/2}$, $\beta=\frac{s-d/2}{s+1-d/2}$ and we use $0<\beta\leq\theta<1$ when $\frac{d}{2}<s\leq\frac{d}{2}+1$. That is, for $\frac{d}{2}-1\leq s\leq \frac{d}{2}+1$,
\eqn \label{bound for decay 2}
\left\|u\right\|_{\dot{B}^{\frac{d}{2}+1}}\left\|\tau\right\|^{l}_{\dot{B}^{\frac{d}{2}+1}}\lesssim \Big(\left\|u\right\|^{l}_{\dot{B}^{\frac{d}{2}-1}}+\left\|\tau\right\|^{l}_{\dot{B}^{\frac{d}{2}}}\Big)\left(\left\|u\right\|_{\dot{B}^{s+2,\frac{d}{2}+1}}+\left\|\tau\right\|^{l}_{\dot{B}^{s+1}}\right).
\een
Therefore, all the terms of right-hand side of (\ref{eq:3.15}) can be absorbed into the left-hand side of (\ref{eq:3.15}) and we have 
\eqn \label{eq:3.18}
\frac{d}{dt}\left(\left\|u\right\|_{\dot{B}^{s,\frac{d}{2}+1}}+\left\|\tau\right\|_{\dot{B}^{s+1,\frac{d}{2}+1}}\right)+\left\|u\right\|_{\dot{B}^{s+2,\frac{d}{2}+2}}+\left\|\tau\right\|_{\dot{B}^{s+1,\frac{d}{2}+1}}\leq 0, \quad \frac{d}{2}-1\leq s\leq \frac{d}{2}+1.
\een
To make (\ref{eq:3.18}) as an ODE of the form $f'+cf^{m}\leq 0$ with $m>1$, we need a lower bound of $\left\|u\right\|^{l}_{\dot{B}^{s+2}}$ in terms of $\left\|u\right\|^{l}_{\dot{B}^{s}}$. Using (\ref{interpolation}) and letting $s=\frac{d}{2}-1+s_{0}$, $s_{0}\in(0,2]$, we have 
\[
\left\|u\right\|^{l}_{\dot{B}^{s}}\leq \Big(\left\|u\right\|^{l}_{\dot{B}^{\frac{d}{2}-1}}\Big)^{\theta}\left(\left\|u\right\|^{l}_{\dot{B}^{s+2}}\right)^{1-\theta}\lesssim \epsilon^{\frac{2}{s_{0}+2}} \left(\left\|u\right\|^{l}_{\dot{B}^{s+2}}\right)^{\frac{s_{0}}{s_{0}+2}}
\]
which yields that
\[
\frac{d}{dt}\left(\left\|u\right\|_{\dot{B}^{s,\frac{d}{2}+1}}+\left\|\tau\right\|_{\dot{B}^{s+1,\frac{d}{2}+1}}\right)+\epsilon^{-\frac{2}{s_{0}}}\left(\left\|u\right\|^{l}_{\dot{B}^{s}}\right)^{1+\frac{2}{s_{0}}} +\left\|u\right\|^{h}_{\dot{B}^{\frac{d}{2}+1}}+\left\|\tau\right\|_{\dot{B}^{s+1,\frac{d}{2}+1}}\leq 0.
\]
Again using the fact that $\left\|(u,\tau)\right\|_{\mathcal{E}_{T}}\lesssim\epsilon$, we have
\[
\frac{d}{dt}\left(\left\|u\right\|_{\dot{B}^{\frac{d}{2}-1+s_{0},\frac{d}{2}+1}}+\left\|\tau\right\|_{\dot{B}^{\frac{d}{2}+s_{0},\frac{d}{2}+1}}\right)+C(\epsilon)\left(\left\|u\right\|_{\dot{B}^{\frac{d}{2}-1+s_{0},\frac{d}{2}+1}}+\left\|\tau\right\|_{\dot{B}^{\frac{d}{2}+s_{0},\frac{d}{2}+1}}\right)^{1+\frac{2}{s_{0}}}\leq 0.
\]
By solving this inequality, we can obtain the decay rates in Theorem \ref{case1 decay}.

\subsection{Proof of Corollary \ref{case1 further decay}}
\textbullet\ If $(u_{0},\tau_{0})\in L^{2}$, using the fact $\left\|(u,\tau)\right\|_{\mathcal{E}_{T}}\lesssim \epsilon$, it is easy to check that $(u,\tau)\in L^{\infty}_{T}L^{2}\cap L^{2}_{T}\dot{H}^{1}\times L^{\infty}_{T}L^{2}\cap L^{2}_{T}L^{2}$ with $\left\|(u,\tau)(t)\right\|_{L^{2}}\leq \left\|(u_{0},\tau_{0})\right\|_{L^{2}}$ for all $t$. Then, (\ref{eq:3.18}) with (\ref{interpolation 2}) gives us 
\[
\left\|u(t)\right\|_{\dot{B}^{s,\frac{d}{2}+1}}+\left\|\tau(t)\right\|_{\dot{B}^{s+1,\frac{d}{2}+1}}\lesssim (1+t)^{-\frac{s}{2}},\quad s\in (\frac{d}{2}-1,\frac{d}{2}+1].
\]
Using (\ref{interpolation 2}) again, we obtain (\ref{case1 further decay of u}). From $\left\|(u,\tau)\right\|_{\mathcal{E}_{T}}\lesssim \epsilon$ and (\ref{case1 further decay of u}), we have 
\[
\frac{d}{dt}\left\|\tau(t)\right\|_{L^{2}}+\frac{\alpha}{2}\left\|\tau(t)\right\|_{L^{2}}\leq C\left\|\nabla u(t)\right\|_{L^{2}}\leq C(1+t)^{-\frac{1}{2}},
\]
\[
\frac{d}{dt}\left\|\tau(t)\right\|_{\dot{B}^{\frac{d}{2}}}+\frac{\alpha}{2}\left\|\tau(t)\right\|_{\dot{B}^{\frac{d}{2}}}\leq C\left\|u(t)\right\|_{\dot{B}^{\frac{d}{2}+1}}\leq C(1+t)^{-\frac{1}{2}-\frac{d}{4}},
\]
which implies
\[
\left\|\tau(t)\right\|_{L^{2}}\lesssim (1+t)^{-\frac{1}{2}},\quad \left\|\tau(t)\right\|_{\dot{B}^{\frac{d}{2}}}\lesssim (1+t)^{-\frac{1}{2}-\frac{d}{4}}.
\]
Then, (\ref{case1 further decay of tau}) is derived by (\ref{interpolation}) and (\ref{interpolation 2}).

\textbullet\ When $(u_{0},\tau_{0})\in \dot{B}^{-\frac{3}{2}}_{2,\infty}\times\dot{B}^{-\frac{1}{2}}_{2,\infty}$, we use a variant of Fourier splitting method \cite{Schonbek} in this framework. To this end, we need the following two estimates arising in low frequency part:
\eqn \label{case1 Fourier split ineq1}
\begin{split}
&\frac{d}{dt}f_{j}+2^{2j}f_{j}\leq X_{j},\qquad j\leq j_{0},\\
f_{j}\simeq \left\|\Delta_{j}u\right\|_{L^{2}}+\left\|\Delta_{j}\Lambda\tau\right\|_{L^{2}},\quad &X_{j}\lesssim \left\|\Delta_{j}(u\cdot\nabla u)\right\|_{L^{2}}+\left\|\Delta_{j}\Lambda(u\cdot\nabla\tau)\right\|_{L^{2}}+\left\|\Delta_{j}\Lambda\mathcal{Q}(\tau,\nabla u)\right\|_{L^{2}},
\end{split}
\een
and 
\eqn \label{case1 Fourier split ineq2}
\begin{split}
\frac{d}{dt}\left(\left\|u\right\|^{l}_{\dot{B}^{0}}+\left\|\tau\right\|^{l}_{\dot{B}^{1}}\right)+\left\|u\right\|^{l}_{\dot{B}^{2}}+\left\|\tau\right\|^{l}_{\dot{B}^{1}} \lesssim \left\|\mathcal{Q}(\tau,\nabla u)\right\|^{l}_{\dot{B}^{1}}+\left(\left\|u\right\|_{\dot{B}^{0}}+\left\|\tau\right\|_{\dot{B}^{1}}\right)\left\|u\right\|_{\dot{B}^{\frac{5}{2}}}.
\end{split}
\een
Let $r(t)$ satisfying $2^{r(t)}=(1+t)^{-\frac{1}{2}}$, then we can choose $t_{1}\geq0$ such that $r(t)< j_{0}$ for all $t\geq t_{1}$. Thus, for $t\geq t_{1}$ we have
\[
\begin{split}
-\left\|u\right\|^{l}_{\dot{B}^{2}}&=-\sum_{j\leq j_{0}}2^{2j}\left\|\Delta_{j}u\right\|_{L^{2}}=-\sum_{j\leq r(t)}2^{2j}\left\|\Delta_{j}u\right\|_{L^{2}}-\sum_{r(t)<j\leq j_{0}}2^{2j}\left\|\Delta_{j}u\right\|_{L^{2}} \\
&\leq -(1+t)^{-1}\sum_{r(t)<j\leq j_{0}}\left\|\Delta_{j}u\right\|_{L^{2}} =-(1+t)^{-1}\left\|u\right\|^{l}_{\dot{B}^{0}}+(1+t)^{-1}\sum_{j\leq r(t)}\left\|\Delta_{j}u\right\|_{L^{2}}.
\end{split}
\]
Together with (\ref{case1 further decay of u}) and (\ref{case1 further decay of tau}), (\ref{case1 Fourier split ineq2}) arrives at 
\eqn \label{case1 further decay estimate}
\begin{split}
&\frac{d}{dt}\left(\left\|u\right\|^{l}_{\dot{B}^{0}}+\left\|\tau\right\|^{l}_{\dot{B}^{1}}\right)+(1+t)^{-1}\left\|u\right\|^{l}_{\dot{B}^{0}}+\left\|\tau\right\|^{l}_{\dot{B}^{1}} \\
&\lesssim (1+t)^{-1}\sum_{j\leq r(t)}\left\|\Delta_{j}u\right\|_{L^{2}}+\left\|u\right\|_{\dot{B}^{0}}\left\|u\right\|_{\dot{B}^{\frac{5}{2}}}+(1+t)^{-\frac{9}{4}}.
\end{split}
\een
To obtain decay of $\sum_{j\leq r(t)}\left\|\Delta_{j}u\right\|_{L^{2}}$, we use Duhamel's formula for (\ref{case1 Fourier split ineq1}):
\[
f_{j}(t)\leq e^{-2^{2j}t}f_{j}(0)+\int^{t}_{0}e^{-2^{2j}(t-t')}X_{j}(t')\, dt',
\]
which implies
\[
\begin{split}
\sum_{j\leq r(t)}\left\|\Delta_{j}u\right\|_{L^{2}} &\lesssim \sum_{j\leq r(t)}f_{j}\lesssim \sum_{j\leq r(t)}2^{\frac{3}{2}j}2^{-\frac{3}{2}j}f_{j}(0)+\sum_{j\leq r(t)}\int^{t}_{0}2^{\frac{5}{2}j}2^{-\frac{5}{2}j}X_{j}(t')\, dt' \\
&\lesssim (1+t)^{-\frac{3}{4}}\left(\left\|u_{0}\right\|^{l}_{\dot{B}^{-\frac{3}{2}}_{2,\infty}}+\left\|\tau_{0}\right\|^{l}_{\dot{B}^{-\frac{1}{2}}_{2,\infty}}\right)+(1+t)^{-\frac{5}{4}}\int^{t}_{0}\sup_{j\leq j_{0}}2^{-\frac{5}{2}j}X_{j}(t')\,dt'.
\end{split}
\]
Since 
\[
\sup_{j\leq j_{0}}2^{-\frac{5}{2}j}X_{j}(t) \lesssim \left\|u\otimes u(t)\right\|_{L^{1}}+\left\|u\cdot\nabla\tau(t)\right\|_{L^{1}}+\left\|\mathcal{Q}(\tau,\nabla u)(t)\right\|_{L^{1}}\lesssim \left\|u(t)\right\|^{2}_{L^{2}}+(1+t)^{-1},
\]
we get 
\[
\sum_{j\leq r(t)}\left\|\Delta_{j}u(t)\right\|_{L^{2}} \lesssim (1+t)^{-\frac{1}{4}}.
\]
From (\ref{case1 further decay estimate}), we obtain
\[
\frac{d}{dt}\left[(1+t)\left(\left\|u(t)\right\|^{l}_{\dot{B}^{0}}+\left\|\tau(t)\right\|^{l}_{\dot{B}^{1}}\right)\right]\lesssim (1+t)^{-\frac{1}{4}},
\]
and equivalently
\[
\left\|u(t)\right\|^{l}_{\dot{B}^{0}}+\left\|\tau(t)\right\|^{l}_{\dot{B}^{1}}\lesssim (1+t)^{-\frac{1}{4}}.
\]
Repeating the above process again with $\left\|u(t)\right\|_{L^{2}}\lesssim \left\|u(t)\right\|_{\dot{B}^{0}}\lesssim (1+t)^{-\frac{1}{4}}$ implies
\[
\left\|u(t)\right\|^{l}_{\dot{B}^{0}}+\left\|\tau(t)\right\|^{l}_{\dot{B}^{1}}\lesssim (1+t)^{-\frac{3}{4}}, 
\]
that is,
\[
\left\|u(t)\right\|_{\dot{B}^{0}}\lesssim (1+t)^{-\frac{3}{4}}.
\]
Since
\[
\left\|u(t)\right\|^{l}_{\dot{B}^{s}}\lesssim \left(\left\|u(t)\right\|^{l}_{\dot{B^{0}}}\right)^{\frac{2}{s+2}}\left(\left\|u(t)\right\|^{l}_{\dot{B}^{s+2}}\right)^{\frac{s}{s+2}}\lesssim (1+t)^{-\frac{3}{4}\cdot\frac{2}{s+2}}\left(\left\|u(t)\right\|^{l}_{\dot{B}^{s+2}}\right)^{\frac{s}{s+2}}
\]
from (\ref{interpolation}), (\ref{eq:3.18}) becomes
\[
\frac{d}{dt}\left(\left\|u(t)\right\|_{\dot{B}^{s,\frac{5}{2}}}+\left\|\tau(t)\right\|_{\dot{B}^{s+1,\frac{5}{2}}}\right)+(1+t)^{\frac{3}{2s}}\left(\left\|u(t)\right\|^{l}_{\dot{B}^{s}}\right)^{1+\frac{2}{s}}+\left\|u(t)\right\|^{h}_{\dot{B}^{\frac{5}{2}}}+\left\|\tau(t)\right\|_{\dot{B}^{s+1,\frac{5}{2}}}\leq 0.
\]
Moreover, from 
\[
\left\|u(t)\right\|^{h}_{\dot{B}^{\frac{5}{2}}}+\left\|\tau(t)\right\|_{\dot{B}^{s+1,\frac{5}{2}}}\lesssim (1+t)^{-\frac{3}{4}},
\]
we have
\[
\frac{d}{dt}\left(\left\|u(t)\right\|_{\dot{B}^{s,\frac{5}{2}}}+\left\|\tau(t)\right\|_{\dot{B}^{s+1,\frac{5}{2}}}\right)+(1+t)^{\frac{3}{2s}}\left(\left\|u(t)\right\|_{\dot{B}^{s,\frac{5}{2}}}+\left\|\tau(t)\right\|_{\dot{B}^{s+1,\frac{5}{2}}}\right)^{1+\frac{2}{s}}\leq 0.
\]
By solving this inequality we get
\[
\left\|u(t)\right\|_{\dot{B}^{s,\frac{5}{2}}}+\left\|\tau(t)\right\|_{\dot{B}^{s+1,\frac{5}{2}}}\lesssim (1+t)^{-\frac{3}{4}-\frac{s}{2}},\quad s\in [\frac{1}{2},\frac{5}{2}],
\]
and using (\ref{interpolation}) again, we obtain (\ref{case1 further decay2 of u}). From $\left\|(u,\tau)\right\|_{\mathcal{E}_{T}}\lesssim \epsilon$ and (\ref{case1 further decay2 of u}), we have 
\[
\frac{d}{dt}\left\|\tau(t)\right\|_{\dot{B}^{0}}+\frac{\alpha}{2}\left\|\tau(t)\right\|_{\dot{B}^{0}}\leq C\left\|u(t)\right\|_{\dot{B}^{1}}\lesssim (1+t)^{-\frac{5}{4}},
\]
and 
\[
\frac{d}{dt}\left\|\tau(t)\right\|_{\dot{B}^{\frac{3}{2}}}+\frac{\alpha}{2}\left\|\tau(t)\right\|_{\dot{B}^{\frac{3}{2}}}\leq C\left\|u(t)\right\|_{\dot{B}^{\frac{5}{2}}}\leq C(1+t)^{-2},
\]
which implies 
\[
\left\|\tau(t)\right\|_{\dot{B}^{0}}\lesssim (1+t)^{-\frac{5}{4}},\quad \left\|\tau(t)\right\|_{\dot{B}^{\frac{3}{2}}}\lesssim (1+t)^{-2}.
\]
Then, (\ref{case1 further decay2 of tau}) is derived by (\ref{interpolation}).

\subsection{Proof of Corollary \ref{case1 decay2}}
As in the proof of Theorem \ref{case1 decay}, if we do not take the time integrations when deriving (\ref{case1 low estimate}) and (\ref{case1 high estimate2}) with $v=u$, $F=0$, and $-G=\mathcal{Q}(\tau,\nabla u)\simeq \tau^{2}+\tau\nabla u$, we arrive at 
\eqn \label{case1 ineq for decay2}
\begin{split}
& \frac{d}{dt}\left(\left\|u\right\|_{\dot{B}^{s,\frac{d}{2}-1}}+\left\|\tau\right\|_{\dot{B}^{s+1,\frac{d}{2}}}\right)+\left\|u\right\|_{\dot{B}^{s+2,\frac{d}{2}+1}}+\left\|\tau\right\|_{\dot{B}^{s+1,\frac{d}{2}}} \\
&\lesssim \left\|\mathcal{Q}(\tau,\nabla u)\right\|_{\dot{B}^{s+1,\frac{d}{2}}}+\left\|\mathbb{B}(\nabla u,\nabla \tau)\right\|_{\dot{B}^{s,\frac{d}{2}-1}}+\left\|u\right\|_{\dot{B}^{\frac{d}{2}+1}}\left(\left\|u\right\|_{\dot{B}^{s,\frac{d}{2}-1}}+\left\|\tau\right\|_{\dot{B}^{\min{(s,\frac{d}{2})}+1,\frac{d}{2}}}\right).
\end{split}
\een
Since
\[
\begin{split}
\left\|\mathcal{Q}(\tau,\nabla u)\right\|_{\dot{B}^{s+1,\frac{d}{2}}}&\lesssim \left\|\tau\right\|_{\dot{B}^{\frac{d}{2}}}\left(\left\|u\right\|_{\dot{B}^{s+2,\frac{d}{2}+1}}+\left\|\tau\right\|_{\dot{B}^{s+1,\frac{d}{2}}}\right)+\left\|u\right\|_{\dot{B}^{\frac{d}{2}+1}}\left\|\tau\right\|_{\dot{B}^{s+1,\frac{d}{2}}},\\
\left\|\mathbb{B}(\nabla u,\nabla\tau)\right\|_{\dot{B}^{s,\frac{d}{2}-1}}&\lesssim \left\|\mathbb{B}(\nabla u,\nabla\tau)\right\|_{\dot{B}^{\min{(s,\frac{d}{2})},\frac{d}{2}-1}}\lesssim \left\|u\right\|_{\dot{B}^{\frac{d}{2}+1}}\left\|\tau\right\|_{\dot{B}^{\min{(s,\frac{d}{2})}+1,\frac{d}{2}}}
\end{split}
\]
for $s>-\frac{d}{2}$, (\ref{case1 ineq for decay2}) can be rewritten as 
\eqn \label{eq:3.20}
\begin{split}
& \frac{d}{dt}\left(\left\|u\right\|_{\dot{B}^{s,\frac{d}{2}-1}}+\left\|\tau\right\|_{\dot{B}^{s+1,\frac{d}{2}}}\right)+\left\|u\right\|_{\dot{B}^{s+2,\frac{d}{2}+1}}+\left\|\tau\right\|_{\dot{B}^{s+1,\frac{d}{2}}} \\
&\lesssim \left\|\tau\right\|_{\dot{B}^{\frac{d}{2}}}\left(\left\|u\right\|_{\dot{B}^{s+2,\frac{d}{2}+1}}+\left\|\tau\right\|_{\dot{B}^{s+1,\frac{d}{2}}}\right) +\left\|u\right\|_{\dot{B}^{\frac{d}{2}+1}}\left(\left\|u\right\|_{\dot{B}^{s,\frac{d}{2}-1}}+\left\|\tau\right\|_{\dot{B}^{\min{(s,\frac{d}{2})}+1,\frac{d}{2}}}\right) \\
&\lesssim \left\|\tau\right\|_{\dot{B}^{\frac{d}{2}}}\left(\left\|u\right\|_{\dot{B}^{s+2,\frac{d}{2}+1}}+\left\|\tau\right\|_{\dot{B}^{s+1,\frac{d}{2}}}\right)+\left\|u\right\|_{\dot{B}^{\frac{d}{2}+1}}\left\|u\right\|^{l}_{\dot{B}^{s}}+\left\|u\right\|_{\dot{B}^{\frac{d}{2}+1}}\left\|\tau\right\|^{l}_{\dot{B}^{\frac{d}{2}+1}} \\
&+\left\|u\right\|_{\dot{B}^{\frac{d}{2}+1}}\left\|u\right\|^{h}_{\dot{B}^{\frac{d}{2}-1}} +\left\|u\right\|_{\dot{B}^{\frac{d}{2}+1}}\left\|\tau\right\|_{\dot{B}^{s+1,\frac{d}{2}}}.
\end{split}
\een
By (\ref{bound for decay 1}), (\ref{bound for decay 2}) and $\left\|(u,\tau)\right\|_{\mathcal{F}_{T}}\lesssim\epsilon$, the first three terms of the right-hand side of (\ref{eq:3.20}) can be absorbed into the left-hand side of (\ref{eq:3.20}) for $\frac{d}{2}-1\leq s\leq\frac{d}{2}+1$. We also bound the remaining terms: 
\[
\begin{split}
\left\|u\right\|_{\dot{B}^{\frac{d}{2}+1}}\left\|u\right\|^{h}_{\dot{B}^{\frac{d}{2}-1}}&=\left\|u\right\|^{l}_{\dot{B}^{\frac{d}{2}+1}}\left\|u\right\|^{h}_{\dot{B}^{\frac{d}{2}-1}}+\left\|u\right\|^{h}_{\dot{B}^{\frac{d}{2}+1}}\left\|u\right\|^{h}_{\dot{B}^{\frac{d}{2}-1}}\lesssim \left\|u\right\|_{\dot{B}^{\frac{d}{2}-1}}\left\|u\right\|^{h}_{\dot{B}^{\frac{d}{2}+1}},\\
\left\|u\right\|_{\dot{B}^{\frac{d}{2}+1}}\left\|\tau\right\|_{\dot{B}^{s+1,\frac{d}{2}}}  &=\left\|u\right\|^{l}_{\dot{B}^{\frac{d}{2}+1}}\left\|\tau\right\|_{\dot{B}^{s+1,\frac{d}{2}}}+\left\|u\right\|^{h}_{\dot{B}^{\frac{d}{2}+1}}\left\|\tau\right\|_{\dot{B}^{s+1,\frac{d}{2}}}  \\
&\lesssim \left\|u\right\|^{l}_{\dot{B}^{\frac{d}{2}-1}}\left\|\tau\right\|_{\dot{B}^{s+1,\frac{d}{2}}}+\left\|\tau\right\|_{\dot{B}^{\frac{d}{2}}}\left\|u\right\|^{h}_{\dot{B}^{\frac{d}{2}+1}}
\end{split}
\]
for $s\geq\frac{d}{2}-1$. Thus, we obtain  
\[
\frac{d}{dt}\left(\left\|u\right\|_{\dot{B}^{s,\frac{d}{2}-1}}+\left\|\tau\right\|_{\dot{B}^{s+1,\frac{d}{2}}}\right)+\left\|u\right\|_{\dot{B}^{s+2,\frac{d}{2}+1}}+\left\|\tau\right\|_{\dot{B}^{s+1,\frac{d}{2}}}\leq0, \quad \frac{d}{2}-1\leq s\leq\frac{d}{2}+1.
\]
As in the proof of Theorem \ref{case1 decay}, letting $s=\frac{d}{2}-1+s_{0}$, $s_{0}\in(0,2]$, we get
\[
\frac{d}{dt}\left(\left\|u\right\|_{\dot{B}^{\frac{d}{2}-1+s_{0},\frac{d}{2}-1}}+\left\|\tau\right\|_{\dot{B}^{\frac{d}{2}+s_{0},\frac{d}{2}}}\right)+\left(\left\|u\right\|_{\dot{B}^{\frac{d}{2}-1+s_{0},\frac{d}{2}-1}}+\left\|\tau\right\|_{\dot{B}^{\frac{d}{2}+s_{0},\frac{d}{2}}}\right)^{1+\frac{2}{s_{0}}}\leq0,
\]
which implies Corollary \ref{case1 decay2}.

\section{Case II: $(\nu_{1}, \nu_{2},\alpha)=(+,0,0)$} \label{case2}
In this section, we deal with (\ref{our model}) with  $(\nu_{1}, \nu_{2},\alpha)=(+,0,0)$: 
\begin{subequations} \label{our model Case II}
\begin{align} 
& u_{t}+u\cdot\nabla u -\nu_{1}\Delta u+\nabla p=\dv \tau, \label{our model Case II a}\\
& \tau_{t}+u\cdot\nabla \tau +\mathcal{Q}(\tau, \nabla u)=D(u), \label{our model Case II b}\\
& \dv u=0,
\end{align}
\end{subequations}
where $\mathcal{Q}(\tau,\nabla u)\simeq\tau\nabla u+(\nabla u)^{2}$. As  in Section \ref{case1}, we set the norms for initial data and solutions: 
\[
\begin{split}
&\left\|(u_{0},\tau_{0})\right\|_{\mathcal{E}_{0}}=\left\|u_{0}\right\|_{\dot{B}^{\frac{d}{2}-1,\frac{d}{2}+1}}+\left\|\tau_{0}\right\|_{\dot{B}^{\frac{d}{2}-1,\frac{d}{2}+1}}, \quad \|(u, \tau)\|_{\mathcal{E}_{T}}=\|(u, \tau)\|_{\mathcal{L}_{T}}+\|(u, \tau)\|_{\mathcal{H}_{T}}\\
&\left\|(u,\tau)\right\|_{\mathcal{L}_{T}}=\left\|u\right\|_{L^{\infty}_{T}\dot{B}^{\frac{d}{2}-1,\frac{d}{2}+1}}+\left\|\tau\right\|_{L^{\infty}_{T}\dot{B}^{\frac{d}{2}-1,\frac{d}{2}+1}},  \left\|(u,\tau)\right\|_{\mathcal{H}_{T}}=\left\|u\right\|_{L^{1}_{T}\dot{B}^{\frac{d}{2}+1,\frac{d}{2}+2}}+\left\|\Lambda^{-1}\mathbb{P}\dv\tau\right\|_{L^{1}_{T}\dot{B}^{\frac{d}{2}+1}}.
\end{split}
\]

\begin{theorem} \label{case2 theorem}\upshape
There exists a constant $\epsilon$ such that if $(u_{0},\tau_{0})\in \mathcal{E}_{0}$ with $\left\|(u_{0},\tau_{0})\right\|_{\mathcal{E}_{0}}\leq \epsilon$, there exists a unique solution $(u,\tau)$ in $\mathcal{E}_{T}$ of (\ref{our model Case II}) such that $\left\|(u,\tau)\right\|_{\mathcal{E}_{T}}\lesssim \left\|(u_{0},\tau_{0})\right\|_{\mathcal{E}_{0}}$ for all $T>0$.
\end{theorem}

To prove Theorem \ref{case2 theorem}, we use the following lemma dealing with the linearized system of (\ref{our model Case II}) with a divergence-free vector field $v$:
\eqn \label{case2 linear eq}
\begin{split} 
& u_{t}+\mathbb{P}(v\cdot\nabla u)-\nu_{1}\Delta u-\mathbb{P}\dv\tau=F, \\
& \tau_{t}+v\cdot\nabla\tau-D(u)=G,\\
&(\mathbb{P}\dv\tau)_{t}+\mathbb{P}\dv(v\cdot\nabla\tau)-\frac{1}{2}\Delta u=\mathbb{P}\dv G.
\end{split}
\een

\begin{lemma} \label{case2 lemma}\upshape
Let $-\frac{d}{2}<s\leq \frac{d}{2}+1$,  $(u_{0}, \tau_{0})\in \dot{B}^{s,\frac{d}{2}+1}$, and $v\in L^{1}_{T}\dot{B}^{\frac{d}{2}+1}$. Then, a solution $(u,\tau)$ of (\ref{case2 linear eq}) satisfies  the following bound
\begin{align*}
&\left\|u\right\|_{L^{\infty}_{T}\dot{B}^{s,\frac{d}{2}+1}}+\left\|\tau\right\|_{L^{\infty}_{T}\dot{B}^{s,\frac{d}{2}+1}}+\left\|u\right\|_{L^{1}_{T}\dot{B}^{s+2,\frac{d}{2}+2}}+\left\|\Lambda^{-1}\mathbb{P}\dv\tau\right\|_{L^{1}_{T}\dot{B}^{s+2,\frac{d}{2}+1}} \lesssim \left\|(u_{0},\tau_{0})\right\|_{\dot{B}^{s,\frac{d}{2}+1}}\\
& +\left\|(F, G)\right\|_{L^{1}_{T}\dot{B}^{s,\frac{d}{2}+1}}+\left\|\mathbb{B}(\nabla v,\nabla \tau)\right\|_{L^{1}_{T}\dot{B}^{s-1,\frac{d}{2}}} +\Big(\left\|u\right\|_{L^{\infty}_{T}\dot{B}^{s,\frac{d}{2}+1}}+\left\|\tau\right\|_{L^{\infty}_{T}\dot{B}^{s,\frac{d}{2}+1}}\Big)\left\|v\right\|_{L^{1}_{T}\dot{B}^{\frac{d}{2}+1}}.
\end{align*}
\end{lemma}

Let $s=\frac{d}{2}-1$, $v=u$, $F=0$, and $-G=\mathcal{Q}(\tau,\nabla u)\simeq \tau\nabla u+(\nabla u)^{2}$. By (\ref{product estimate}) and (\ref{decompose divergence free}),    
\[
\begin{split}
\left\|\mathcal{Q}(\tau,\nabla u)\right\|_{L^{1}_{T}\dot{B}^{\frac{d}{2}-1,\frac{d}{2}+1}}&\lesssim \left\|(u,\tau)\right\|_{\mathcal{L}_{T}}\left\|(u,\tau)\right\|_{\mathcal{H}_{T}},\\
\left\|\mathbb{B}(\nabla u,\nabla \tau)\right\|_{L^{1}_{T}\dot{B}^{\frac{d}{2}-2,\frac{d}{2}}} &\lesssim \left\|\tau\nabla u\right\|_{L^{1}_{T}\dot{B}^{\frac{d}{2}-1}}+\left\|\nabla u\nabla\tau\right\|_{L^{1}_{T}\dot{B}^{\frac{d}{2}}}\lesssim \left\|\tau\right\|_{L^{\infty}_{T}\dot{B}^{\frac{d}{2}-1,\frac{d}{2}+1}}\left\|u\right\|_{L^{1}_{T}\dot{B}^{\frac{d}{2}+1}}
\end{split}
\]
and so Lemma \ref{case2 lemma} yields $\left\|(u,\tau)\right\|_{\mathcal{E}_{T}}\lesssim \left\|(u_{0},\tau_{0})\right\|_{\mathcal{E}_{0}}$ for all $T>0$ when $\left\|(u_{0},\tau_{0})\right\|_{\mathcal{E}_{0}}$ is sufficiently small.

\begin{remark}\upshape \label{case2 decay remark}
In Case II, we have the  dissipative effect   to the equation of $\mathbb{P}\dv\tau$.  As in (\ref{case1 ineq for decay}), if we do not take time integrations in the proof of (\ref{case2 low estimate for decay}) and (\ref{case2 high estimate for decay}) with $v=u$, $F=0$ and $-G=\mathcal{Q}(\tau,\nabla u)\simeq \tau\nabla u+(\nabla u)^{2}$, we can obtain the following
\begin{align*}
& \frac{d}{dt}\left(\left\|u\right\|_{\dot{B}^{s,\frac{d}{2}+1}}+\left\|\Lambda^{-1}\mathbb{P}\dv \tau\right\|_{\dot{B}^{s,\frac{d}{2}+1}}\right)+\left\|u\right\|_{\dot{B}^{s+2,\frac{d}{2}+1}}+\left\|\Lambda^{-1}\mathbb{P}\dv \tau\right\|_{\dot{B}^{s+2,\frac{d}{2}+1}} \\
& \lesssim \left\|\mathcal{Q}(\tau,\nabla u)\right\|_{\dot{B}^{s,\frac{d}{2}+1}}+\left\|\mathbb{B}(\nabla u,\nabla\tau)\right\|_{\dot{B}^{s-1,\frac{d}{2}}} +\left\|u\right\|_{\dot{B}^{\frac{d}{2}+1}}\left(\left\|u\right\|_{\dot{B}^{s,\frac{d}{2}+1}}+\left\|\Lambda^{-1}\mathbb{P}\dv \tau\right\|_{\dot{B}^{s,\frac{d}{2}+1}}\right).
\end{align*}
However in this form, the left-hand side is not enough to control the right-hand side and hence we cannot derive the decay rates of $\mathbb{P}\dv\tau$.
\end{remark}

\begin{remark}\upshape
Using the bounds in the proof of Lemma \ref{case2 lemma}, we can deduce the following bound
\eqn \label{case2 uniqueness estimate} 
\begin{split}
&\left\|u\right\|_{L^{\infty}_{T}\dot{B}^{\frac{d}{2}-1}}+\left\|\tau\right\|_{L^{\infty}_{T}\dot{B}^{\frac{d}{2}-1,\frac{d}{2}}}+\left\|u\right\|_{L^{1}_{T}\dot{B}^{\frac{d}{2}+1}}+\left\|\Lambda^{-1}\mathbb{P}\dv\tau\right\|_{L^{1}_{T}\dot{B}^{\frac{d}{2}+1,\frac{d}{2}}} \lesssim \left\|u_{0}\right\|_{\dot{B}^{\frac{d}{2}-1}}+\left\|\tau_{0}\right\|_{\dot{B}^{\frac{d}{2}-1,\frac{d}{2}}}\\
& +\left\|F\right\|_{L^{1}_{T}\dot{B}^{\frac{d}{2}-1}} +\left\|G\right\|_{L^{1}_{T}\dot{B}^{\frac{d}{2}-1,\frac{d}{2}}} +\Big(\left\|u\right\|_{L^{\infty}_{T}\dot{B}^{\frac{d}{2}-1}}+\left\|\tau\right\|_{L^{\infty}_{T}\dot{B}^{\frac{d}{2}-1,\frac{d}{2}}}\Big)\left\|v\right\|_{L^{1}_{T}\dot{B}^{\frac{d}{2}+1}}.
\end{split}
\een
Let $v=u$, $F=0$ and $-G=\mathcal{Q}(\tau,\nabla u)\simeq \tau\nabla u$ and
\eqn \label{case2 uniqueness space}
\mathcal{F}_{0}=\dot{B}^{\frac{d}{2}-1}\times\dot{B}^{\frac{d}{2}-1,\frac{d}{2}},\qquad \mathcal{F}_{T}=L^{\infty}_{T}\dot{B}^{\frac{d}{2}-1}\cap L^{1}_{T}\dot{B}^{\frac{d}{2}+1}\times L^{\infty}_{T}\dot{B}^{\frac{d}{2}-1,\frac{d}{2}}.
\een
Then, $\left\|(u,\tau)\right\|_{\mathcal{F}_{T}} \lesssim \left\|(u_{0},\tau_{0})\right\|_{\mathcal{F}_{0}}$ for all $T>0$ when $\left\|(u_{0},\tau_{0})\right\|_{\mathcal{F}_{0}}$ is sufficiently small.  So, we can reprove the results in \cite{Chen-Hao}.
\end{remark}

\subsection{Proof of Lemma \ref{case2 lemma}}
\subsubsection{\bf Low frequency part}
From (\ref{case2 linear eq}), we obtain
\eqn \label{case2 eq1}
\frac{1}{2}\frac{d}{dt}\left\|\Delta_{j}u\right\|^{2}_{L^{2}}+\nu_{1}\left\|\Delta_{j}\Lambda u\right\|^{2}_{L^{2}}-(\Delta_{j}u,\Delta_{j}\dv \tau)=(\Delta_{j}F,\Delta_{j}u)-(\Delta_{j}(v\cdot\nabla u),\Delta_{j}u),
\een
\eqn \label{case2 eq2}
\begin{split}
\frac{1}{2}\frac{d}{dt}\left\|\Delta_{j}\Lambda^{-1}\mathbb{P}\dv \tau\right\|^{2}_{L^{2}}+\frac{1}{2}(\Delta_{j}u,\Delta_{j}\dv \tau)&=(\Delta_{j}\Lambda^{-1}\mathbb{P}\dv G,\Delta_{j}\Lambda^{-1}\mathbb{P}\dv \tau) \\
& -(\Delta_{j}\Lambda^{-1}\mathbb{P}\dv(v\cdot\nabla\tau),\Delta_{j}\Lambda^{-1}\mathbb{P}\dv \tau),
\end{split}
\een
and 
\eqn \label{case2 eq3}
\begin{split}
&\frac{d}{dt}(\Delta_{j}u,\Delta_{j}\mathbb{P}\dv \tau)-\left\|\Delta_{j}\mathbb{P}\dv \tau\right\|^{2}_{L^{2}}+\frac{1}{2}\left\|\Delta_{j}\Lambda u\right\|^{2}_{L^{2}}+\nu_{1}(\Delta_{j}\Lambda^{2}u,\Delta_{j}\mathbb{P}\dv \tau) \\
&=(\Delta_{j}F,\Delta_{j}\mathbb{P}\dv \tau)+(\Delta_{j}\mathbb{P}\dv G,\Delta_{j}u)-(\Delta_{j}(v\cdot\nabla u),\Delta_{j}\mathbb{P}\dv \tau)-(\Delta_{j}\mathbb{P}\dv(v\cdot\nabla\tau),\Delta_{j}u).
\end{split}
\een
Computing $\frac{1}{2}(\ref{case2 eq1})+(\ref{case2 eq2})-K_{1}(\ref{case2 eq3})$, we have
\[
\frac{1}{2}\frac{d}{dt}f^{2}_{j}+h^{2}_{j} =\widetilde{F}^{1}_{j}-\widetilde{F}^{2}_{j}-\widetilde{F}^{3}_{j}-\widetilde{F}^{4}_{j} \qquad\text{for}\ j \leq j_{0},
\]
where 
\[
\begin{split}
f^{2}_{j}&=\frac{1}{2}\left\|\Delta_{j}u\right\|^{2}_{L^{2}}+\left\|\Delta_{j}\Lambda^{-1}\mathbb{P}\dv \tau\right\|^{2}_{2}-2K_{1}(\Delta_{j}u,\Delta_{j}\mathbb{P}\dv \tau), \\
h^{2}_{j}&=\frac{\nu_{1}-K_{1}}{2}\left\|\Delta_{j}\Lambda u\right\|^{2}_{L^{2}}+K_{1}\left\|\Delta_{j}\mathbb{P}\dv \tau\right\|^{2}_{L^{2}}-\nu_{1}K_{1}(\Delta_{j}\Lambda^{2}u,\Delta_{j}\mathbb{P}\dv \tau),\\
\widetilde{F}^{1}_{j}&=\frac{1}{2}(\Delta_{j}F,\Delta_{j}u)+(\Delta_{j}\Lambda^{-1}\mathbb{P}\dv G,\Delta_{j}\Lambda^{-1}\mathbb{P}\dv \tau)-K_{1}(\Delta_{j}F,\Delta_{j}\mathbb{P}\dv \tau)-K_{1}(\Delta_{j}\mathbb{P}\dv G,\Delta_{j}u), \\
\widetilde{F}^{2}_{j}&=\frac{1}{2}(\Delta_{j}(v\cdot\nabla u),\Delta_{j}u)-K_{1}\left[(\Delta_{j}(v\cdot\nabla u),\Delta_{j}\mathbb{P}\dv \tau)+(\Delta_{j}(v\cdot\nabla\mathbb{P}\dv \tau),\Delta_{j}u)\right], \\
\widetilde{F}^{3}_{j}&=(\Delta_{j}\Lambda^{-1}(v\cdot\nabla\mathbb{P}\dv \tau),\Delta_{j}\Lambda^{-1}\mathbb{P}\dv \tau), \\
\widetilde{F}^{4}_{j}&=-K_{1}(\Delta_{j}\mathbb{B}(\nabla v,\nabla\tau),\Delta_{j}u)+(\Delta_{j}\Lambda^{-1}\mathbb{B}(\nabla v,\nabla \tau),\Delta_{j}\Lambda^{-1}\mathbb{P}\dv \tau).
\end{split}
\]
We apply Proposition \ref{decompose prop} again to $\mathbb{P}\dv(v\cdot\nabla\tau)$. For a sufficiently small $K_{1}>0$ and for each $j_{0}$, we note that  
\[
f^{2}_{j}\simeq \left\|\Delta_{j}u\right\|^{2}_{L^{2}}+\left\|\Delta_{j}\Lambda^{-1}\mathbb{P}\dv \tau\right\|^{2}_{L^{2}}, \quad h^{2}_{j}\simeq \left\|\Delta_{j}\Lambda u\right\|^{2}_{L^{2}}+\left\|\Delta_{j}\mathbb{P}\dv \tau\right\|^{2}_{L^{2}}\simeq 2^{2j}f^{2}_{j}.
\]
Also using (\ref{hybrid Besov embedding}), (\ref{convection term estimate}) and (\ref{convection term estimate 2}), we bound the nonlinear terms as  
\begin{align*}
& \left|\widetilde{F}^{1}_{j}\right|+\left|\widetilde{F}^{4}_{j}\right| \lesssim c_{j}2^{-sj}\left(\left\|(F, G)\right\|^{l}_{\dot{B}^{s}}+\left\|\mathbb{B}(\nabla v,\nabla \tau)\right\|^{l}_{\dot{B}^{s-1}}\right)f_{j}, \\
& \left|\widetilde{F}^{2}_{j}\right| \lesssim c_{j}2^{-sj}\left\|v\right\|_{\dot{B}^{\frac{d}{2}+1}}\left(\left\|u\right\|_{\dot{B}^{s,\frac{d}{2}+1}}+\left\|\Lambda^{-1}\mathbb{P}\dv \tau\right\|_{\dot{B}^{s,\frac{d}{2}+1}}\right) f_{j}, \\
& \left|\widetilde{F}^{3}_{j}\right| \lesssim c_{j}2^{-sj}\left\|v\right\|_{\dot{B}^{\frac{d}{2}+1}}\left\|\Lambda^{-1}\mathbb{P}\dv \tau\right\|_{\dot{B}^{s,\frac{d}{2}+1}} f_{j},
\end{align*}
whenever $j\leq j_{0}$ and $-\frac{d}{2}<s\leq \frac{d}{2}+1$. From these bounds, we obtain
\eqn \label{case2 low estimate for decay}
\begin{split}
&\left\|u\right\|^{l}_{L^{\infty}_{T}\dot{B}^{s}}+\left\|\Lambda^{-1}\mathbb{P}\dv \tau\right\|^{l}_{L^{\infty}_{T}\dot{B}^{s}}+\left\|u\right\|^{l}_{L^{1}_{T}\dot{B}^{s+2}}+\left\|\Lambda^{-1}\mathbb{P}\dv \tau\right\|^{l}_{L^{1}_{T}\dot{B}^{s+2}}  \\
&\lesssim \left\|u_{0}\right\|^{l}_{\dot{B}^{s}}+\left\|\Lambda^{-1}\mathbb{P}\dv \tau_{0}\right\|^{l}_{\dot{B}^{s}}+\left\|(F, G)\right\|^{l}_{L^{1}_{T}\dot{B}^{s}}+\left\|\mathbb{B}(\nabla v,\nabla \tau)\right\|^{l}_{L^{1}_{T}\dot{B}^{s-1}} \\
&+\Big(\left\|u\right\|_{L^{\infty}_{T}\dot{B}^{s,\frac{d}{2}+1}}+\left\|\Lambda^{-1}\mathbb{P}\dv \tau\right\|_{L^{\infty}_{T}\dot{B}^{s,\frac{d}{2}+1}}\Big)\left\|v\right\|_{L^{1}_{T}\dot{B}^{\frac{d}{2}+1}}.
\end{split}
\een
From (\ref{case2 linear eq}), we also have   
\eqn  \label{case2 eqs}
\begin{split}
& \frac{1}{2}\frac{d}{dt}\left\|\Delta_{j}u\right\|^{2}_{L^{2}}+\nu_{1}\left\|\Delta_{j}\Lambda u\right\|^{2}_{L^{2}}-(\Delta_{j}u,\Delta_{j}\mathbb{P}\dv \tau)=(\Delta_{j}F,\Delta_{j}u)-(\Delta_{j}(v\cdot\nabla u),\Delta_{j}u), \\
& \frac{1}{2}\frac{d}{dt}\left\|\Delta_{j}\tau\right\|^{2}_{L^{2}}+(\Delta_{j}u,\Delta_{j}\mathbb{P}\dv \tau)=(\Delta_{j}G,\Delta_{j}\tau)-(\Delta_{j}(v\cdot\nabla\tau),\Delta_{j}\tau)
\end{split}
\een
and so we obtain
\begin{align*}
& \frac{1}{2}\frac{d}{dt}\left(\left\|\Delta_{j}u\right\|^{2}_{L^{2}}+\left\|\Delta_{j}\tau\right\|^{2}_{L^{2}}\right)+\nu_{1}\left\|\Delta_{j}\Lambda u\right\|^{2}_{L^{2}} \\
& \lesssim c_{j}2^{-sj}\Big[\left\|(F,G)\right\|^{l}_{\dot{B}^{s}}+\left\|v\right\|_{\dot{B}^{\frac{d}{2}+1}}\left(\left\|u\right\|_{\dot{B}^{s,\frac{d}{2}+1}}+\left\|\tau\right\|_{\dot{B}^{s,\frac{d}{2}+1}}\right)\Big]\left(\left\|\Delta_{j}u\right\|_{L^{2}}+\left\|\Delta_{j}\tau\right\|_{L^{2}}\right).
\end{align*}
Thus, we have 
\[
\begin{split}
\left\|u\right\|^{l}_{L^{\infty}_{T}\dot{B}^{s}}+\left\|\tau\right\|^{l}_{L^{\infty}_{T}\dot{B}^{s}} &\lesssim \left\|(u_{0},\tau_{0})\right\|^{l}_{\dot{B}^{s}}+\left\|(F,G)\right\|^{l}_{L^{1}_{T}\dot{B}^{s}} +\Big(\left\|u\right\|_{L^{\infty}_{T}\dot{B}^{s,\frac{d}{2}+1}}+\left\|\tau\right\|_{L^{\infty}_{T}\dot{B}^{s,\frac{d}{2}+1}}\Big)\left\|v\right\|_{L^{1}_{T}\dot{B}^{\frac{d}{2}+1}}.
\end{split}
\]
By combining this with (\ref{case2 low estimate for decay}), we obtain
\eqn \label{case2 low estimate}
\begin{split}
&\left\|u\right\|^{l}_{L^{\infty}_{T}\dot{B}^{s}}+\left\|\tau\right\|^{l}_{L^{\infty}_{T}\dot{B}^{s}}+\left\|u\right\|^{l}_{L^{1}_{T}\dot{B}^{s+2}}+\left\|\Lambda^{-1}\mathbb{P}\dv \tau\right\|^{l}_{L^{1}_{T}\dot{B}^{s+2}}  \lesssim \left\|(u_{0},\tau_{0})\right\|^{l}_{\dot{B}^{s}}\\
&+\left\|(F,G)\right\|^{l}_{L^{1}_{T}\dot{B}^{s}}+\left\|\mathbb{B}(\nabla v,\nabla \tau)\right\|^{l}_{L^{1}_{T}\dot{B}^{s-1}} +\Big(\left\|u\right\|_{L^{\infty}_{T}\dot{B}^{s,\frac{d}{2}+1}}+\left\|\tau\right\|_{L^{\infty}_{T}\dot{B}^{s,\frac{d}{2}+1}}\Big)\left\|v\right\|_{L^{1}_{T}\dot{B}^{\frac{d}{2}+1}}.
\end{split}
\een

\subsubsection{\bf High frequency part}
From  (\ref{case2 linear eq}), we obtain
\eqn \label{case2 eq4}
\frac{1}{2}\frac{d}{dt}\left\|\Delta_{j}\Lambda u\right\|^{2}_{L^{2}}+\nu_{1}\left\|\Delta_{j}\Lambda^{2}u\right\|^{2}_{L^{2}}-(\Delta_{j}\Lambda u,\Delta_{j}\Lambda\dv \tau)=(\Delta_{j}\Lambda F,\Delta_{j}\Lambda u)-(\Delta_{j}\Lambda(v\cdot\nabla u),\Delta_{j}\Lambda u) \quad 
\een
and
\eqn \label{case2 eq5}
\frac{1}{2}\frac{d}{dt}\left\|\Delta_{j}\mathbb{P}\dv \tau\right\|^{2}_{L^{2}}+\frac{1}{2}(\Delta_{j}\Lambda u,\Delta_{
j}\Lambda\dv \tau)=(\Delta_{j}\mathbb{P}\dv G,\Delta_{j}\mathbb{P}\dv \tau)-(\Delta_{j}\mathbb{P}\dv(v\cdot\nabla\tau),\Delta_{j}\mathbb{P}\dv \tau). \quad 
\een
By $(\ref{case2 eq4})+(2+2\nu_{1}K_{2})(\ref{case2 eq5})-K_{2}(\ref{case2 eq3})$, we have
\[
\frac{1}{2}\frac{d}{dt}f^{2}_{j}+h^{2}_{j}=\widetilde{G}^{1}_{j}-\widetilde{G}^{2}_{j}-\widetilde{G}^{3}_{j}, \qquad\text{for}\ j\geq j_{0}+1,
\]
where 
\[
\begin{split}
&f^{2}_{j}=\left\|\Delta_{j}\Lambda u\right\|^{2}_{L^{2}}+(2+2\nu_{1}K_{2})\left\|\Delta_{j}\mathbb{P}\dv \tau\right\|^{2}_{2}-2K_{2}(\Delta_{j}u,\Delta_{j}\mathbb{P}\dv \tau), \\
&h^{2}_{j}=\nu_{1}\left\|\Delta_{j}\Lambda^{2}u\right\|^{2}_{L^{2}}-\frac{K_{2}}{2}\left\|\Delta_{j}\Lambda u\right\|^{2}_{L^{2}}+K_{2}\left\|\Delta_{j}\mathbb{P}\dv \tau\right\|^{2}_{L^{2}},\\
&\widetilde{G}^{1}_{j}=(\Delta_{j}\Lambda F,\Delta_{j}\Lambda u)+(2+2\nu_{1}K_{2})(\Delta_{j}\mathbb{P}\dv G,\Delta_{j}\mathbb{P}\dv \tau) -K_{2}(\Delta_{j}F,\Delta_{j}\mathbb{P}\dv \tau)-K_{2}(\Delta_{j}\mathbb{P}\dv G,\Delta_{j}u), \\
&\widetilde{G}^{2}_{j}=(\Delta_{j}\Lambda(v\cdot\nabla u),\Delta_{j}\Lambda u)+(2+2\nu_{1}K_{2})(\Delta_{j}(v\cdot\nabla\mathbb{P}\dv \tau),\Delta_{j}\mathbb{P}\dv \tau)\\
 &\quad\quad -K_{2}\left[(\Delta_{j}(v\cdot\nabla u),\Delta_{j}\mathbb{P}\dv \tau)+(\Delta_{j}(v\cdot\nabla\mathbb{P}\dv \tau),\Delta_{j}u)\right], \\
&\widetilde{G}^{3}_{j}=(2+2\nu_{1}K_{2})(\Delta_{j}\mathbb{B}(\nabla v,\nabla\tau),\Delta_{j}\mathbb{P}\dv \tau)-K_{2}(\Delta_{j}\mathbb{B}(\nabla v,\nabla\tau),\Delta_{j}u).
\end{split}
\]
For a sufficiently small $K_{2}>0$ and for each $j_{0}$, we have 
\[
f^{2}_{j}\simeq\left\|\Delta_{j}\Lambda u\right\|^{2}_{L^{2}}+\left\|\Delta_{j}\mathbb{P}\dv \tau\right\|^{2}_{L^{2}}, \quad h^{2}_{j}\simeq \left\|\Delta_{j}\Lambda^{2}u\right\|^{2}_{L^{2}}+\left\|\Delta_{j}\mathbb{P}\dv \tau\right\|^{2}_{L^{2}} \gtrsim f^{2}_{j}.
\]
Also using (\ref{hybrid Besov embedding}), (\ref{convection term estimate}) and (\ref{convection term estimate 2}), we bound the nonlinear terms as 
\begin{align*}
\left|\widetilde{G}^{1}_{j}\right|+\left|\widetilde{G}^{3}_{j}\right| &\lesssim c_{j}2^{-\frac{d}{2}j}\Big(\left\|(F,G)\right\|^{h}_{\dot{B}^{\frac{d}{2}+1}}+\left\|\mathbb{B}(\nabla v,\nabla\tau)\right\|^{h}_{\dot{B}^{\frac{d}{2}}}\Big) f_{j}, \\
\left|\widetilde{G}^{2}_{j}\right| &\lesssim c_{j}2^{-\frac{d}{2}j}\left\|v\right\|_{\dot{B}^{\frac{d}{2}+1}}\Big(\left\|u\right\|_{\dot{B}^{s,\frac{d}{2}+1}}+\left\|\mathbb{P}\dv \tau\right\|_{\dot{B}^{s,\frac{d}{2}}}\Big)f_{j}
\end{align*}
whenever $j\geq j_{0}+1$. From these bounds, we derive 
\eqn \label{case2 high estimate for decay}
\begin{split}
& \left\|u\right\|^{h}_{L^{\infty}_{T}\dot{B}^{\frac{d}{2}+1}}+\left\|\mathbb{P}\dv \tau\right\|^{h}_{L^{\infty}_{T}\dot{B}^{\frac{d}{2}}}+\left\|u\right\|^{h}_{L^{1}_{T}\dot{B}^{\frac{d}{2}+1}}+\left\|\mathbb{P}\dv \tau\right\|^{h}_{L^{1}_{T}\dot{B}^{\frac{d}{2}}} \\
& \lesssim \left\|u_{0}\right\|^{h}_{\dot{B}^{\frac{d}{2}+1}}+\left\|\mathbb{P}\dv \tau_{0}\right\|^{h}_{\dot{B}^{\frac{d}{2}}}+\left\|(F,G)\right\|^{h}_{L^{1}_{T}\dot{B}^{\frac{d}{2}+1}}+\left\|\mathbb{B}(\nabla v,\nabla\tau)\right\|^{h}_{L^{1}_{T}\dot{B}^{\frac{d}{2}}} \\
& +\Big(\left\|u\right\|_{L^{\infty}_{T}\dot{B}^{s,\frac{d}{2}+1}}+\left\|\mathbb{P}\dv \tau\right\|_{L^{\infty}_{T}\dot{B}^{s,\frac{d}{2}}}\Big)\left\|v\right\|_{L^{1}_{T}\dot{B}^{\frac{d}{2}+1}}.
\end{split}
\een
From (\ref{case2 eqs}), we also obtain 
\begin{align*}
& \left\|u\right\|^{h}_{L^{\infty}_{T}\dot{B}^{\frac{d}{2}}}+\left\|u\right\|^{h}_{L^{1}_{T}\dot{B}^{\frac{d}{2}+2}}\lesssim \left\|u_{0}\right\|^{h}_{\dot{B}^{\frac{d}{2}}}+\left\|\mathbb{P}\dv \tau\right\|^{h}_{L^{1}_{T}\dot{B}^{\frac{d}{2}}}+\left\|F\right\|^{h}_{L^{1}_{T}\dot{B}^{\frac{d}{2}}}+\left\|u\right\|_{L^{\infty}_{T}\dot{B}^{s,\frac{d}{2}+1}}\left\|v\right\|_{L^{1}_{T}\dot{B}^{\frac{d}{2}+1}}, \\
& \left\|\tau\right\|^{h}_{L^{\infty}_{T}\dot{B}^{\frac{d}{2}+1}}\lesssim \left\|\tau_{0}\right\|^{h}_{\dot{B}^{\frac{d}{2}+1}}+\left\|u\right\|^{h}_{L^{1}_{T}\dot{B}^{\frac{d}{2}+2}}+\left\|G\right\|^{h}_{L^{1}_{T}\dot{B}^{\frac{d}{2}+1}}+\left\|\tau\right\|_{L^{\infty}_{T}\dot{B}^{s,\frac{d}{2}+1}}\left\|v\right\|_{L^{1}_{T}\dot{B}^{\frac{d}{2}+1}}.
\end{align*}
Combining these three inequalities, we arrive at 
\eqn \label{case2 high estimate}
\begin{split}
& \left\|u\right\|^{h}_{L^{\infty}_{T}\dot{B}^{\frac{d}{2}+1}}+\left\|\tau\right\|^{h}_{L^{\infty}_{T}\dot{B}^{\frac{d}{2}+1}}+\left\|u\right\|^{h}_{L^{1}_{T}\dot{B}^{\frac{d}{2}+2}}+\left\|\Lambda^{-1}\mathbb{P}\dv \tau\right\|^{h}_{L^{1}_{T}\dot{B}^{\frac{d}{2}+1}} \lesssim \left\|(u_{0}, \tau_{0})\right\|^{h}_{\dot{B}^{\frac{d}{2}+1}}\\
& +\left\|(F,G)\right\|^{h}_{L^{1}_{T}\dot{B}^{\frac{d}{2}+1}}+\left\|\mathbb{B}(\nabla v,\nabla\tau)\right\|^{h}_{L^{1}_{T}\dot{B}^{\frac{d}{2}}}  +\Big(\left\|u\right\|_{L^{\infty}_{T}\dot{B}^{s,\frac{d}{2}+1}}+\left\|\tau\right\|_{L^{\infty}_{T}\dot{B}^{s,\frac{d}{2}+1}}\Big)\left\|v\right\|_{L^{1}_{T}\dot{B}^{\frac{d}{2}+1}}.
\end{split}
\een
By  (\ref{case2 low estimate}) and (\ref{case2 high estimate}), we complete the proof of Lemma \ref{case2 lemma}.

\subsection{Uniqueness}
We now assume that $(u_{1},\tau_{1})$ and $(u_{2},\tau_{2})$ are two solutions of (\ref{our model Case II}) with the same initial data. Let $(u, \tau)=(u_{1}-u_{2},\tau_{1}-\tau_{2})$. Then, $(u, \tau)$ satisfies $\dv u_{2}=0$ and
\eqn \label{case2 difference eq}
\begin{split}
& u_{t}+\mathbb{P}(u_{2}\cdot\nabla u)-\nu_{1}\Delta u-\mathbb{P}\dv \tau=-\mathbb{P}(u\cdot\nabla u_{1})=\delta F,\\
& \tau_{t}+u_{2}\cdot\nabla\tau-D(u)=-u\cdot\nabla\tau_{1}-\left[\mathcal{Q}(\tau_{1},\nabla u_{1})-\mathcal{Q}(\tau_{2},\nabla u_{2})\right]=\delta G.
\end{split}
\een
We now show the uniqueness in $\mathcal{F}_{T}$ defined in (\ref{case2 uniqueness space}). By applying (\ref{case2 uniqueness estimate}) to (\ref{case2 difference eq}), we have
\[
\left\|(u,\tau)\right\|_{\mathcal{F}_{T}}\lesssim \left\|\delta F\right\|_{L^{1}_{T}\dot{B}^{\frac{d}{2}-1}}+\left\|\delta G\right\|_{L^{1}_{T}\dot{B}^{\frac{d}{2}-1,\frac{d}{2}}}+\left\|u_{2}\right\|_{L^{1}_{T}\dot{B}^{\frac{d}{2}+1}}\left\|(u,\tau)\right\|_{\mathcal{F}_{T}}.
\]
We now look at each of the nonlinear terms in $\mathcal{F}_{T}$. First of all, 
\[
\begin{split}
\left\|u\cdot\nabla u_{1}\right\|_{L^{1}_{T}\dot{B}^{\frac{d}{2}-1}} +\left\|u\cdot\nabla\tau_{1}\right\|_{L^{1}_{T}\dot{B}^{\frac{d}{2}-1,\frac{d}{2}}} &\lesssim \left\|u\right\|_{L^{\infty}_{T}\dot{B}^{\frac{d}{2}-1}}\left\|u_{1}\right\|_{L^{1}_{T}\dot{B}^{\frac{d}{2}+1}}+\left\|u\right\|_{L^{2}_{T}\dot{B}^{\frac{d}{2}}}\left\|\tau_{1}\right\|_{L^{\infty}_{T}\dot{B}^{\frac{d}{2},\frac{d}{2}+1}}T^{1/2} \\
&\lesssim \left(1+T^{1/2}\right)\left\|(u_{1},\tau_{1})\right\|_{\mathcal{E}_{T}}\left\|(u,\tau)\right\|_{\mathcal{F}_{T}}.
\end{split}
\]
Since $\mathcal{Q}(\tau_{1},\nabla u_{1})-\mathcal{Q}(\tau_{2},\nabla u_{2})\simeq \tau\nabla u_{1}+\tau_{2}\nabla u+(\nabla u_{1}+\nabla u_{2})\nabla u$, we get
\[
\left\|\mathcal{Q}(\tau_{1},\nabla u_{1})-\mathcal{Q}(\tau_{2},\nabla u_{2})\right\|_{L^{1}_{T}\dot{B}^{\frac{d}{2}-1,\frac{d}{2}}}\lesssim \left(\left\|(u_{1},\tau_{1})\right\|_{\mathcal{E}_{T}}+\left\|(u_{2},\tau_{2})\right\|_{\mathcal{E}_{T}}\right)\left\|(u,\tau)\right\|_{\mathcal{F}_{T}}.
\]
By all these bounds, we deduce that  
\[
\left\|(u, \tau)\right\|_{\mathcal{F}_{T}}\lesssim (1+T^{1/2})\Big(\left\|(u_{1},\tau_{1})\right\|_{\mathcal{E}_{T}}+\left\|(u_{2},\tau_{2})\right\|_{\mathcal{E}_{T}}\Big)\left\|(u, \tau)\right\|_{\mathcal{F}_{T}}.
\]
Since $\left\|(u_{1},\tau_{1})\right\|_{\mathcal{E}_{T}}+\left\|(u_{2},\tau_{2})\right\|_{\mathcal{E}_{T}}\lesssim\epsilon$, $(u,\tau)=0$ in $\mathcal{F}_{T}$ for some finite time $T>0$. By repeating this process, we can conclude the uniqueness of solutions globally-in-time.

\section{Case III: $(\nu_{1}, \nu_{2},\alpha)=(0,+,+)$} \label{case3}
In this section, we deal with (\ref{our model}) with $(\nu_{1}, \nu_{2},\alpha)=(0,+,+)$:
\begin{subequations} \label{our model Case III}
\begin{align} 
& u_{t}+u\cdot\nabla u +\nabla p=\dv \tau, \label{our model a}\\
& \tau_{t}+u\cdot\nabla \tau+\alpha\tau-\nu_{2}\Delta\tau+\mathcal{Q}(\tau, \nabla u)=D(u), \label{our model b}\\
& \dv u=0
\end{align}
\end{subequations}
with $\mathcal{Q}(\tau,\nabla u)\simeq\tau^{2}+\tau\nabla u+(\nabla u)^{2}$. Let 
\[
\begin{split}
&\left\|(u_{0},\tau_{0})\right\|_{\mathcal{E}_{0}}=\left\|u_{0}\right\|_{\dot{B}^{\frac{d}{2}-1,\frac{d}{2}+1}}+\left\|\tau_{0}\right\|_{\dot{B}^{\frac{d}{2}}}, \quad \left\|(u,\tau)\right\|_{\mathcal{E}_{T}}=\left\|(u,\tau)\right\|_{\mathcal{L}_{T}}+\left\|(u,\tau)\right\|_{\mathcal{H}_{T}}\\
&\left\|(u,\tau)\right\|_{\mathcal{L}_{T}}=\left\|u\right\|_{L^{\infty}_{T}\dot{B}^{\frac{d}{2}-1,\frac{d}{2}+1}}+\left\|\tau\right\|_{L^{\infty}_{T}\dot{B}^{\frac{d}{2}}}, \quad \left\|(u,\tau)\right\|_{\mathcal{H}_{T}}=\left\|u\right\|_{L^{1}_{T}\dot{B}^{\frac{d}{2}+1}}+\left\|\tau\right\|_{L^{1}_{T}\dot{B}^{\frac{d}{2},\frac{d}{2}+2}}.
\end{split}
\]

\begin{theorem} \label{case3 theorem}\upshape
There exists  $\epsilon>0$ such that if $(u_{0},\tau_{0})\in \mathcal{E}_{0}$ with $\left\|(u_{0},\tau_{0})\right\|_{\mathcal{E}_{0}}\leq \epsilon$, there exists a unique solution $(u,\tau)$ in $\mathcal{E}_{T}$ of (\ref{our model Case III}) such that $\left\|(u,\tau)\right\|_{\mathcal{E}_{T}} \lesssim \left\|(u_{0},\tau_{0})\right\|_{\mathcal{E}_{0}}$ for all $T>0$.
\end{theorem}

\begin{theorem} \label{case3 decay} \upshape
The solution of Theorem \ref{case3 theorem} has the following decay rates:
\[
\left\|u(t)\right\|_{\dot{B}^{\frac{d}{2}-1+s_{0},\frac{d}{2}+1}}+\left\|\tau(t)\right\|_{\dot{B}^{\frac{d}{2}+s_{0},\frac{d}{2}}}\leq C(\left\|(u_{0},\tau_{0})\right\|_{\mathcal{E}_{0}}, s_{0})(1+t)^{-\frac{s_{0}}{2}}, \quad s_{0}\in(0,2].
\]
\end{theorem}

Similar to Case I, we also derive the further decay rates of the solution of Theorem \ref{case3 theorem}.

\begin{corollary} \label{case3 further decay}
Let $(u,\tau)\in\mathcal{E}_{T}$ be the solution of Theorem \ref{case3 theorem} with $(u_{0},\tau_{0})\in\mathcal{E}_{0}$.
\item{\textbullet} Assume that $(u_{0},\tau_{0})\in L^{2}$. Then, $(u,\tau)$ is in $L^{\infty}_{T}L^{2}\times L^{\infty}_{T}L^{2}\cap L^{2}_{T}H^{1}$ and satisfies
\begin{align*}
&\left\|u(t)\right\|_{\dot{B}^{s}}\lesssim (1+t)^{-\frac{s}{2}},\quad s\in(0,\frac{d}{2}+1], \\
\left\|\tau(t)\right\|_{L^{2}}\lesssim (1+t)^{-\frac{1}{2}},\qquad &\left\|\tau(t)\right\|_{\dot{B}^{s}}\lesssim (1+t)^{-\frac{1}{2}-\frac{s}{2}},\quad s\in (0,\frac{d}{2}].
\end{align*}
\item{\textbullet} Let $d=3$ and assume that $(u_{0},\tau_{0})\in \dot{B}^{-\frac{3}{2}}_{2,\infty}\times\dot{B}^{-\frac{1}{2}}_{2,\infty}$. Then, $(u,\tau)$ has the further decay rates:
\begin{align*}
&\left\|u(t)\right\|_{\dot{B}^{s}}\lesssim (1+t)^{-\frac{3}{4}-\frac{s}{2}},\quad s\in[0,\frac{5}{2}],\\
&\left\|\tau(t)\right\|_{\dot{B}^{s}}\lesssim (1+t)^{-\frac{5}{4}-\frac{s}{2}},\quad s\in[0,\frac{3}{2}].
\end{align*}
\end{corollary}

The proof of Corollary \ref{case3 further decay} is exactly same with that of Corollary \ref{case1 further decay} since the estimates of low frequency part are equivalent. Also, the comments in Remark \ref{case1 further decay remark} about the results in \cite{Huang} and present paper still hold in Case III.

As before, to prove Theorem \ref{case3 theorem} and Theorem \ref{case3 decay}, the following linearized equation with a divergence-free vector field $v$ is used
\eqn \label{case3 linear eq}
\begin{split} 
& u_{t}+\mathbb{P}(v\cdot\nabla u)-\mathbb{P}\dv \tau=F, \\
& \tau_{t}+v\cdot\nabla\tau+\alpha\tau-\nu_{2}\Delta\tau-D(u)=G,\\
& (\mathbb{P}\dv \tau)_{t}+\mathbb{P}\dv(v\cdot\nabla\tau)+\alpha\mathbb{P}\dv \tau-\nu_{2}\Delta(\mathbb{P}\dv \tau)-\frac{1}{2}\Delta u=\mathbb{P}\dv G.
\end{split}
\een

\begin{lemma} \label{case3 lemma}\upshape
Let $-\frac{d}{2}-1<s_{1}\leq \frac{d}{2}+1$, $-\frac{d}{2}<s_{2}\leq\frac{d}{2}$, $u_{0}\in \dot{B}^{s_{1},s_{2}+1}$, $\tau_{0}\in \dot{B}^{s_{1}+1,s_{2}}$ and $v\in L^{1}_{T}\dot{B}^{\frac{d}{2}+1}$. Then, a solution $(u,\tau)$ of (\ref{case3 linear eq}) satisfies the following bound
\begin{align*}
&\left\|u\right\|_{L^{\infty}_{T}\dot{B}^{s_{1},s_{2}+1}}+\left\|\tau\right\|_{L^{\infty}_{T}\dot{B}^{s_{1}+1,s_{2}}}+\left\|u\right\|_{L^{1}_{T}\dot{B}^{s_{1}+2,s_{2}+1}}+\left\|\tau\right\|_{L^{1}_{T}\dot{B}^{s_{1}+1,s_{2}+2}} \\
& \lesssim \left\|u_{0}\right\|_{\dot{B}^{s_{1},s_{2}+1}}+\left\|\tau_{0}\right\|_{\dot{B}^{s_{1}+1,s_{2}}}+\left\|F\right\|_{L^{1}_{T}\dot{B}^{s_{1},s_{2}+1}}+\left\|G\right\|_{L^{1}_{T}\dot{B}^{s_{1}+1,s_{2}}} +\left\|\mathbb{B}(\nabla v,\nabla \tau)\right\|_{L^{1}_{T}\dot{B}^{s_{1},s_{2}-1}} \\
& +\Big(\left\|u\right\|_{L^{\infty}_{T}\dot{B}^{s_{1},s_{2}+1}}+\left\|\tau\right\|_{L^{\infty}_{T}\dot{B}^{\min{(s_{1},\frac{d}{2})}+1,s_{2}}}\Big)\left\|v\right\|_{L^{1}_{T}\dot{B}^{\frac{d}{2}+1}}.
\end{align*}
\end{lemma}

Let $(s_{1},s_{2})=\left(\frac{d}{2}-1,\frac{d}{2}\right)$, $v=u$, $F=0$ and $-G=\mathcal{Q}(\tau, \nabla u)\simeq \tau^{2}+\tau\nabla u+(\nabla u)^{2}$. Since
\[
\left\|\mathcal{Q}(\tau,\nabla u)\right\|_{L^{1}_{T}\dot{B}^{\frac{d}{2}}}+\left\|\mathbb{B}(\nabla u,\nabla\tau)\right\|_{L^{1}_{T}\dot{B}^{\frac{d}{2}-1}}\lesssim \left\|(u,\tau)\right\|_{\mathcal{L}_{T}} \left\|(u,\tau)\right\|_{\mathcal{H}_{T}}, 
\]

Lemma \ref{case3 lemma} yields  the desired result $\left\|(u,\tau)\right\|_{\mathcal{E}_{T}}\lesssim \left\|(u_{0},\tau_{0})\right\|_{\mathcal{E}_{0}}$ for all $T>0$ when $\left\|(u_{0},\tau_{0})\right\|_{\mathcal{E}_{0}}$ is sufficiently small.

\subsection{Proof of Lemma \ref{case3 lemma}}
\subsubsection{\bf Low frequency part}
From  (\ref{case3 linear eq}), we obtain
\eqn \label{case3 eq1}
\frac{1}{2}\frac{d}{dt}\left\|\Delta_{j}u\right\|^{2}_{L^{2}}-(\Delta_{j}u,\Delta_{j}\dv \tau)=(\Delta_{j}F,\Delta_{j}u)-(\Delta_{j}(v\cdot\nabla u),\Delta_{j}u), 
\een
\eqn \label{case3 eq2}
\begin{split}
\frac{1}{2}\frac{d}{dt}\left\|\Delta_{j}\Lambda\tau\right\|^{2}_{L^{2}}+\alpha\left\|\Delta_{j}\Lambda\tau\right\|^{2}_{L^{2}}+\nu_{2}\left\|\Delta_{j}\Lambda^{2}\tau\right\|^{2}_{L^{2}}+(\Delta_{j}u,\Delta_{j}\Lambda^{2}\dv \tau) \\=(\Delta_{j}\Lambda G,\Delta_{j}\Lambda\tau)-(\Delta_{j}\Lambda(v\cdot\nabla\tau),\Delta_{j}\Lambda\tau),
\end{split}
\een
and 
\eqn \label{case3 eq3}
\begin{split}
&\frac{d}{dt}(\Delta_{j}u,\Delta_{j}\mathbb{P}\dv \tau)-\left\|\Delta_{j}\mathbb{P}\dv \tau\right\|^{2}_{L^{2}}+\frac{1}{2}\left\|\Delta_{j}\Lambda u\right\|^{2}_{L^{2}}+\alpha(\Delta_{j}u,\Delta_{j}\dv \tau)+\nu_{2}(\Delta_{j}u,\Delta_{j}\Lambda^{2}\dv \tau) \\
&=(\Delta_{j}F,\Delta_{j}\mathbb{P}\dv \tau)+(\Delta_{j}\mathbb{P}\dv G,\Delta_{j}u)-(\Delta_{j}(v\cdot\nabla u),\Delta_{j}\mathbb{P}\dv \tau)-(\Delta_{j}\mathbb{P}\dv(v\cdot\nabla\tau),\Delta_{j}u).
\end{split}
\een
Computing $\alpha(\ref{case3 eq1})+K_{1}(\ref{case3 eq2})+(\ref{case3 eq3})$, we have
\[
\frac{1}{2}\frac{d}{dt}f^{2}_{j}+h^{2}_{j} =\widetilde{F}^{1}_{j}-\widetilde{F}^{2}_{j}-\widetilde{F}^{3}_{j}-\widetilde{F}^{4}_{j} \qquad\text{for}\ j \leq j_{1},
\]
where 
\[
\begin{split}
&f^{2}_{j}=\alpha\left\|\Delta_{j}u\right\|^{2}_{L^{2}}+K_{1}\left\|\Delta_{j}\Lambda\tau\right\|^{2}_{2}+2(\Delta_{j}u,\Delta_{j}\mathbb{P}\dv \tau), \\
&h^{2}_{j}=\alpha K_{1}\left\|\Delta_{j}\Lambda\tau\right\|^{2}_{L^{2}}+\nu_{2}K_{1}\left\|\Delta_{j}\Lambda^{2}\tau\right\|^{2}_{L^{2}}-\left\|\Delta_{j}\mathbb{P}\dv \tau\right\|^{2}_{L^{2}}+\frac{1}{2}\left\|\Delta_{j}\Lambda u\right\|^{2}_{L^{2}} +(K_{1}+\nu_{2})(\Delta_{j}u,\Delta_{j}\Lambda^{2}\dv \tau),\\
&\widetilde{F}^{1}_{j}=\alpha(\Delta_{j}F,\Delta_{j}u)+K_{1}(\Delta_{j}\Lambda G,\Delta_{j}\Lambda\tau)+(\Delta_{j}F,\Delta_{j}\mathbb{P}\dv \tau)+(\Delta_{j}\mathbb{P}\dv G,\Delta_{j}u), \\
&\widetilde{F}^{2}_{j}=\alpha(\Delta_{j}(v\cdot\nabla u),\Delta_{j}u)+(\Delta_{j}(v\cdot\nabla u),\Delta_{j}\mathbb{P}\dv \tau)+(\Delta_{j}(v\cdot\nabla\mathbb{P}\dv \tau),\Delta_{j}u), \\
&\widetilde{F}^{3}_{j}= K_{1}(\Delta_{j}\Lambda(v\cdot\nabla\tau),\Delta_{j}\Lambda\tau), \\
& \widetilde{F}^{4}_{j}= (\Delta_{j}\mathbb{B}(\nabla v,\nabla\tau),\Delta_{j}u).
\end{split}
\]
By choosing $K_{1}=\frac{2C^{2}_{0}}{\alpha}$ and sufficiently small $j_{1}$, we obtain 
\[
f^{2}_{j}\simeq \left\|\Delta_{j}u\right\|^{2}_{L^{2}}+\left\|\Delta_{j}\Lambda\tau\right\|^{2}_{L^{2}}, \quad h^{2}_{j}\simeq \left\|\Delta_{j}\Lambda u\right\|^{2}_{L^{2}}+\left\|\Delta_{j}\Lambda\tau\right\|^{2}_{L^{2}} \gtrsim 2^{2j}f^{2}_{j}.
\]
Using (\ref{convection term estimate}) and (\ref{convection term estimate 2}), we also bound the nonlinear terms  
\begin{align*}
&\left|\widetilde{F}^{1}_{j}\right|+\left|\widetilde{F}^{4}_{j}\right| \lesssim c_{j}2^{-s_{1}j}\left(\left\|F\right\|^{l}_{\dot{B}^{s_{1}}}+\left\|G\right\|^{l}_{\dot{B}^{s_{1}+1}}+\left\|\mathbb{B}(\nabla v,\nabla \tau)\right\|^{l}_{\dot{B}^{s_{1}}}\right)f_{j}, \\
&\left|\widetilde{F}^{2}_{j}\right| \lesssim c_{j}2^{-s_{1}j}\left\|v\right\|_{\dot{B}^{\frac{d}{2}+1}}\left(\left\|u\right\|_{\dot{B}^{s_{1},s_{2}+1}}+\left\|\tau\right\|_{\dot{B}^{s_{1}+1,s_{2}}}\right) f_{j},  \\
&\left|\widetilde{F}^{3}_{j}\right| \lesssim c_{j}2^{-\min{(s_{1},\frac{d}{2})}j}\left\|v\right\|_{\dot{B}^{\frac{d}{2}+1}}\left\|\tau\right\|_{\dot{B}^{\min{(s_{1},\frac{d}{2})}+1,s_{2}}}\left\|\Delta_{j}\Lambda\tau\right\|_{L^{2}} \lesssim c_{j}2^{-s_{1}j}\left\|v\right\|_{\dot{B}^{\frac{d}{2}+1}}\left\|\tau\right\|_{\dot{B}^{\min{(s_{1},\frac{d}{2})}+1,s_{2}}}f_{j}
\end{align*}
whenever $j\leq j_{1}$ and $-\frac{d}{2}-1<s_{1}\leq\frac{d}{2}+1$, $-\frac{d}{2}<s_{2}\leq \frac{d}{2}$. From these bounds, we deduce that 
\begin{align*}
&\left\|u\right\|^{l}_{L^{\infty}_{T}\dot{B}^{s_{1}}}+\left\|\tau\right\|^{l}_{L^{\infty}_{T}\dot{B}^{s_{1}+1}}+\left\|u\right\|^{l}_{L^{1}_{T}\dot{B}^{s_{1}+2}}+\left\|\tau\right\|^{l}_{L^{1}_{T}\dot{B}^{s_{1}+3}}  \\
&\lesssim \left\|u_{0}\right\|^{l}_{\dot{B}^{s_{1}}}+\left\|\tau_{0}\right\|^{l}_{\dot{B}^{s_{1}+1}}+\left\|F\right\|^{l}_{L^{1}_{T}\dot{B}^{s_{1}}}+\left\|G\right\|^{l}_{L^{1}_{T}\dot{B}^{s_{1}+1}}+\left\|\mathbb{B}(\nabla v,\nabla\tau)\right\|^{l}_{L^{1}_{T}\dot{B}^{s_{1}}} \\
&+\Big(\left\|u\right\|_{L^{\infty}_{T}\dot{B}^{s_{1},s_{2}+1}}+\left\|\tau\right\|_{L^{\infty}_{T}\dot{B}^{\min{(s_{1},\frac{d}{2})}+1,s_{2}}}\Big)\left\|v\right\|_{L^{1}_{T}\dot{B}^{\frac{d}{2}+1}}. 
\end{align*}
 From (\ref{case3 eq2}), we also have
\begin{align*}
\left\|\tau\right\|^{l}_{L^{\infty}_{T}\dot{B}^{s_{1}+1}}+\left\|\tau\right\|^{l}_{L^{1}_{T}\dot{B}^{s_{1}+1}} \lesssim \left\|\tau_{0}\right\|^{l}_{\dot{B}^{s_{1}+1}}+\left\|u\right\|^{l}_{L^{1}_{T}\dot{B}^{s_{1}+2}}+\left\|G\right\|^{l}_{L^{1}_{T}\dot{B}^{s_{1}+1}}+\left\|\tau\right\|_{L^{\infty}_{T}\dot{B}^{\min{(s_{1},\frac{d}{2})}+1,s_{2}}}\left\|v\right\|_{L^{1}_{T}\dot{B}^{\frac{d}{2}+1}}.
\end{align*}
Combining these two inequalities with (\ref{hybrid Besov embedding}), we arrive at 
\eqn \label{case3 low estimate}
\begin{split}
&\left\|u\right\|^{l}_{L^{\infty}_{T}\dot{B}^{s_{1}}}+\left\|\tau\right\|^{l}_{L^{\infty}_{T}\dot{B}^{s_{1}+1}}+\left\|u\right\|^{l}_{L^{1}_{T}\dot{B}^{s_{1}+2}}+\left\|\tau\right\|^{l}_{L^{1}_{T}\dot{B}^{s_{1}+1}} \lesssim \left\|u_{0}\right\|^{l}_{\dot{B}^{s_{1}}}+\left\|\tau_{0}\right\|^{l}_{\dot{B}^{s_{1}+1}} +\left\|F\right\|^{l}_{L^{1}_{T}\dot{B}^{s_{1}}}\\
&+\left\|G\right\|^{l}_{L^{1}_{T}\dot{B}^{s_{1}+1}}+\left\|\mathbb{B}(\nabla v,\nabla\tau)\right\|^{l}_{L^{1}_{T}\dot{B}^{s_{1}}} +\Big(\left\|u\right\|_{L^{\infty}_{T}\dot{B}^{s_{1},s_{2}+1}}+\left\|\tau\right\|_{L^{\infty}_{T}\dot{B}^{\min{(s_{1},\frac{d}{2})}+1,s_{2}}}\Big)\left\|v\right\|_{L^{1}_{T}\dot{B}^{\frac{d}{2}+1}}.
\end{split}
\een

\subsubsection{\bf High frequency part}
Using (\ref{case3 linear eq}), we obtain
\eqn \label{case3 eq4}
\frac{1}{2}\frac{d}{dt}\left\|\Delta_{j}\Lambda u\right\|^{2}_{L^{2}}-(\Delta_{j}u,\Delta_{j}\Lambda^{2}\dv \tau)=(\Delta_{j}\Lambda F,\Delta_{j}\Lambda u)-(\Delta_{j}\Lambda(v\cdot\nabla u),\Delta_{j}\Lambda u),
\een
and
\eqn \label{case3 eq5}
\begin{split}
\frac{1}{2}\frac{d}{dt}\left\|\Delta_{j}\tau\right\|^{2}_{L^{2}}+\alpha\left\|\Delta_{j}\tau\right\|^{2}_{L^{2}}+\nu_{2}\left\|\Delta_{j}\Lambda\tau\right\|^{2}_{L^{2}}+(\Delta_{j}u,\Delta_{j}\dv \tau) =(\Delta_{j}G,\Delta_{j}\tau)-(\Delta_{j}(v\cdot\nabla\tau),\Delta_{j}\tau).
\end{split}
\een
By taking $\nu_{2}(\ref{case3 eq4})+K_{2}(\ref{case3 eq5})+(\ref{case3 eq3})$, we have
\[
\frac{1}{2}\frac{d}{dt}f^{2}_{j}+h^{2}_{j}=\widetilde{G}^{1}_{j}-\widetilde{G}^{2}_{j}-\widetilde{G}^{3}_{j}, \qquad\text{for}\ j\geq j_{2}+1,
\]
where 
\begin{align*}
f^{2}_{j} &=\nu_{2}\left\|\Delta_{j}\Lambda u\right\|^{2}_{L^{2}}+K_{2}\left\|\Delta_{j}\tau\right\|^{2}_{L^{2}}+2(\Delta_{j}u,\Delta_{j}\mathbb{P}\dv \tau), \\
h^{2}_{j}&=\alpha K_{2}\left\|\Delta_{j}\tau\right\|^{2}_{L^{2}}+\nu_{2} K_{2}\left\|\Delta_{j}\Lambda\tau\right\|^{2}_{L^{2}}-\left\|\Delta_{j}\mathbb{P}\dv \tau\right\|^{2}_{L^{2}}+\frac{1}{2}\left\|\Delta_{j}\Lambda u\right\|^{2}_{L^{2}} +(K_{2}+\alpha)(\Delta_{j}u,\Delta_{j}\dv \tau),\\
\widetilde{G}^{1}_{j} &=\nu_{2}(\Delta_{j}\Lambda F,\Delta_{j}\Lambda u)+K_{2}(\Delta_{j}G,\Delta_{j}\tau)+(\Delta_{j}F,\Delta_{j}\mathbb{P}\dv \tau)+(\Delta_{j}\mathbb{P}\dv G,\Delta_{j}u), \\
\widetilde{G}^{2}_{j}&=\nu_{2}(\Delta_{j}\Lambda(v\cdot\nabla u),\Delta_{j}\Lambda u)+K_{2}(\Delta_{j}(v\cdot\nabla\tau),\Delta_{j}\tau) +(\Delta_{j}(v\cdot\nabla u),\Delta_{j}\mathbb{P}\dv \tau)\\
&+(\Delta_{j}(v\cdot\nabla\mathbb{P}\dv \tau),\Delta_{j}u), \\
\widetilde{G}^{3}_{j} &=(\Delta_{j}\mathbb{B}(\nabla v,\nabla\tau),\Delta_{j}u).
\end{align*}
By choosing $K_{2}=\frac{2C^{2}_{0}}{\nu_{2}}$ and sufficiently large $j_{2}$, we notice that  
\[
f^{2}_{j}\simeq\left\|\Delta_{j}\Lambda u\right\|^{2}_{L^{2}}+\left\|\Delta_{j}\tau\right\|^{2}_{L^{2}}, \quad h^{2}_{j}\simeq \left\|\Delta_{j}\Lambda u\right\|^{2}_{L^{2}}+\left\|\Delta_{j}\Lambda\tau\right\|^{2}_{L^{2}} \gtrsim f^{2}_{j}.
\]
Moreover, (\ref{convection term estimate}) and (\ref{convection term estimate 2}) are also used to bound the nonlinear terms as follows: 
\begin{align*}
\left|\widetilde{G}^{1}_{j}\right|+\left|\widetilde{G}^{3}_{j}\right| &\lesssim c_{j}2^{-s_{2}j}\Big(\left\|F\right\|^{h}_{\dot{B}^{s_{2}+1}}+\left\|G\right\|^{h}_{\dot{B}^{s_{2}}}+\left\|\mathbb{B}(\nabla v,\nabla \tau)\right\|^{h}_{\dot{B}^{s_{2}-1}}\Big)f_{j}, \\
\left|\widetilde{G}^{2}_{j}\right| &\lesssim c_{j}2^{-s_{2}j}\left\|v\right\|_{\dot{B}^{\frac{d}{2}+1}}\Big(\left\|u\right\|_{\dot{B}^{s_{1},s_{2}+1}}+\left\|\tau\right\|_{\dot{B}^{\frac{d}{2}+1,s_{2}}}\Big)f_{j}
\end{align*}
whenever $j\geq j_{2}+1$ and $-\frac{d}{2}-1<s_{1}\leq \frac{d}{2}+1$, $-\frac{d}{2}<s_{2}\leq \frac{d}{2}$. From these bounds, we obtain 
\begin{align*}
& \left\|u\right\|^{h}_{L^{\infty}_{T}\dot{B}^{s_{2}+1}}+\left\|\tau\right\|^{h}_{L^{\infty}_{T}\dot{B}^{s_{2}}}+\left\|u\right\|^{h}_{L^{1}_{T}\dot{B}^{s_{2}+1}}+\left\|\tau\right\|^{h}_{L^{1}_{T}\dot{B}^{s_{2}}} \lesssim \left\|u_{0}\right\|^{h}_{\dot{B}^{s_{2}+1}}+\left\|\tau_{0}\right\|^{h}_{\dot{B}^{s_{2}}} +\left\|F\right\|^{h}_{L^{1}_{T}\dot{B}^{s_{2}+1}}\\
& +\left\|G\right\|^{h}_{L^{1}_{T}\dot{B}^{s_{2}}}+\left\|\mathbb{B}(\nabla v,\nabla \tau)\right\|^{h}_{L^{1}_{T}\dot{B}^{s_{2}-1}} +\Big(\left\|u\right\|_{L^{\infty}_{T}\dot{B}^{s_{1},s_{2}+1}}+\left\|\tau\right\|_{L^{\infty}_{T}\dot{B}^{\frac{d}{2}+1,s_{2}}}\Big)\left\|v\right\|_{L^{1}_{T}\dot{B}^{\frac{d}{2}+1}}.
\end{align*}
From (\ref{case3 eq5}), we also bound $\tau$ as 
\begin{align*}
\left\|\tau\right\|^{h}_{L^{\infty}_{T}\dot{B}^{s_{2}}}+\left\|\tau\right\|^{h}_{L^{1}_{T}\dot{B}^{s_{2}+2}} &\lesssim \left\|\tau_{0}\right\|^{h}_{\dot{B}^{s_{2}}}+\left\|u\right\|^{h}_{L^{1}_{T}\dot{B}^{s_{2}+1}}+\left\|G\right\|^{h}_{L^{1}_{T}\dot{B}^{s_{2}}}+\left\|\tau\right\|_{L^{\infty}_{T}\dot{B}^{\frac{d}{2}+1,s_{2}}}\left\|v\right\|_{L^{1}_{T}\dot{B}^{\frac{d}{2}+1}}.
\end{align*}
Combining these two inequalities, we obtain 
\eqn \label{case3 high estimate}
\begin{split}
&\left\|u\right\|^{h}_{L^{\infty}_{T}\dot{B}^{s_{2}+1}}+\left\|\tau\right\|^{h}_{L^{\infty}_{T}\dot{B}^{s_{2}}}+\left\|u\right\|^{h}_{L^{1}_{T}\dot{B}^{s_{2}+1}}+\left\|\tau\right\|^{h}_{L^{1}_{T}\dot{B}^{s_{2}+2}} \lesssim \left\|u_{0}\right\|^{h}_{\dot{B}^{s_{2}+1}}+\left\|\tau_{0}\right\|^{h}_{\dot{B}^{s_{2}}}+\left\|F\right\|^{h}_{L^{1}_{T}\dot{B}^{s_{2}+1}}  \\
&+\left\|G\right\|^{h}_{L^{1}_{T}\dot{B}^{s_{2}}}+\left\|\mathbb{B}(\nabla v,\nabla \tau)\right\|^{h}_{L^{1}_{T}\dot{B}^{s_{2}-1}}+\Big(\left\|u\right\|_{L^{\infty}_{T}\dot{B}^{s_{1},s_{2}+1}}+\left\|\tau\right\|_{L^{\infty}_{T}\dot{B}^{\frac{d}{2}+1,s_{2}}}\Big)\left\|v\right\|_{L^{1}_{T}\dot{B}^{\frac{d}{2}+1}}.
\end{split}
\een

\subsubsection{\bf Intermediate frequency part}
Unlike the proof of Lemma \ref{case1 lemma}, there is a gap between low and high frequency parts because we choose sufficiently small $j_{1}$ and sufficiently large $j_{2}$ to obtain (\ref{case3 low estimate}) and (\ref{case3 high estimate}), respectively. We now fill up the gap, $j_{1}+1\leq j \leq j_{2}$. In fact,  this part is very easy to handle since $2^{j}\simeq 1$, but we provide some details to complete the proof. From $(\ref{case3 eq1})+(\ref{case3 eq5})+K(\ref{case3 eq3})$, we obtain
\[
\frac{1}{2}\frac{d}{dt}f^{2}_{j}+h^{2}_{j} =\widetilde{H}^{1}_{j}-\widetilde{H}^{2}_{j}-\widetilde{H}^{3}_{j} \qquad\text{for}\ j_{1}+1\leq j \leq j_{2},
\]
where 
\begin{align*}
f^{2}_{j}&=\left\|\Delta_{j}u\right\|^{2}_{L^{2}}+\left\|\Delta_{j}\tau\right\|^{2}_{L^{2}}+2K(\Delta_{j}u,\Delta_{j}\mathbb{P}\dv \tau), \\
h^{2}_{j}&=\alpha\left\|\Delta_{j}\tau\right\|^{2}_{L^{2}}+\nu_{2}\left\|\Delta_{j}\Lambda\tau\right\|^{2}_{L^{2}}-K\left\|\Delta_{j}\mathbb{P}\dv \tau\right\|^{2}_{L^{2}}+\frac{K}{2}\left\|\Delta_{j}\Lambda u\right\|^{2}_{L^{2}} +\alpha K(\Delta_{j}u,\Delta_{j}\dv \tau)\\
&+\nu_{2} K(\Delta_{j}u,\Delta_{j}\Lambda^{2}\dv \tau)\\
\widetilde{H}^{1}_{j}&=(\Delta_{j}F,\Delta_{j}u)+(\Delta_{j}G,\Delta_{j}\tau)+K(\Delta_{j}F,\Delta_{j}\mathbb{P}\dv \tau)+K(\Delta_{j}\mathbb{P}\dv G,\Delta_{j}u), \\
\widetilde{H}^{2}_{j}&=(\Delta_{j}(v\cdot\nabla u),\Delta_{j}u)+(\Delta_{j}(v\cdot\nabla\tau),\Delta_{j}\tau)+K\left[(\Delta_{j}(v\cdot\nabla u),\Delta_{j}\mathbb{P}\dv \tau)+(\Delta_{j}(v\cdot\nabla\mathbb{P}\dv \tau),\Delta_{j}u)\right], \\
\widetilde{H}^{3}_{j}&=K(\Delta_{j}\mathbb{B}(\nabla v,\nabla\tau),\Delta_{j}u).
\end{align*}
For sufficiently small $K>0$, we obtain
\[
f^{2}_{j}\simeq h^{2}_{j}\simeq \left\|\Delta_{j}u\right\|^{2}_{L^{2}}+\left\|\Delta_{j}\tau\right\|^{2}_{L^{2}},
\]
and by (\ref{convection term estimate}) and (\ref{convection term estimate 2}) with a slight abuse of the fact $2^{j}\simeq 1$, we have 
\begin{align*}
\left|\widetilde{H}^{1}_{j}\right|+\left|\widetilde{H}^{3}_{j}\right| &\lesssim \left(\left\|\Delta_{j}F\right\|_{L^{2}}+\left\|\Delta_{j}G\right\|_{L^{2}}+\left\|\Delta_{j}\mathbb{B}(\nabla v,\nabla \tau)\right\|_{L^{2}}\right)f_{j}, \\
\left|\widetilde{H}^{2}_{j}\right| &\lesssim c_{j}\left\|v\right\|_{\dot{B}^{\frac{d}{2}+1}}\left(\left\|u\right\|_{\dot{B}^{s_{1},s_{2}+1}}+\left\|\tau\right\|_{\dot{B}^{\frac{d}{2}+1,s_{2}}}\right)f_{j},
\end{align*}
wherever $j_{1}+1\leq j\leq j_{2}$. Thus, we obtain 
\eqn \label{case3 intermediate estimate}
\begin{split}
&\left\|u\right\|^{I}_{L^{\infty}_{T}\dot{B}^{s_{1}}}+\left\|\tau\right\|^{I}_{L^{\infty}_{T}\dot{B}^{s_{1}+1}}+\left\|u\right\|^{I}_{L^{1}_{T}\dot{B}^{s_{1}+2}}+\left\|\tau\right\|^{I}_{L^{1}_{T}\dot{B}^{s_{1}+1}} \lesssim \left\|u_{0}\right\|^{I}_{\dot{B}^{s_{1}}}+\left\|\tau_{0}\right\|^{I}_{\dot{B}^{s_{1}+1}}+\left\|F\right\|^{I}_{L^{1}_{T}\dot{B}^{s_{1}}} \\
&+\left\|G\right\|^{I}_{L^{1}_{T}\dot{B}^{s_{1}+1}}+\left\|\mathbb{B}(\nabla v,\nabla \tau)\right\|^{I}_{L^{1}_{T}\dot{B}^{s_{1}}}+\Big(\left\|u\right\|_{L^{\infty}_{T}\dot{B}^{s_{1},s_{2}+1}}+\left\|\tau\right\|_{L^{\infty}_{T}\dot{B}^{\frac{d}{2}+1,s_{2}}}\Big)\left\|v\right\|_{L^{1}_{T}\dot{B}^{\frac{d}{2}+1}},
\end{split}
\een
where 
\[
\left\|f\right\|^{I}_{\dot{B}^{s}}:=\sum^{j_{2}}_{j=j_{1}+1}2^{js}\left\|\Delta_{j}f\right\|_{L^{2}}\simeq \sum^{j_{2}}_{j=j_{1}+1}\left\|\Delta_{j}f\right\|_{L^{2}}.
\]
By (\ref{case3 low estimate}), (\ref{case3 high estimate}) and (\ref{case3 intermediate estimate}) with taking $j_{0}=j_{2}$, we complete the proof of Lemma \ref{case3 lemma}.

\subsection{Uniqueness} \label{sec:5.2}
If $(u_{1},\tau_{1})$ and $(u_{2},\tau_{2})$ are two solutions of (\ref{our model Case III}) with the same initial data, let $(u,\tau)= (u_{1}-u_{2},\tau_{1}-\tau_{2})$. Then, $(u,\tau)$ satisfies $\dv u_{2}=0$ and
\eqn \label{case3 difference eq}
\begin{split}
& u_{t}+\mathbb{P}(u_{2}\cdot\nabla u)-\mathbb{P}\dv \tau=-\mathbb{P}(u\cdot\nabla u_{1})=\delta F, \\
& \tau_{t}+u_{2}\cdot\nabla\tau+\alpha\tau-\nu_{2}\Delta\tau-D(u)=-u\cdot\nabla\tau_{1}-\left[\mathcal{Q}(\tau_{1},\nabla u_{1})-\mathcal{Q}(\tau_{2},\nabla u_{2})\right]=\delta G.
\end{split}
\een
We now show the uniqueness of solutions by estimating $(u,\tau)$ in 
\[
\mathcal{F}_{T}:=L^{\infty}_{T}\dot{B}^{\frac{d}{2}-1,\frac{d}{2}}\cap L^{1}_{T}\dot{B}^{\frac{d}{2}+1,\frac{d}{2}}\times L^{\infty}_{T}\dot{B}^{\frac{d}{2},\frac{d}{2}-1}\cap L^{1}_{T}\dot{B}^{\frac{d}{2},\frac{d}{2}+1}.
\]
Lemma \ref{case3 lemma} with $(s_{1},s_{2})=\left(\frac{d}{2}-1,\frac{d}{2}-1\right)$ yields
\eqn \label{case3 difference estimate}
\begin{split}
\left\|(u,\tau)\right\|_{\mathcal{F}_{T}} &\lesssim \left\|\delta F\right\|_{L^{1}_{T}\dot{B}^{\frac{d}{2}-1,\frac{d}{2}}}+\left\|\delta G\right\|_{L^{1}_{T}\dot{B}^{\frac{d}{2},\frac{d}{2}-1}}+\left\|\mathbb{B}(\nabla u_{2},\nabla\tau)\right\|_{L^{1}_{T}\dot{B}^{\frac{d}{2}-1,\frac{d}{2}-2}} \\
&+\Big(\left\|u\right\|_{L^{\infty}_{T}\dot{B}^{\frac{d}{2}-1,\frac{d}{2}}}+\left\|\tau\right\|_{L^{\infty}_{T}\dot{B}^{\frac{d}{2},\frac{d}{2}-1}}\Big)\left\|u_{2}\right\|_{L^{1}_{T}\dot{B}^{\frac{d}{2}+1}}.
\end{split}
\een
Using (\ref{hybrid Besov embedding}) and (\ref{product estimate}) repeatedly, we bound 
\begin{align*}
& \left\|u\cdot\nabla u_{1}\right\|_{L^{1}_{T}\dot{B}^{\frac{d}{2}-1,\frac{d}{2}}} \lesssim \left\|u\right\|_{L^{\infty}_{T}\dot{B}^{\frac{d}{2}-1,\frac{d}{2}}}\left\|u_{1}\right\|_{L^{1}_{T}\dot{B}^{\frac{d}{2}+1}}, \\
& \left\|u\cdot\nabla\tau_{1}\right\|_{L^{1}_{T}\dot{B}^{\frac{d}{2},\frac{d}{2}-1}} \lesssim \left\|u\right\|_{L^{\infty}_{T}\dot{B}^{\frac{d}{2},\frac{d}{2}-1}}\left\|\tau_{1}\right\|_{L^{1}_{T}\dot{B}^{\frac{d}{2}+1}}, \\
& \left\|\mathbb{B}(\nabla u_{2},\nabla\tau)\right\|_{L^{1}_{T}\dot{B}^{\frac{d}{2}-1,\frac{d}{2}-2}} \lesssim \left\|\mathbb{B}(\nabla u_{2},\nabla\tau)\right\|_{L^{1}_{T}\dot{B}^{\frac{d}{2}-1}} \lesssim \left\|u_{2}\right\|_{L^{\infty}_{T}\dot{B}^{\frac{d}{2}+1}}\left\|\tau\right\|_{L^{1}_{T}\dot{B}^{\frac{d}{2}}}.
\end{align*}
Since $\mathcal{Q}(\tau_{1},\nabla u_{1})-\mathcal{Q}(\tau_{2},\nabla u_{2}) \simeq (\tau_{1}+\tau_{2})\tau+\tau\nabla u_{1}+\tau_{2}\nabla u+(\nabla u_{1}+\nabla u_{2})\nabla u$, we have 
\begin{align*}
\left\|\mathcal{Q}(\tau_{1},\nabla u_{1})-\mathcal{Q}(\tau_{2},\nabla u_{2})\right\|_{L^{1}_{T}\dot{B}^{\frac{d}{2},\frac{d}{2}-1}} &\lesssim \Big(\left\|\tau_{1}\right\|_{L^{1}_{T}\dot{B}^{\frac{d}{2}}}+\left\|\tau_{2}\right\|_{L^{1}_{T}\dot{B}^{\frac{d}{2}}}+\left\|u_{1}\right\|_{L^{1}_{T}\dot{B}^{\frac{d}{2}+1}}\Big)\left\|\tau\right\|_{L^{\infty}_{T}\dot{B}^{\frac{d}{2},\frac{d}{2}-1}} \\
&+\Big(\left\|\tau_{2}\right\|_{L^{\infty}_{T}\dot{B}^{\frac{d}{2}}}+\left\|u_{1}\right\|_{L^{\infty}_{T}\dot{B}^{\frac{d}{2}+1}}+\left\|u_{2}\right\|_{L^{\infty}_{T}\dot{B}^{\frac{d}{2}+1}}\Big)\left\|u\right\|_{L^{1}_{T}\dot{B}^{\frac{d}{2}+1,\frac{d}{2}}}.
\end{align*}
Combining (\ref{case3 difference estimate}) and the bound of nonlinear terms, and using  (\ref{hybrid Besov embedding}), we deduce that 
\[
\left\|(u, \tau)\right\|_{\mathcal{F}_{T}}\lesssim \Big(\left\|(u_{1},\tau_{1})\right\|_{\mathcal{E}_{T}}+\left\|(u_{2},\tau_{2})\right\|_{\mathcal{E}_{T}}\Big)\left\|(u, \tau)\right\|_{\mathcal{F}_{T}}.
\]
Since $\left\|(u_{1},\tau_{1})\right\|_{\mathcal{E}_{T}}+\left\|(u_{2},\tau_{2})\right\|_{\mathcal{E}_{T}}\lesssim \epsilon$, $(u,\tau)=0$ in $\mathcal{F}_{T}$ which completes the uniqueness part.


\subsection{Proof of Theorem  \ref{case3 decay}}
A similar argument of Theorem \ref{case1 decay} can be used to prove Theorem  \ref{case3 decay}.
Indeed, if we do not take time integrations in the proof of Lemma \ref{case3 lemma} with $(s_{1},s_{2})=(s,\frac{d}{2})$ and $v=u$, $F=0$, $-G=\mathcal{Q}(\tau,\nabla u)\simeq\tau^{2}+\tau\nabla u+(\nabla u)^{2}$, we derive  
\eqn \label{case3 ineq for decay}
\begin{split}
& \frac{d}{dt}\Big(\left\|u\right\|_{\dot{B}^{s,\frac{d}{2}+1}}+\left\|\tau\right\|_{\dot{B}^{s+1,\frac{d}{2}}}\Big)+\left\|u\right\|_{\dot{B}^{s+2,\frac{d}{2}+1}}+\left\|\tau\right\|_{\dot{B}^{s+1,\frac{d}{2}+2}} \\
&\lesssim \left\|\mathcal{Q}(\tau,\nabla u)\right\|_{\dot{B}^{s+1,\frac{d}{2}}}+\left\|\mathbb{B}(\nabla u,\nabla \tau)\right\|_{\dot{B}^{s,\frac{d}{2}-1}}+\left\|u\right\|_{\dot{B}^{\frac{d}{2}+1}}\left(\left\|u\right\|_{\dot{B}^{s,\frac{d}{2}+1}}+\left\|\tau\right\|_{\dot{B}^{\min{(s,\frac{d}{2})}+1,\frac{d}{2}}}\right).
\end{split}
\een
From (\ref{hybrid Besov embedding}) and (\ref{product estimate}), we have 
\[
\begin{split}
\left\|\mathcal{Q}(\tau,\nabla u)\right\|_{\dot{B}^{s+1,\frac{d}{2}}}&\lesssim \Big(\left\|u\right\|_{\dot{B}^{\frac{d}{2}+1}}+\left\|\tau\right\|_{\dot{B}^{\frac{d}{2}}}\Big)\left(\left\|u\right\|_{\dot{B}^{s+2,\frac{d}{2}+1}}+\left\|\tau\right\|_{\dot{B}^{s+1,\frac{d}{2}}}\right),\\
\left\|\mathbb{B}(\nabla u,\nabla \tau)\right\|_{\dot{B}^{s,\frac{d}{2}-1}}&\lesssim \left\|\mathbb{B}(\nabla u,\nabla\tau)\right\|_{\dot{B}^{\min{(s,\frac{d}{2})},\frac{d}{2}-1}}\lesssim \left\|u\right\|_{\dot{B}^{\frac{d}{2}+1}}\left\|\tau\right\|_{\dot{B}^{\min{(s,\frac{d}{2})}+1,\frac{d}{2}}}
\end{split}
\]
for $s>-\frac{d}{2}$. Thus, using (\ref{hybrid Besov embedding}) we can rearrange (\ref{case3 ineq for decay})  to get
\begin{align*}
& \frac{d}{dt}\left(\left\|u\right\|_{\dot{B}^{s,\frac{d}{2}+1}}+\left\|\tau\right\|_{\dot{B}^{s+1,\frac{d}{2}}}\right)+\left\|u\right\|_{\dot{B}^{s+2,\frac{d}{2}+1}}+\left\|\tau\right\|_{\dot{B}^{s+1,\frac{d}{2}+2}} \\
&\lesssim \left\|\tau\right\|_{\dot{B}^{\frac{d}{2}}}\left(\left\|u\right\|_{\dot{B}^{s+2,\frac{d}{2}+1}}+\left\|\tau\right\|_{\dot{B}^{s+1,\frac{d}{2}}}\right)+\left\|u\right\|_{\dot{B}^{\frac{d}{2}+1}}\left(\left\|u\right\|_{\dot{B}^{s,\frac{d}{2}+1}}+\left\|\tau\right\|_{\dot{B}^{\min{(s,\frac{d}{2})}+1,\frac{d}{2}}}\right) \\
&\lesssim \Big(\left\|u\right\|_{\dot{B}^{\frac{d}{2}+1}}+\left\|\tau\right\|_{\dot{B}^{\frac{d}{2}}}\Big)\left(\left\|u\right\|_{\dot{B}^{s+2,\frac{d}{2}+1}}+\left\|\tau\right\|_{\dot{B}^{s+1,\frac{d}{2}}}\right)+\left\|u\right\|_{\dot{B}^{\frac{d}{2}+1}}\left\|u\right\|^{l}_{\dot{B}^{s}}+\left\|u\right\|_{\dot{B}^{\frac{d}{2}+1}}\left\|\tau\right\|^{l}_{\dot{B}^{\frac{d}{2}+1}}.
\end{align*}
Using (\ref{bound for decay 1}) and (\ref{bound for decay 2}), we obtain that for $\frac{d}{2}-1\leq s\leq \frac{d}{2}+1$
\begin{align*}
& \frac{d}{dt}\left(\left\|u\right\|_{\dot{B}^{s,\frac{d}{2}+1}}+\left\|\tau\right\|_{\dot{B}^{s+1,\frac{d}{2}}}\right)+\left\|u\right\|_{\dot{B}^{s+2,\frac{d}{2}+1}}+\left\|\tau\right\|_{\dot{B}^{s+1,\frac{d}{2}+2}} \\
&\lesssim \Big(\left\|u\right\|_{\dot{B}^{\frac{d}{2}-1,\frac{d}{2}+1}}+\left\|\tau\right\|_{\dot{B}^{\frac{d}{2}}}\Big)\left(\left\|u\right\|_{\dot{B}^{s+2,\frac{d}{2}+1}}+\left\|\tau\right\|_{\dot{B}^{s+1,\frac{d}{2}}}\right) \lesssim \left\|(u,\tau)\right\|_{\mathcal{E}}\left(\left\|u\right\|_{\dot{B}^{s+2,\frac{d}{2}+1}}+\left\|\tau\right\|_{\dot{B}^{s+1,\frac{d}{2}}}\right).
\end{align*}
Since $\left\|(u,\tau)\right\|_{\mathcal{E}}\lesssim\epsilon$, we obtain 
\[
\frac{d}{dt}\left(\left\|u\right\|_{\dot{B}^{s,\frac{d}{2}+1}}+\left\|\tau\right\|_{\dot{B}^{s+1,\frac{d}{2}}}\right)+\left\|u\right\|_{\dot{B}^{s+2,\frac{d}{2}+1}}+\left\|\tau\right\|_{\dot{B}^{s+1,\frac{d}{2}+2}}\leq 0.
\]
As in the proof of Theorem \ref{case1 decay}, using (\ref{interpolation}) and letting $s=\frac{d}{2}-1+s_{0}$, $s_{0}\in(0,2]$, we get
\[
\frac{d}{dt}\left(\left\|u\right\|_{\dot{B}^{\frac{d}{2}-1+s_{0},\frac{d}{2}+1}}+\left\|\tau\right\|_{\dot{B}^{\frac{d}{2}+s_{0},\frac{d}{2}}}\right)+C\left(\left\|u\right\|_{\dot{B}^{\frac{d}{2}-1+s_{0},\frac{d}{2}+1}}+\left\|\tau\right\|_{\dot{B}^{\frac{d}{2}+s_{0},\frac{d}{2}}}\right)^{1+\frac{2}{s_{0}}}\leq 0,
\]
which implies Theorem \ref{case3 decay}.

\section{Case IV: $(\nu_{1}, \nu_{2},\alpha)=(0,+,0)$} \label{case4}
In this section, we consider (\ref{our model}) with $(\nu_{1}, \nu_{2},\alpha)=(0,+,0)$:
\eqn \label{our model Case IV}
\begin{split} 
& u_{t}+u\cdot\nabla u +\nabla p=\dv \tau, \label{our model a}\\
& \tau_{t}+u\cdot\nabla \tau-\nu_{2}\Delta\tau+\mathcal{Q}(\tau, \nabla u)=D(u), \label{our model b}\\
& \dv u=0,
\end{split}
\een
where $\mathcal{Q}(\tau, \nabla u)\simeq \tau\nabla u+(\nabla u)^{2}$. Let 
\[
\begin{split}
&\|(u_{0}, \tau_{0})\|_{\mathcal{E}_{0}}=\left\|u_{0}\right\|_{\dot{B}^{\frac{d}{2}-1,\frac{d}{2}+1}}+\left\|\tau_{0}\right\|_{\dot{B}^{\frac{d}{2}-1,\frac{d}{2}}}, \quad \|(u, \tau)\|_{\mathcal{E}_{T}}=\|(u, \tau)\|_{\mathcal{L}_{T}}+\|(u, \tau)\|_{\mathcal{H}_{T}}\\
& \|(u, \tau)\|_{\mathcal{L}_{T}}=\left\|u\right\|_{L^{\infty}_{T}\dot{B}^{\frac{d}{2}-1,\frac{d}{2}+1}}+\left\|\tau\right\|_{L^{\infty}_{T}\dot{B}^{\frac{d}{2}-1,\frac{d}{2}}}, \quad \|(u, \tau)\|_{\mathcal{H}_{T}}=\left\|u\right\|_{L^{1}_{T}\dot{B}^{\frac{d}{2}+1}}+\left\|\tau\right\|_{L^{1}_{T}\dot{B}^{\frac{d}{2}+1,\frac{d}{2}+2}}.
\end{split}
\]

\begin{theorem} \label{case4 theorem}\upshape
There exists $\epsilon>0$ such that if $(u_{0},\tau_{0})\in \mathcal{E}_{0}$ with $\left\|(u_{0},\tau_{0})\right\|_{\mathcal{E}_{0}}\leq \epsilon$, there exists a unique solution $(u,\tau)$ in $\mathcal{E}_{T}$ of (\ref{our model Case IV}) such that $\left\|(u,\tau)\right\|_{\mathcal{E}_{T}} \lesssim \left\|(u_{0},\tau_{0})\right\|_{\mathcal{E}_{0}}$ for all $T>0$.
\end{theorem}

\begin{theorem} \label{case4 decay}\upshape
The solutions of Theorem \ref{case4 theorem} have the following decay rates:
\[
\left\|u(t)\right\|_{\dot{B}^{\frac{d}{2}-1+s_{0},\frac{d}{2}+1}}+\left\|\tau(t)\right\|_{\dot{B}^{\frac{d}{2}-1+s_{0},\frac{d}{2}}}\leq C(\left\|(u_{0},\tau_{0})\right\|_{\mathcal{E}_{0}},s_{0})(1+t)^{-\frac{s_{0}}{2}}, \quad s_{0}\in(0,2].
\]
\end{theorem}

We also derive the further decay rates of the solution of Theorem \ref{case4 theorem}.

\begin{corollary} \label{case4 further decay}
Let $(u,\tau)\in\mathcal{E}_{T}$ be the solution of Theorem \ref{case4 theorem} with $(u_{0},\tau_{0})\in\mathcal{E}_{0}$.
\item{\textbullet} Assume that $(u_{0},\tau_{0})\in L^{2}$. Then, $(u,\tau)$ is in $L^{\infty}_{T}L^{2}\times L^{\infty}_{T}L^{2}\cap L^{2}_{T}\dot{H}^{1}$ and satisfies
\eqn \label{case4 further decay1}
\begin{split}
&\left\|u(t)\right\|_{\dot{B}^{s}}\lesssim (1+t)^{-\frac{s}{2}},\quad s\in(0,\frac{d}{2}+1], \\
&\left\|\tau(t)\right\|_{\dot{B}^{s}}\lesssim (1+t)^{-\frac{s}{2}},\quad s\in (0,\frac{d}{2}]. 
\end{split}
\een
\item{\textbullet} Let $d=3$ and assume that $(u_{0},\tau_{0})\in \dot{B}^{-\frac{3}{2}}_{2,\infty}$. Then, $(u,\tau)$ has the further decay rates:
\eqn \label{case4 further decay2}
\begin{split}
&\left\|u(t)\right\|_{\dot{B}^{s}}\lesssim (1+t)^{-\frac{3}{4}-\frac{s}{2}},\quad s\in[0,\frac{5}{2}], \\
&\left\|\tau(t)\right\|_{\dot{B}^{s}}\lesssim (1+t)^{-\frac{3}{4}-\frac{s}{2}},\quad s\in[0,\frac{3}{2}]. 
\end{split}
\een
\end{corollary}

\begin{remark}\upshape \label{case4 further decay remark}
Decay rates of Oldroyd type model with fractional diffusive stress $(-\Delta)^{\beta}\tau$ is considered in \cite{Wang}, and the results related to (\ref{our model Case IV}) in \cite{Wang} ($\beta=1$) is following: in three dimension, for $(u_{0},\tau_{0})\in L^{1}\cap H^{3}$ with $\left\|(u_{0},\tau_{0})\right\|_{L^{1}\cap H^{3}}\leq \epsilon$,
\eqn \label{Wang result}
\begin{split}
& \left\|u(t)\right\|_{L^{2}}\lesssim (1+t)^{-\frac{3}{4}},\quad \left\|u(t)\right\|_{L^{\infty}}\lesssim (1+t)^{-\frac{3}{2}},\quad\left\|\nabla u(t)\right\|_{L^{2}}\lesssim (1+t)^{-\frac{5}{4}}, \\
& \left\|\nabla u(t)\right\|_{L^{\infty}}\lesssim (1+t)^{-2},\quad \left\|\mathbb{P}\nabla\cdot\tau(t)\right\|_{L^{2}}\lesssim (1+t)^{-\frac{5}{4}}.
\end{split}
\een
In \cite{Wang} the authors assume sufficiently small initial data in $L^{1}\cap H^{3}$, but Corollary \ref{case4 further decay} only need the smallness assumption of $(u_{0},\tau_{0})$ in $\mathcal{E}_{0}$, which is larger than $H^{3}$. Moreover, (\ref{case4 further decay2}) covers (\ref{Wang result}) since $\dot{B}^{k}\hookrightarrow \dot{H}^{k}$ and $\dot{B}^{\frac{3}{2}}\hookrightarrow L^{\infty}$. We also obtain the decay rates of the entire stress tensor $\tau$, not $\mathbb{P}\nabla\cdot\tau$.
\end{remark}

To prove Theorem \ref{case4 theorem} and Theorem \ref{case4 decay}, we begin with the following linearized equations with a divergence-free vector field $v$
\eqn \label{case4 linear eq}
\begin{split} 
& u_{t}+\mathbb{P}(v\cdot\nabla u)-\mathbb{P}\dv \tau=F, \\
& \tau_{t}+v\cdot\nabla\tau-\nu_{2}\Delta\tau-D(u)=G, \\
& (\mathbb{P}\dv \tau)_{t}+\mathbb{P}\dv(v\cdot\nabla\tau)-\nu_{2}\Delta(\mathbb{P}\dv \tau)-\frac{1}{2}\Delta u=\mathbb{P}\dv G.
\end{split}
\een

\begin{lemma} \label{case4 lemma}\upshape
Let $-\frac{d}{2}-1<s_{1}\leq\frac{d}{2}+1$, $-\frac{d}{2}<s_{2}\leq \frac{d}{2}$, $u_{0}\in \dot{B}^{s_{1},s_{2}+1}$, $\tau_{0}\in \dot{B}^{s_{1},s_{2}}$, and $v\in L^{1}_{T}\dot{B}^{\frac{d}{2}+1}$. Then, a solution $(u,\tau)$ of (\ref{case4 linear eq}) satisfies the following bound
\begin{align*}
&\left\|u\right\|_{L^{\infty}_{T}\dot{B}^{s_{1},s_{2}+1}}+\left\|\tau\right\|_{L^{\infty}_{T}\dot{B}^{s_{1},s_{2}}}+\left\|u\right\|_{L^{1}_{T}\dot{B}^{s_{1}+2,s_{2}+1}}+\left\|\tau\right\|_{L^{1}_{T}\dot{B}^{s_{1}+2,s_{2}+2}} \\
& \lesssim \left\|u_{0}\right\|_{\dot{B}^{s_{1},s_{2}+1}}+\left\|\tau_{0}\right\|_{\dot{B}^{s_{1},s_{2}}}+\left\|F\right\|_{L^{1}_{T}\dot{B}^{s_{1},s_{2}+1}}+\left\|G\right\|_{L^{1}_{T}\dot{B}^{s_{1},s_{2}}} +\left\|\mathbb{B}(\nabla v,\nabla\tau)\right\|_{L^{1}_{T}\dot{B}^{s_{1},s_{2}-1}}  \\
&+\left(\left\|u\right\|_{L^{\infty}_{T}\dot{B}^{s_{1},s_{2}+1}}+\left\|\tau\right\|_{L^{\infty}_{T}\dot{B}^{s_{1},s_{2}}}\right)\left\|v\right\|_{L^{1}_{T}\dot{B}^{\frac{d}{2}+1}}.
\end{align*}
\end{lemma}

Let $(s_{1},s_{2})=\left(\frac{d}{2}-1,\frac{d}{2}\right)$, $v=u$, $F=0$ and $-G=\mathcal{Q}(\tau, \nabla u)\simeq \tau\nabla u+(\nabla u)^{2}$. Since
\[
\left\|\mathcal{Q}(\tau,\nabla u)\right\|_{L^{1}_{T}\dot{B}^{\frac{d}{2}-1,\frac{d}{2}}}+\left\|\mathbb{B}(\nabla u,\nabla\tau)\right\|_{L^{1}_{T}\dot{B}^{\frac{d}{2}-1}}\lesssim \left\|(u,\tau)\right\|_{\mathcal{L}_{T}}\left\|(u,\tau)\right\|_{\mathcal{H}_{T}},
\]
Lemma \ref{case4 lemma} provides the bound $\left\|(u,\tau)\right\|_{\mathcal{E}_{T}}\lesssim \left\|(u_{0},\tau_{0})\right\|_{\mathcal{E}_{0}}$ for all $T>0$ when $\left\|(u_{0},\tau_{0})\right\|_{\mathcal{E}_{0}}$ is sufficiently small.

\subsection{Proof of Lemma \ref{case4 lemma}}


\subsubsection{\bf Low frequency part}
First of all, using (\ref{case4 linear eq}), we obtain
\begin{align}
& \frac{1}{2}\frac{d}{dt}\left\|\Delta_{j}u\right\|^{2}_{L^{2}}-(\Delta_{j}u,\Delta_{j}\dv \tau)=(\Delta_{j}F,\Delta_{j}u)-(\Delta_{j}(v\cdot\nabla u),\Delta_{j}u), \label{case4 eq1}\\
& \frac{1}{2}\frac{d}{dt}\left\|\Delta_{j}\tau\right\|^{2}_{L^{2}}+\nu_{2}\left\|\Delta_{j}\Lambda\tau\right\|^{2}_{L^{2}}+(\Delta_{j}u,\Delta_{j}\dv \tau)=(\Delta_{j}G,\Delta_{j}\tau)-(\Delta_{j}(v\cdot\nabla\tau),\Delta_{j}\tau), \label{case4 eq2}
\end{align}
and
\eqn \label{case4 eq3}
\begin{split}
&\frac{d}{dt}(\Delta_{j}u,\Delta_{j}\mathbb{P}\dv \tau)-\left\|\Delta_{j}\mathbb{P}\dv \tau\right\|^{2}_{L^{2}}+\frac{1}{2}\left\|\Delta_{j}\Lambda u\right\|^{2}_{L^{2}}+\nu_{2}(\Delta_{j}u,\Delta_{j}\Lambda^{2}\dv \tau) \\
&=(\Delta_{j}F,\Delta_{j}\mathbb{P}\dv \tau)+(\Delta_{j}\mathbb{P}\dv G,\Delta_{j}u)-(\Delta_{j}(v\cdot\nabla u),\Delta_{j}\mathbb{P}\dv \tau)-(\Delta_{j}\mathbb{P}\dv(v\cdot\nabla\tau),\Delta_{j}u).
\end{split}
\een
By proceeding to add  (\ref{case4 eq1}), (\ref{case4 eq2}), and $K_{1}(\ref{case4 eq3})$, we obtain 
\[
\frac{1}{2}\frac{d}{dt}f^{2}_{j}+h^{2}_{j}= \widetilde{F}^{1}_{j}-\widetilde{F}^{2}_{j}-\widetilde{F}^{3}_{j},\qquad\text{for}\ j\leq j_{0},
\]
where
\begin{align*}
&f^{2}_{j} =\left\|\Delta_{j}u\right\|^{2}_{L^{2}}+\left\|\Delta_{j}\tau\right\|^{2}_{L^{2}}+2K_{1}(\Delta_{j}u,\Delta_{j}\mathbb{P}\dv \tau), \\
&h^{2}_{j} =\nu_{2}\left\|\Delta_{j}\Lambda\tau\right\|^{2}_{L^{2}}-K_{1}\left\|\Delta_{j}\mathbb{P}\dv \tau\right\|^{2}_{L^{2}}+\frac{K_{1}}{2}\left\|\Delta_{j}\Lambda u\right\|^{2}_{L^{2}}+\nu_{2}K_{1}(\Delta_{j}\Lambda u,\Delta_{j}\Lambda\dv \tau),\\
&\widetilde{F}^{1}_{j} = (\Delta_{j}F,\Delta_{j}u)+(\Delta_{j}G,\Delta_{j}\tau)+K_{1}(\Delta_{j}F,\Delta_{j}\mathbb{P}\dv \tau)+K_{1}(\Delta_{j}\mathbb{P}\dv G,\Delta_{j}u), \\
&\widetilde{F}^{2}_{j} =(\Delta_{j}(v\cdot\nabla u),\Delta_{j}u)+(\Delta_{j}(v\cdot\nabla\tau),\Delta_{j}\tau) +K_{1}\left[(\Delta_{j}(v\cdot\nabla u),\Delta_{j}\mathbb{P}\dv \tau)+(\Delta_{j}(v\cdot\nabla\mathbb{P}\dv \tau),\Delta_{j}u)\right], \\
&\widetilde{F}^{3}_{j} = K_{1}(\Delta_{j}\mathbb{B}(\nabla v,\nabla\tau),\Delta_{j}u).
\end{align*}
For a sufficiently small $K_{1}>0$ and for each $j_{0}$, we have 
\[
f^{2}_{j}\simeq \left\|\Delta_{j}u\right\|^{2}_{L^{2}}+\left\|\Delta_{j}\tau\right\|^{2}_{L^{2}}, \quad h^{2}_{j}\simeq \left\|\Delta_{j}\Lambda u\right\|^{2}_{L^{2}}+\left\|\Delta_{j}\Lambda\tau\right\|^{2}_{L^{2}}\simeq 2^{2j}f^{2}_{j}.
\]
Also using (\ref{hybrid Besov embedding}), (\ref{convection term estimate}) and (\ref{convection term estimate 2}), we bound $\widetilde{F}_{j}$ terms
\begin{align*}
\left|\widetilde{F}^{1}_{j}\right|+\left|\widetilde{F}^{3}_{j}\right| &\lesssim c_{j}2^{-s_{1}j}\left(\left\|(F,G)\right\|^{l}_{\dot{B}^{s_{1}}}+\left\|\mathbb{B}(\nabla v,\nabla\tau)\right\|^{l}_{\dot{B}^{s_{1}}}\right)f_{j}, \\
\left|\widetilde{F}^{2}_{j}\right| &\lesssim c_{j}2^{-s_{1}j}\left\|v\right\|_{\dot{B}^{\frac{d}{2}+1}}\left(\left\|u\right\|_{\dot{B}^{s_{1},s_{2}+1}}+\left\|\tau\right\|_{\dot{B}^{s_{1},s_{2}}}\right)
\end{align*}
whenever $j\leq j_{0}$ and $-\frac{d}{2}-1<s_{1}\leq\frac{d}{2}+1$, $-\frac{d}{2}<s_{2}\leq\frac{d}{2}$. From these bounds, we can obtain 
\eqn \label{case4 low estimate}
\begin{split}
\left\|u\right\|^{l}_{L^{\infty}_{T}\dot{B}^{s_{1}}}+\left\|\tau\right\|^{l}_{L^{\infty}_{T}\dot{B}^{s_{1}}}&+\left\|u\right\|^{l}_{L^{1}_{T}\dot{B}^{s_{1}+2}}+\left\|\tau\right\|^{l}_{L^{1}_{T}\dot{B}^{s_{1}+2}} \lesssim \left\|(u_{0}, \tau_{0})\right\|^{l}_{\dot{B}^{s_{1}}}+\left\|(F,G)\right\|^{l}_{L^{1}_{T}\dot{B}^{s_{1}}} \\
&+\left\|\mathbb{B}(\nabla v,\nabla\tau)\right\|^{l}_{L^{1}_{T}\dot{B}^{s_{1}}}+\left(\left\|u\right\|_{L^{\infty}_{T}\dot{B}^{s_{1},s_{2}+1}}+\left\|\tau\right\|_{L^{\infty}_{T}\dot{B}^{s_{1},s_{2}}}\right)\left\|v\right\|_{L^{1}_{T}\dot{B}^{\frac{d}{2}+1}}.
\end{split}
\een

\subsubsection{\bf High frequency part}
Also using (\ref{case4 linear eq}), we obtain
\eqn \label{case4 eq4}
\frac{1}{2}\frac{d}{dt}\left\|\Delta_{j}\Lambda u\right\|^{2}_{L^{2}}-(\Delta_{j}\Lambda u,\Delta_{j}\Lambda\dv \tau)=(\Delta_{j}\Lambda F,\Delta_{j}\Lambda u)-(\Delta_{j}\Lambda(v\cdot\nabla u),\Delta_{j}\Lambda u).
\een
By computing $\nu_{2}(\ref{case4 eq4})+K_{2}(\ref{case4 eq2})+(\ref{case4 eq3})$, we have
\[
\frac{1}{2}\frac{d}{dt}f^{2}_{j}+h^{2}_{j}=\widetilde{G}^{1}_{j}-\widetilde{G}^{2}_{j}-\widetilde{G}^{3}_{j}, \qquad\text{for}\ j\geq j_{0}+1,
\]
where
\begin{align*}
f^{2}_{j} &=\nu_{2}\left\|\Delta_{j}\Lambda u\right\|^{2}_{L^{2}}+K_{2}\left\|\Delta_{j}\tau\right\|^{2}_{L^{2}}+2(\Delta_{j}u,\Delta_{j}\mathbb{P}\dv \tau), \\
h^{2}_{j} &=\nu_{2} K_{2}\left\|\Delta_{j}\Lambda\tau\right\|^{2}_{L^{2}}-\left\|\Delta_{j}\mathbb{P}\dv \tau\right\|^{2}_{L^{2}}+\frac{1}{2}\left\|\Delta_{j}\Lambda u\right\|^{2}_{L^{2}}+K_{2}(\Delta_{j}u,\Delta_{j}\dv \tau),\\
\widetilde{G}^{1}_{j} &=\nu_{2}(\Delta_{j}\Lambda F,\Delta_{j}\Lambda u)+K_{2}(\Delta_{j}G,\Delta_{j}\tau)+(\Delta_{j}F,\Delta_{j}\mathbb{P}\dv \tau)+(\Delta_{j}\mathbb{P}\dv G,\Delta_{j}u), \\
\widetilde{G}^{2}_{j}&=\nu_{2}(\Delta_{j}\Lambda(v\cdot\nabla u),\Delta_{j}\Lambda u)+K_{2}(\Delta_{j}(v\cdot\nabla\tau),\Delta_{j}\tau) +(\Delta_{j}(v\cdot\nabla u),\Delta_{j}\mathbb{P}\dv \tau)\\
&+(\Delta_{j}(v\cdot\nabla\mathbb{P}\dv \tau),\Delta_{j}u), \\
\widetilde{G}^{3}_{j} &=(\Delta_{j}\mathbb{B}(\nabla v,\nabla\tau),\Delta_{j}u).
\end{align*}
By choosing $K_{2}=\frac{2C^{2}_{0}}{\nu_{2}}$ and sufficiently large $j_{0}$, we obtain 
\[
f^{2}_{j}\simeq\left\|\Delta_{j}\Lambda u\right\|^{2}_{L^{2}}+\left\|\Delta_{j}\tau\right\|^{2}_{L^{2}}, \quad h^{2}_{j}\simeq \left\|\Delta_{j}\Lambda u\right\|^{2}_{L^{2}}+\left\|\Delta_{j}\Lambda\tau\right\|^{2}_{L^{2}} \gtrsim f^{2}_{j}.
\]
Using (\ref{convection term estimate}) and (\ref{convection term estimate 2}), we estimate the nonlinear terms as 
\begin{align*}
\left|\widetilde{G}^{1}_{j}\right|+\left|\widetilde{G}^{3}_{j}\right| &\lesssim c_{j}2^{-s_{2}j}\left(\left\|F\right\|^{h}_{\dot{B}^{s_{2}+1}}+\left\|G\right\|^{h}_{\dot{B}^{s_{2}}}+\left\|\mathbb{B}(\nabla v,\nabla \tau)\right\|^{h}_{\dot{B}^{s_{2}-1}}\right)f_{j}, \\
\left|\widetilde{G}^{2}_{j}\right| &\lesssim c_{j}2^{-s_{2}j}\left\|v\right\|_{\dot{B}^{\frac{d}{2}+1}}\left(\left\|u\right\|_{\dot{B}^{s_{1},s_{2}+1}}+\left\|\tau\right\|_{\dot{B}^{s_{1},s_{2}}}\right)f_{j}
\end{align*}
whenever $j\geq j_{0}+1$ and $-\frac{d}{2}-1<s_{1}\leq \frac{d}{2}+1$, $-\frac{d}{2}<s_{2}\leq \frac{d}{2}$. So, we obtain 
\begin{align*}
\left\|u\right\|^{h}_{L^{\infty}_{T}\dot{B}^{s_{2}+1}}&+\left\|\tau\right\|^{h}_{L^{\infty}_{T}\dot{B}^{s_{2}}}+\left\|u\right\|^{h}_{L^{1}_{T}\dot{B}^{s_{2}+1}}+\left\|\tau\right\|^{h}_{L^{1}_{T}\dot{B}^{s_{2}}}  \lesssim \left\|u_{0}\right\|^{h}_{\dot{B}^{s_{2}+1}}+\left\|\tau_{0}\right\|^{h}_{\dot{B}^{s_{2}}}+\left\|F\right\|^{h}_{L^{1}_{T}\dot{B}^{s_{2}+1}} \\
&+\left\|G\right\|^{h}_{L^{1}_{T}\dot{B}^{s_{2}}}+\left\|\mathbb{B}(\nabla v,\nabla \tau)\right\|^{h}_{L^{1}_{T}\dot{B}^{s_{2}-1}}+\left(\left\|u\right\|_{L^{\infty}_{T}\dot{B}^{s_{1},s_{2}+1}}+\left\|\tau\right\|_{L^{\infty}_{T}\dot{B}^{s_{1},s_{2}}}\right)\left\|v\right\|_{L^{1}_{T}\dot{B}^{\frac{d}{2}+1}}.
\end{align*}
Also from (\ref{case4 eq2}), we obtain 
\[
\left\|\tau\right\|^{h}_{L^{\infty}_{T}\dot{B}^{s_{2}}}+\left\|\tau\right\|^{h}_{L^{1}_{T}\dot{B}^{s_{2}+2}} \lesssim \left\|\tau_{0}\right\|^{h}_{\dot{B}^{s_{2}}}+\left\|u\right\|^{h}_{L^{1}_{T}\dot{B}^{s_{2}+1}}+\left\|G\right\|^{h}_{L^{1}_{T}\dot{B}^{s_{2}}}+\left\|\tau\right\|_{L^{\infty}_{T}\dot{B}^{s_{1},s_{2}}}\left\|v\right\|_{L^{1}_{T}\dot{B}^{\frac{d}{2}+1}}.
\]
Combining these two inequalities, we deduce that 
\eqn \label{case4 high estimate}
\begin{split}
\left\|u\right\|^{h}_{L^{\infty}_{T}\dot{B}^{s_{2}+1}}&+\left\|\tau\right\|^{h}_{L^{\infty}_{T}\dot{B}^{s_{2}}}+\left\|u\right\|^{h}_{L^{1}_{T}\dot{B}^{s_{2}+1}}+\left\|\tau\right\|^{h}_{L^{1}_{T}\dot{B}^{s_{2}+2}} \lesssim \left\|u_{0}\right\|^{h}_{\dot{B}^{s_{2}+1}}+\left\|\tau_{0}\right\|^{h}_{\dot{B}^{s_{2}}}+\left\|F\right\|^{h}_{L^{1}_{T}\dot{B}^{s_{2}+1}} \\
&+\left\|G\right\|^{h}_{L^{1}_{T}\dot{B}^{s_{2}}} +\left\|\mathbb{B}(\nabla v,\nabla \tau)\right\|^{h}_{L^{1}_{T}\dot{B}^{s_{2}-1}}+\left(\left\|u\right\|_{L^{\infty}_{T}\dot{B}^{s_{1},s_{2}+1}}+\left\|\tau\right\|_{L^{\infty}_{T}\dot{B}^{s_{1},s_{2}}}\right)\left\|v\right\|_{L^{1}_{T}\dot{B}^{\frac{d}{2}+1}}.
\end{split}
\een

By (\ref{case4 low estimate}) and (\ref{case4 high estimate}), we complete the proof of Lemma \ref{case4 lemma}.

\subsection{Uniqueness}
We assume that $(u_{1},\tau_{1})$ and $(u_{2},\tau_{2})$ are two solutions of (\ref{our model Case IV}) with the same initial data. Let $(u,\tau)= (u_{1}-u_{2},\tau_{1}-\tau_{2})$. Then $(u,\tau)$ satisfies $\dv u_{2}=0$ and
\eqn \label{difference eq}
\begin{split}
& u_{t}+\mathbb{P}(u_{2}\cdot\nabla u)-\mathbb{P}\dv \tau=-\mathbb{P}(u\cdot\nabla u_{1})=\delta F, \\
& \tau_{t}+u_{2}\cdot\nabla\tau-\nu_{2}\Delta\tau-D(u)=-u\cdot\nabla\tau_{1}-\left[\mathcal{Q}(\tau_{1},\nabla u_{1})-\mathcal{Q}(\tau_{2},\nabla u_{2})\right]=\delta G.
\end{split}
\een
We now show the uniqueness of solutions by estimating $(u,\tau)$ in 
\[
\mathcal{F}_{T}=L^{\infty}_{T}\dot{B}^{\frac{d}{2}-1,\frac{d}{2}}\cap L^{1}_{T}\dot{B}^{\frac{d}{2}+1,\frac{d}{2}}\times L^{\infty}_{T}\dot{B}^{\frac{d}{2}-1}\cap L^{1}_{T}\dot{B}^{\frac{d}{2}+1}.
\]
Lemma \ref{case4 lemma} with $(s_{1}, s_{2})=(\frac{d}{2}-1, \frac{d}{2}-1)$ yields 
\eqn \label{case4 difference estimate}
\begin{split}
\left\|(u,\tau)\right\|_{\mathcal{F}_{T}} &\lesssim \left\|\delta F\right\|_{L^{1}_{T}\dot{B}^{\frac{d}{2}-1,\frac{d}{2}}}+\left\|\delta G\right\|_{L^{1}_{T}\dot{B}^{\frac{d}{2}-1}} +\left\|\mathbb{B}(\nabla u_{2},\nabla\tau)\right\|_{L^{1}_{T}\dot{B}^{\frac{d}{2}-1,\frac{d}{2}-2}} \\
&+\Big(\left\|u\right\|_{L^{\infty}_{T}\dot{B}^{\frac{d}{2}-1,\frac{d}{2}}}+\left\|\tau\right\|_{L^{\infty}_{T}\dot{B}^{\frac{d}{2}-1}}\Big) \left\|u_{2}\right\|_{L^{1}_{T}\dot{B}^{\frac{d}{2}+1}}.
\end{split}
\een
Using (\ref{hybrid Besov embedding}) and (\ref{product estimate}) repeatedly, we also estimate the following terms 
\begin{align*}
&\left\|u\cdot\nabla u_{1}\right\|_{L^{1}_{T}\dot{B}^{\frac{d}{2}-1,\frac{d}{2}}} \lesssim \left\|u\right\|_{L^{\infty}_{T}\dot{B}^{\frac{d}{2}-1,\frac{d}{2}}}\left\|u_{1}\right\|_{L^{1}_{T}\dot{B}^{\frac{d}{2}+1}}, \quad \left\|u\cdot\nabla\tau_{1}\right\|_{L^{1}_{T}\dot{B}^{\frac{d}{2}-1}} \lesssim \left\|u\right\|_{L^{\infty}_{T}\dot{B}^{\frac{d}{2}-1}}\left\|\tau_{1}\right\|_{L^{1}_{T}\dot{B}^{\frac{d}{2}+1}}, \\
&\left\|\mathbb{B}(\nabla u_{2},\nabla\tau)\right\|_{L^{1}_{T}\dot{B}^{\frac{d}{2}-1,\frac{d}{2}-2}} \lesssim \left\|\mathbb{B}(\nabla u_{2},\nabla\tau)\right\|_{L^{1}_{T}\dot{B}^{\frac{d}{2}-1}}\lesssim \left\|u_{2}\right\|_{L^{\infty}_{T}\dot{B}^{\frac{d}{2}}}\left\|\tau\right\|_{L^{1}_{T}\dot{B}^{\frac{d}{2}+1}}.
\end{align*}
Since $\mathcal{Q}(\tau_{1},\nabla u_{1})-\mathcal{Q}(\tau_{2},\nabla u_{2}) \simeq \tau\nabla u_{1}+\tau_{2} \nabla u+(\nabla u_{1}+\nabla u_{2})\nabla u$, we finally have 
\begin{align*}
&\left\|\mathcal{Q}(\tau_{1},\nabla u_{1})-\mathcal{Q}(\tau_{2},\nabla u_{2})\right\|_{L^{1}_{T}\dot{B}^{\frac{d}{2}-1}} \lesssim \left\|\tau\right\|_{L^{\infty}_{T}\dot{B}^{\frac{d}{2}-1}}\left\|u_{1}\right\|_{L^{1}_{T}\dot{B}^{\frac{d}{2}+1}} +\left\|\tau_{2}\right\|_{L^{2}_{T}\dot{B}^{\frac{d}{2}}}\left\|u\right\|_{L^{2}_{T}\dot{B}^{\frac{d}{2}}} \\
&+\Big(\left\|u_{1}\right\|_{L^{1}_{T}\dot{B}^{\frac{d}{2}+1}}+\left\|u_{2}\right\|_{L^{1}_{T}\dot{B}^{\frac{d}{2}+1}}\Big)\left\|u\right\|_{L^{\infty}_{T}\dot{B}^{\frac{d}{2}}} \lesssim \left(\left\|(u_{1},\tau_{1})\right\|_{\mathcal{E}_{T}}+\left\|(u_{2},\tau_{2})\right\|_{\mathcal{E}_{T}}\right)\left\|(u,\tau)\right\|_{\mathcal{F}_{T}}.
\end{align*}
So, (\ref{case4 difference estimate}) and these estimates  imply
\[
\left\|(u,\tau)\right\|_{\mathcal{F}_{T}}\lesssim \left(\left\|(u_{1},\tau_{1})\right\|_{\mathcal{E}_{T}}+\left\|(u_{2},\tau_{2})\right\|_{\mathcal{E}_{T}}\right)\left\|(u,\tau)\right\|_{\mathcal{F}_{T}}.
\]
Since $\left\|(u_{1},\tau_{1})\right\|_{\mathcal{E}_{T}}+\left\|(u_{2},\tau_{2})\right\|_{\mathcal{E}_{T}}\lesssim \epsilon$,  $(u,\tau)=0$ in $\mathcal{F}_{T}$. This completes the uniqueness part.

\subsection{Proof of Theorem \ref{case4 decay}}
If we do not take time integrations in the proof of Lemma \ref{case4 lemma} with $(s_{1},s_{2})=(s,\frac{d}{2})$, $v=u$, $F=0$, $-G=\mathcal{Q}(\tau,\nabla u)\simeq \tau\nabla u+(\nabla u)^{2}$, we can obtain the following
\eqn \label{case4 ineq for decay}
\begin{split}
& \frac{d}{dt}\left(\left\|u\right\|_{\dot{B}^{s,\frac{d}{2}+1}}+\left\|\tau\right\|_{\dot{B}^{s,\frac{d}{2}}}\right)+\left\|u\right\|_{\dot{B}^{s+2,\frac{d}{2}+1}}+\left\|\tau\right\|_{\dot{B}^{s+2,\frac{d}{2}+2}} \\
& \lesssim \left\|\mathcal{Q}(\tau,\nabla u)\right\|_{\dot{B}^{s,\frac{d}{2}}}+\left\|\mathbb{B}(\nabla u,\nabla\tau)\right\|_{\dot{B}^{s,\frac{d}{2}-1}}+\left\|u\right\|_{\dot{B}^{\frac{d}{2}+1}}\left(\left\|u\right\|_{\dot{B}^{s,\frac{d}{2}+1}}+\left\|\tau\right\|_{\dot{B}^{s,\frac{d}{2}}}\right).
\end{split}
\een
From (\ref{hybrid Besov embedding}) and (\ref{product estimate}), for $-\frac{d}{2}<s\leq \frac{d}{2}+1$, the nonlinear terms are bounded as 
\[
\begin{split}
&\left\|\mathcal{Q}(\tau,\nabla u)\right\|_{\dot{B}^{s,\frac{d}{2}}} \lesssim \left\|u\right\|_{\dot{B}^{\frac{d}{2}+1}}\left(\left\|u\right\|_{\dot{B}^{s+1,\frac{d}{2}+1}}+\left\|\tau\right\|_{\dot{B}^{s,\frac{d}{2}}}\right)+\left\|\tau\right\|_{\dot{B}^{\frac{d}{2}}}\left\|u\right\|_{\dot{B}^{s+1,\frac{d}{2}+1}},\\
&\left\|\mathbb{B}(\nabla u,\nabla\tau)\right\|_{\dot{B}^{s,\frac{d}{2}-1}} \lesssim\left\|\mathbb{B}(\nabla u,\nabla \tau)\right\|_{\dot{B}^{\min{(s,\frac{d}{2})},\frac{d}{2}-1}}\lesssim \left\|u\right\|_{\dot{B}^{\frac{d}{2}+1}}\left\|\tau\right\|_{\dot{B}^{\min{(s,\frac{d}{2})+1,\frac{d}{2}}}} \lesssim\left\|u\right\|_{\dot{B}^{\frac{d}{2}+1}}\left\|\tau\right\|_{\dot{B}^{s,\frac{d}{2}}}
\end{split}
\]
Thus, using (\ref{hybrid Besov embedding}), we can rewrite (\ref{case4 ineq for decay}) as 
\eqn \label{eq:6.12}
\begin{split}
& \frac{d}{dt}\left(\left\|u\right\|_{\dot{B}^{s,\frac{d}{2}+1}}+\left\|\tau\right\|_{\dot{B}^{s,\frac{d}{2}}}\right)+\left\|u\right\|_{\dot{B}^{s+2,\frac{d}{2}+1}}+\left\|\tau\right\|_{\dot{B}^{s+2,\frac{d}{2}+2}} \\
& \lesssim \left\|u\right\|_{\dot{B}^{\frac{d}{2}+1}}\left(\left\|u\right\|_{\dot{B}^{s,\frac{d}{2}+1}}+\left\|\tau\right\|_{\dot{B}^{s,\frac{d}{2}}}\right)+\left\|\tau\right\|_{\dot{B}^{\frac{d}{2}}}\left\|u\right\|_{\dot{B}^{s+1,\frac{d}{2}+1}} \\
& \lesssim \left\|u\right\|_{\dot{B}^{\frac{d}{2}+1}}\Big(\left\|u\right\|^{h}_{\dot{B}^{\frac{d}{2}+1}}+\left\|\tau\right\|^{h}_{\dot{B}^{\frac{d}{2}}}\Big)+\left\|\tau\right\|_{\dot{B}^{\frac{d}{2}}}\left\|u\right\|^{h}_{\dot{B}^{\frac{d}{2}+1}} +\left\|u\right\|_{\dot{B}^{\frac{d}{2}+1}}\left\|u\right\|^{l}_{\dot{B}^{s}}\\
& +\left\|u\right\|_{\dot{B}^{\frac{d}{2}+1}}\left\|\tau\right\|^{l}_{\dot{B}^{s}}+\left\|\tau\right\|_{\dot{B}^{\frac{d}{2}}}\left\|u\right\|^{l}_{\dot{B}^{s+1}}.
\end{split}
\een
Since $\left\|(u,\tau)\right\|_{\mathcal{E}}\lesssim\epsilon$ and (\ref{bound for decay 1}), for $s\geq \frac{d}{2}-1$ the first three terms of the right-hand side of (\ref{eq:6.12}) can be absorbed into the left-hand side of (\ref{eq:6.12}). As (\ref{bound for decay 1}), we also bound $\left\|u\right\|_{\dot{B}^{\frac{d}{2}+1}}\left\|\tau\right\|^{l}_{\dot{B}^{s}}$ as follows
\begin{align*}
\left\|u\right\|_{\dot{B}^{\frac{d}{2}+1}}\left\|\tau\right\|^{l}_{\dot{B}^{s}} &=\left\|u\right\|^{l}_{\dot{B}^{\frac{d}{2}+1}}\left\|\tau\right\|^{l}_{\dot{B}^{s}}+\left\|u\right\|^{h}_{\dot{B}^{\frac{d}{2}+1}}\left\|\tau\right\|^{l}_{\dot{B}^{s}} \\
& \lesssim \left(\left\|u\right\|^{l}_{\dot{B}^{\frac{d}{2}-1}}\left\|\tau\right\|^{l}_{\dot{B}^{s+2}}\right)^{1-\theta}\left(\left\|\tau\right\|^{l}_{\dot{B}^{\frac{d}{2}-1}}\left\|u\right\|^{l}_{\dot{B}^{s+2}}\right)^{\theta}+\left\|u\right\|^{h}_{\dot{B}^{\frac{d}{2}+1}}\left\|\tau\right\|^{l}_{\dot{B}^{s}} \\
& \lesssim \left\|u\right\|^{l}_{\dot{B}^{\frac{d}{2}-1}}\left\|\tau\right\|^{l}_{\dot{B}^{s+2}}+\left\|\tau\right\|^{l}_{\dot{B}^{\frac{d}{2}-1}}\left\|u\right\|^{l}_{\dot{B}^{s+2}}+\left\|u\right\|^{h}_{\dot{B}^{\frac{d}{2}+1}}\left\|\tau\right\|^{l}_{\dot{B}^{s}} \\
& \lesssim \left(\left\|u\right\|^{l}_{\dot{B}^{\frac{d}{2}-1}}+\left\|\tau\right\|^{l}_{\dot{B}^{\frac{d}{2}-1}}\right)\left(\left\|u\right\|_{\dot{B}^{s+2,\frac{d}{2}+1}}+\left\|\tau\right\|^{l}_{\dot{B}^{s+2}}\right)
\end{align*}
where $0<\theta=\frac{2}{s+3-d/2}\leq1$. Similarly, we can bound $\left\|\tau\right\|_{\dot{B}^{\frac{d}{2}}}\left\|u\right\|^{l}_{\dot{B}^{s+1}}$:
\begin{align*}
\left\|\tau\right\|_{\dot{B}^{\frac{d}{2}}}\left\|u\right\|^{l}_{\dot{B}^{s+1}} &=\left\|\tau\right\|^{l}_{\dot{B}^{\frac{d}{2}}}\left\|u\right\|^{l}_{\dot{B}^{s+1}}+\left\|\tau\right\|^{h}_{\dot{B}^{\frac{d}{2}}}\left\|u\right\|^{l}_{\dot{B}^{s+1}} \\
& \lesssim \left(\left\|u\right\|^{l}_{\dot{B}^{\frac{d}{2}-1}}\left\|\tau\right\|^{l}_{\dot{B}^{s+2}}\right)^{\gamma}\left(\left\|\tau\right\|^{l}_{\dot{B}^{\frac{d}{2}-1}}\left\|u\right\|^{l}_{\dot{B}^{s+2}}\right)^{1-\gamma}+\left\|\tau\right\|^{h}_{\dot{B}^{\frac{d}{2}}}\left\|u\right\|^{l}_{\dot{B}^{s+1}} \\
& \lesssim \left\|u\right\|^{l}_{\dot{B}^{\frac{d}{2}-1}}\left\|\tau\right\|^{l}_{\dot{B}^{s+2}}+\left\|\tau\right\|^{l}_{\dot{B}^{\frac{d}{2}-1}}\left\|u\right\|^{l}_{\dot{B}^{s+2}}+\left\|\tau\right\|^{h}_{\dot{B}^{\frac{d}{2}}}\left\|u\right\|^{l}_{\dot{B}^{s+1}} \\
& \lesssim \left(\left\|u\right\|^{l}_{\dot{B}^{\frac{d}{2}-1}}+\left\|\tau\right\|^{l}_{\dot{B}^{\frac{d}{2}-1}}\right)\left(\left\|u\right\|^{l}_{\dot{B}^{s+2}}+\left\|\tau\right\|_{\dot{B}^{s+2,\frac{d}{2}}}\right)
\end{align*}
where $0<\gamma=\frac{1}{s+3-d/2}<1$. Since all the terms of right-hand side of (\ref{eq:6.12}) can be absorbed into the left-hand side of (\ref{eq:6.12}), we derive 
\eqn \label{eq:6.15}
\frac{d}{dt}\left(\left\|u\right\|_{\dot{B}^{s,\frac{d}{2}+1}}+\left\|\tau\right\|_{\dot{B}^{s,\frac{d}{2}}}\right)+\left\|u\right\|_{\dot{B}^{s+2,\frac{d}{2}+1}}+\left\|\tau\right\|_{\dot{B}^{s+2,\frac{d}{2}+2}}\leq 0
\een
for $\frac{d}{2}-1\leq s\leq \frac{d}{2}+1$. As in the proof of Theorem \ref{case1 decay}, using (\ref{interpolation}) and letting $s=\frac{d}{2}-1+s_{0}$, $s_{0}\in (0,2]$, we get
\[
\frac{d}{dt}\left(\left\|u\right\|_{\dot{B}^{\frac{d}{2}-1+s_{0},\frac{d}{2}+1}}+\left\|\tau\right\|_{\dot{B}^{\frac{d}{2}-1+s_{0},\frac{d}{2}}}\right)+\left(\left\|u\right\|_{\dot{B}^{\frac{d}{2}-1+s_{0},\frac{d}{2}+1}}+\left\|\tau\right\|_{\dot{B}^{\frac{d}{2}-1+s_{0},\frac{d}{2}}}\right)^{1+\frac{2}{s_{0}}} \leq 0.
\]
from which we prove Theorem \ref{case4 decay}.

\subsection{Proof of Corollary \ref{case4 further decay}}
We follow the same argument in the proof of Corollary \ref{case1 further decay}, so we present a brief proof here.

\textbullet\ If $(u_{0},\tau_{0})\in L^{2}$, then using the fact $\left\|(u,\tau)\right\|\lesssim \epsilon$, it is easy to show $(u,\tau)\in L^{\infty}_{T}L^{2}\times L^{\infty}_{T}L^{2}\cap L^{2}_{T}\dot{H}^{1}$ with $\left\|(u,\tau)(t)\right\|_{L^{2}}\leq C(\epsilon,\left\|(u_{0},\tau_{0})\right\|_{L^{2}})$ for all $t$. Thus, (\ref{eq:6.15}) with (\ref{interpolation 2}) gives us
\[
\left\|u(t)\right\|_{\dot{B}^{s,\frac{d}{2}+1}}+\left\|\tau(t)\right\|_{\dot{B}^{s,\frac{d}{2}}}\lesssim (1+t)^{-\frac{s}{2}},\quad s\in(\frac{d}{2}-1,\frac{d}{2}+1],
\]
which implies (\ref{case4 further decay1}) using (\ref{interpolation 2}) again.

\textbullet\ When $(u_{0},\tau_{0})\in\dot{B}^{-\frac{3}{2}}_{2,\infty}$, we use the following two estimates arising in low frequency part:
\eqn \label{case4 Fourier split ineq1}
\begin{split}
& \frac{d}{dt}f_{j}+2^{2j}f_{j}\leq X_{j},\quad j\leq j_{0}, \\
f_{j}\simeq \left\|\Delta_{j}u\right\|_{L^{2}}+\left\|\Delta_{j}\tau\right\|_{L^{2}},\quad &X_{j}\lesssim \left\|\Delta_{j}(u\cdot\nabla u)\right\|_{L^{2}}+\left\|\Delta_{j}(u\cdot\nabla\tau)\right\|_{L^{2}}+\left\|\Delta_{j}\mathcal{Q}(\tau,\nabla u)\right\|_{L^{2}},
\end{split}
\een
and 
\eqn \label{case4 Fourier split ineq2}
\frac{d}{dt}\left\|(u,\tau)\right\|^{l}_{\dot{B}^{0}}+\left\|(u,\tau)\right\|^{l}_{\dot{B}^{2}}\lesssim \left\|\mathcal{Q}(\tau,\nabla u)\right\|^{l}_{\dot{B}^{0}}+\left\|(u,\tau)\right\|_{\dot{B}^{0}}\left\|u\right\|_{\dot{B}^{\frac{5}{2}}}.
\een
Let $r(t)$ and $t_{1}\geq0$ as in proof of Corollary \ref{case1 further decay}, then for $t\geq t_{1}$ we have
\[
-\left\|(u,\tau)\right\|^{l}_{\dot{B}^{2}}\leq -(1+t)^{-1}\left\|(u,\tau)\right\|^{l}_{\dot{B}^{0}}+(1+t)^{-1}\sum_{j\leq r(t)}\left\|(\Delta_{j}u,\Delta_{j}\tau)\right\|_{L^{2}}.
\]
The inequality (\ref{case4 Fourier split ineq2}) arrives at
\[
\frac{d}{dt}\left\|(u,\tau)\right\|^{l}_{\dot{B}^{0}}+(1+t)^{-1}\left\|(u,\tau)\right\|^{l}_{\dot{B}^{0}}\lesssim (1+t)^{-1}\sum_{j\leq r(t)}\left\|(\Delta_{j}u,\Delta_{j}\tau)\right\|_{L^{2}}+\left\|(u,\tau)\right\|_{\dot{B}^{0}}\left\|u\right\|_{\dot{B}^{\frac{5}{2}}}.
\]
Using Duhamel's formula for (\ref{case4 Fourier split ineq1}), we obtain
\[
\sum_{j\leq r(t)}\left\|(\Delta_{j}u,\Delta_{j}\tau)(t)\right\|_{L^{2}}\lesssim \sum_{j\leq r(t)}f_{j}(t)\lesssim (1+t)^{-\frac{3}{4}}\left\|(u_{0},\tau_{0})\right\|^{l}_{\dot{B}^{-\frac{3}{2}}_{2,\infty}}+(1+t)^{-\frac{3}{4}}\int^{t}_{0}\sup_{j\leq j_{0}}2^{-\frac{3}{2}j}X_{j}(t')\, dt'.
\]
Since
\[
\begin{split}
\sup_{j\leq j_{0}}2^{-\frac{3}{2}j}X_{j}(t) &\lesssim \left\|(u\cdot\nabla u)(t)\right\|_{L^{1}}+\left\|(u\cdot\nabla\tau)(t)\right\|_{L^{1}}+\left\|\mathcal{Q}(\tau,\nabla u)(t)\right\|_{L^{1}} \\
&\lesssim \left\|(u,\tau)(t)\right\|_{L^{2}}\left\|(\nabla u,\nabla\tau)(t)\right\|_{L^{2}}+\left\|\nabla u(t)\right\|^{2}_{L^{2}}\lesssim (1+t)^{-\frac{1}{2}},
\end{split}
\]
we get
\[
\sum_{j\leq r(t)}\left\|(\Delta_{j}u,\Delta_{j}\tau)(t)\right\|_{L^{2}}\lesssim (1+t)^{-\frac{1}{4}}.
\]
Thus, we obtain 
\[
\frac{d}{dt}\left[(1+t)\left\|(u,\tau)(t)\right\|^{l}_{\dot{B}^{0}}\right]\lesssim (1+t)^{-\frac{1}{4}},
\]
and then
\[
\left\|(u,\tau)(t)\right\|^{l}_{\dot{B}^{0}}\lesssim (1+t)^{-\frac{1}{4}}.
\]
By repeating the above process from this, we end up obtaining
\[
\left\|(u,\tau)(t)\right\|_{\dot{B}^{0}}\lesssim (1+t)^{-\frac{3}{4}}.
\]
From (\ref{eq:6.15}) and (\ref{interpolation}), we derive (\ref{case4 further decay2}).

\section{Case V: $(\nu_{1}, \nu_{2},\alpha)=(0,0,+)$} \label{case5}
We finally deal with (\ref{our model}) with $(\nu_{1}, \nu_{2},\alpha)=(0,0,+)$:
\eqn \label{our model Case V}
\begin{split} 
& u_{t}+u\cdot\nabla u +\nabla p=\dv \tau, \label{our model a}\\
& \tau_{t}+u\cdot\nabla \tau+\alpha\tau+\mathcal{Q}(\tau, \nabla u)= D(u), \label{our model b}\\
& \dv u=0
\end{split}
\een
with $\mathcal{Q}(\nabla u,\tau)\simeq\tau^{2}$. Let 
\[
\begin{split}
&\left\|(u_{0},\tau_{0})\right\|_{\mathcal{E}_{0}}=\left\|u_{0}\right\|_{\dot{B}^{\frac{d}{2}-1,\frac{d}{2}+1}}+\left\|\tau_{0}\right\|_{\dot{B}^{\frac{d}{2},\frac{d}{2}+1}}, \quad \|(u, \tau)\|_{\mathcal{E}_{T}}=\|(u, \tau)\|_{\mathcal{L}_{T}}+\|(u, \tau)\|_{\mathcal{H}_{T}},\\
&\left\|(u,\tau)\right\|_{\mathcal{L}_{T}}=\left\|u\right\|_{L^{\infty}_{T}\dot{B}^{\frac{d}{2}-1,\frac{d}{2}+1}}+\left\|\tau\right\|_{L^{\infty}_{T}\dot{B}^{\frac{d}{2},\frac{d}{2}+1}}, \quad \left\|(u,\tau)\right\|_{\mathcal{H}_{T}}=\left\|u\right\|_{L^{1}_{T}\dot{B}^{\frac{d}{2}+1}}+\left\|\tau\right\|_{L^{1}_{T}\dot{B}^{\frac{d}{2},\frac{d}{2}+1}}.
\end{split}
\]

\begin{theorem} \label{case5 theorem}\upshape
There exists $\epsilon>0$ such that if $(u_{0},\tau_{0})\in \mathcal{E}_{0}$ with $\left\|(u_{0},\tau_{0})\right\|_{\mathcal{E}_{0}}\leq \epsilon$, there exists a unique solution $(u,\tau)$ in $\mathcal{E}_{T}$ of (\ref{our model Case V}) such that $\left\|(u,\tau)\right\|_{\mathcal{E}_{T}}\lesssim \left\|(u_{0},\tau_{0})\right\|_{\mathcal{E}_{0}}$ for all $T>0$.
\end{theorem}

\begin{theorem} \label{case5 decay}\upshape
The solutions of Theorem \ref{case5 theorem} have the following decay rates:
\[
\left\|u(t)\right\|_{\dot{B}^{\frac{d}{2}-1+s_{0},\frac{d}{2}+1}}+\left\|\tau(t)\right\|_{\dot{B}^{\frac{d}{2}+s_{0},\frac{d}{2}+1}}\leq C(\left\|(u_{0},\tau_{0})\right\|_{\mathcal{E}_{0}},s_{0})(1+t)^{-\frac{s_{0}}{2}}, \quad s_{0}\in(0,2].
\]
\end{theorem}

We also derive the further decay rates of the solution of Theorem \ref{case5 theorem}, but we omit the proof since it is exactly same with that of Corollary \ref{case1 further decay}.

\begin{corollary} \label{case5 further decay}
Let $(u,\tau)\in\mathcal{E}_{T}$ be the solution of Theorem \ref{case5 theorem} with $(u_{0},\tau_{0})\in\mathcal{E}_{0}$.
\item{\textbullet} Assume that $(u_{0},\tau_{0})\in L^{2}$. Then, $(u,\tau)$ is in $L^{\infty}_{T}L^{2}\times L^{\infty}_{T}L^{2}\cap L^{2}_{T}L^{2}$ and satisfies
\begin{align*}
&\left\|u(t)\right\|_{\dot{B}^{s}}\lesssim (1+t)^{-\frac{s}{2}},\quad s\in(0,\frac{d}{2}+1],\\
\left\|\tau(t)\right\|_{L^{2}}\lesssim (1+t)^{-\frac{1}{2}},\qquad &\left\|\tau(t)\right\|_{\dot{B}^{s}}\lesssim 
\begin{cases}
(1+t)^{-\frac{1}{2}-\frac{s}{2}},\quad s\in (0,\frac{d}{2}],\\ 
(1+t)^{-\frac{1}{2}-\frac{d}{4}},\quad s\in (\frac{d}{2},\frac{d}{2}+1].
\end{cases}
\end{align*}
\item{\textbullet} Let $d=3$ and assume that $(u_{0},\tau_{0})\in \dot{B}^{-\frac{3}{2}}_{2,\infty}\times\dot{B}^{-\frac{1}{2}}_{2,\infty}$. Then, $(u,\tau)$ has the further decay rates:
\begin{align*}
&\left\|u(t)\right\|_{\dot{B}^{s}}\lesssim (1+t)^{-\frac{3}{4}-\frac{s}{2}},\quad s\in[0,\frac{5}{2}],\\
&\left\|\tau(t)\right\|_{\dot{B}^{s}}\lesssim
\begin{cases}
(1+t)^{-\frac{5}{4}-\frac{s}{2}},\quad s\in[0,\frac{3}{2}], \\
(1+t)^{-2},\quad s\in (\frac{3}{2},\frac{5}{2}].
\end{cases}
\end{align*}
\end{corollary}

To prove Theorem \ref{case5 theorem} and Theorem \ref{case5 decay}, just as in the previous cases, we use the  linearized equations with a divergence-free vector field $v$:
\eqn \label{case5 linear eq}
\begin{split} 
& u_{t}+\mathbb{P}(v\cdot\nabla u)-\mathbb{P}\dv \tau=F, \\
& \tau_{t}+v\cdot\nabla\tau+\alpha \tau-D(u)=G,\\
& (\mathbb{P}\dv \tau)_{t}+\mathbb{P}\dv(v\cdot\nabla\tau)+\alpha\mathbb{P}\dv \tau-\frac{1}{2}\Delta u=\mathbb{P}\dv G.
\end{split}
\een

\begin{lemma} \label{case5 lemma}\upshape
Let $-\frac{d}{2}-1<s\leq\frac{d}{2}+1$, $u_{0}\in \dot{B}^{s,\frac{d}{2}+1}$, $\tau_{0}\in \dot{B}^{s+1,\frac{d}{2}+1}$, and $v\in L^{1}_{T}\dot{B}^{\frac{d}{2}+1}$. Then, a solution $(u,\tau)$ of (\ref{case5 linear eq}) satisfies the following bound
\begin{align*}
& \left\|u\right\|_{L^{\infty}_{T}\dot{B}^{s,\frac{d}{2}+1}}+\left\|\tau\right\|_{L^{\infty}_{T}\dot{B}^{s+1,\frac{d}{2}+1}}+\left\|u\right\|_{L^{1}_{T}\dot{B}^{s+2,\frac{d}{2}+1}}+\left\|\tau\right\|_{L^{1}_{T}\dot{B}^{s+1,\frac{d}{2}+1}} \\
& \lesssim \left\|u_{0}\right\|_{\dot{B}^{s,\frac{d}{2}+1}}+\left\|\tau_{0}\right\|_{\dot{B}^{s+1,\frac{d}{2}+1}}+\left\|F\right\|_{L^{1}_{T}\dot{B}^{s,\frac{d}{2}+1}}+\left\|G\right\|_{L^{1}_{T}\dot{B}^{s+1,\frac{d}{2}+1}} +\left\|\mathbb{B}(\nabla v,\nabla \tau)\right\|_{L^{1}_{T}\dot{B}^{s,\frac{d}{2}-1}} \\
& +\Big(\left\|u\right\|_{L^{\infty}_{T}\dot{B}^{s,\frac{d}{2}+1}}+\left\|\tau\right\|_{L^{\infty}_{T}\dot{B}^{\min{(s,\frac{d}{2})}+1,\frac{d}{2}+1}}\Big)\left\|v\right\|_{L^{1}_{T}\dot{B}^{\frac{d}{2}+1}}.
\end{align*}
\end{lemma}

Let $s=\frac{d}{2}-1$, $v=u$, $F=0$ and $-G=\mathcal{Q}(\tau,\nabla u)\simeq\tau^{2}$. Since
\[
\left\|\mathcal{Q}(\tau,\nabla u)\right\|_{L^{1}_{T}\dot{B}^{\frac{d}{2},\frac{d}{2}+1}}+\left\|\mathbb{B}(\nabla u,\nabla\tau)\right\|_{L^{1}_{T}\dot{B}^{\frac{d}{2}-1}}\lesssim \left\|(u,\tau)\right\|_{\mathcal{L}_{T}} \left\|(u,\tau)\right\|_{\mathcal{H}_{T}},
\]
we see that Lemma \ref{case5 lemma} implies $\left\|(u,\tau)\right\|_{\mathcal{E}_{T}}\lesssim \left\|(u_{0},\tau_{0})\right\|_{\mathcal{E}_{0}}$ for all $T>0$ when $\left\|(u_{0},\tau_{0})\right\|_{\mathcal{E}_{0}}$ is sufficiently small.

\subsection{Proof of Lemma \ref{case5 lemma}}

\subsubsection{\bf Low frequency part}
From (\ref{case5 linear eq}), we first obtain
\begin{align}
& \frac{1}{2}\frac{d}{dt}\left\|\Delta_{j}u\right\|^{2}_{L^{2}}-(\Delta_{j}u,\Delta_{j}\dv \tau)=(\Delta_{j}F,\Delta_{j}u)-(\Delta_{j}(v\cdot\nabla u),\Delta_{j}u), \label{case5 eq1}\\
& \frac{1}{2}\frac{d}{dt}\left\|\Delta_{j}\Lambda\tau\right\|^{2}_{L^{2}}+\alpha\left\|\Delta_{j}\Lambda\tau\right\|^{2}_{L^{2}}+(\Delta_{j}\Lambda u,\Delta_{j}\Lambda\dv \tau)=(\Delta_{j}\Lambda G,\Delta_{j}\Lambda\tau)-(\Delta_{j}\Lambda(v\cdot\nabla\tau),\Delta_{j}\Lambda\tau), \label{case5 eq2}
\end{align}
and
\eqn \label{case5 eq3}
\begin{split}
&\frac{d}{dt}(\Delta_{j}u,\Delta_{j}\mathbb{P}\dv \tau)-\left\|\Delta_{j}\mathbb{P}\dv \tau\right\|^{2}_{L^{2}}+\frac{1}{2}\left\|\Delta_{j}\Lambda u\right\|^{2}_{L^{2}}+\alpha(\Delta_{j}u,\Delta_{j}\dv \tau) \\
&=(\Delta_{j}F,\Delta_{j}\mathbb{P}\dv \tau)+(\Delta_{j}\mathbb{P}\dv G,\Delta_{j}u)-(\Delta_{j}(v\cdot\nabla u),\Delta_{j}\mathbb{P}\dv \tau)-(\Delta_{j}\mathbb{P}\dv(v\cdot\nabla\tau),\Delta_{j}u).
\end{split}
\een
Then, $\alpha(\ref{case5 eq1})+K_{1}(\ref{case5 eq2})+(\ref{case5 eq3})$ gives us 
\[
\frac{1}{2}\frac{d}{dt}f^{2}_{j}+h^{2}_{j}= \widetilde{F}^{1}_{j}-\widetilde{F}^{2}_{j}-\widetilde{F}^{3}_{j}-\widetilde{F}^{4}_{j},\qquad\text{for}\ j\leq j_{0},
\]
where
\begin{align*}
f^{2}_{j} &=\alpha\left\|\Delta_{j}u\right\|^{2}_{L^{2}}+K_{1}\left\|\Delta_{j}\Lambda\tau\right\|^{2}_{L^{2}}+2(\Delta_{j}u,\Delta_{j}\mathbb{P}\dv \tau), \\
h^{2}_{j} &=\alpha K_{1}\left\|\Delta_{j}\Lambda\tau\right\|^{2}_{L^{2}}-\left\|\Delta_{j}\mathbb{P}\dv \tau\right\|^{2}_{L^{2}}+\frac{1}{2}\left\|\Delta_{j}\Lambda u\right\|^{2}_{L^{2}}+K_{1}(\Delta_{j}\Lambda u,\Delta_{j}\Lambda\dv \tau),\\
\widetilde{F}^{1}_{j} &= \alpha(\Delta_{j}F,\Delta_{j}u)+K_{1}(\Delta_{j}\Lambda G,\Delta_{j}\Lambda\tau)+(\Delta_{j}F,\Delta_{j}\mathbb{P}\dv \tau)+(\Delta_{j}\mathbb{P}\dv G,\Delta_{j}u), \\
\widetilde{F}^{2}_{j} &=\alpha(\Delta_{j}(v\cdot\nabla u),\Delta_{j}u) +\left[(\Delta_{j}(v\cdot\nabla u),\Delta_{j}\mathbb{P}\dv \tau)+(\Delta_{j}(v\cdot\nabla\mathbb{P}\dv \tau),\Delta_{j}u)\right], \\
\widetilde{F}^{3}_{j} &= K_{1}(\Delta_{j}\Lambda(v\cdot\nabla\tau),\Delta_{j}\Lambda\tau), \\
\widetilde{F}^{4}_{j} &= (\Delta_{j}\mathbb{B}(\nabla v,\nabla\tau),\Delta_{j}u).
\end{align*}
By choosing $K_{1}=\frac{2C^{2}_{0}}{\alpha}$ and sufficiently small $j_{0}$, we obtain 
\[
f^{2}_{j}\simeq \left\|\Delta_{j}u\right\|^{2}_{L^{2}}+\left\|\Delta_{j}\Lambda\tau\right\|^{2}_{L^{2}}, \quad h^{2}_{j}\simeq \left\|\Delta_{j}\Lambda u\right\|^{2}_{L^{2}}+\left\|\Delta_{j}\Lambda\tau\right\|^{2}_{L^{2}}\gtrsim 2^{2j}f^{2}_{j}.
\]
By using (\ref{convection term estimate}) and (\ref{convection term estimate 2}), we next estimate the nonlinear terms:
\begin{align*}
&\left|\widetilde{F}^{1}_{j}\right|+\left|\widetilde{F}^{4}_{j}\right| \lesssim c_{j}2^{-sj}\left(\left\|F\right\|^{l}_{\dot{B}^{s}}+\left\|G\right\|^{l}_{\dot{B}^{s+1}}+\left\|\mathbb{B}(\nabla v,\nabla\tau)\right\|^{l}_{\dot{B}^{s}}\right)f_{j}, \\
&\left|\widetilde{F}^{2}_{j}\right| \lesssim c_{j}2^{-sj}\left\|v\right\|_{\dot{B}^{\frac{d}{2}+1}}\left(\left\|u\right\|_{\dot{B}^{s, \frac{d}{2}+1}}+\left\|\tau\right\|_{\dot{B}^{s+1,\frac{d}{2}+1}}\right)f_{j}, \\
&\left|\widetilde{F}^{3}_{j}\right| \lesssim c_{j}2^{-\min{(s,\frac{d}{2})}j}\left\|v\right\|_{\dot{B}^{\frac{d}{2}+1}}\left\|\tau\right\|_{\dot{B}^{\min{(s,\frac{d}{2})}+1,\frac{d}{2}+1}}\left\|\Delta_{j}\Lambda\tau\right\|_{L^{2}} \lesssim c_{j}2^{-sj}\left\|v\right\|_{\dot{B}^{\frac{d}{2}+1}}\left\|\tau\right\|_{\dot{B}^{\min{(s,\frac{d}{2})}+1,\frac{d}{2}+1}}f_{j}
\end{align*}
whenever $j\leq j_{0}$ and $-\frac{d}{2}-1<s\leq\frac{d}{2}+1$. From these bounds, we obtain 
\begin{align*}
\left\|u\right\|^{l}_{L^{\infty}_{T}\dot{B}^{s}}&+\left\|\tau\right\|^{l}_{L^{\infty}_{T}\dot{B}^{s+1}}+\left\|u\right\|^{l}_{L^{1}_{T}\dot{B}^{s+2}}+\left\|\tau\right\|^{l}_{L^{1}_{T}\dot{B}^{s+3}} \lesssim \left\|u_{0}\right\|^{l}_{\dot{B}^{s}}+ \left\|\tau_{0}\right\|^{l}_{\dot{B}^{s+1}}+\left\|F\right\|^{l}_{L^{1}_{T}\dot{B}^{s}} \\
&+\left\|G\right\|^{l}_{L^{1}_{T}\dot{B}^{s+1}}+\left\|\mathbb{B}(\nabla v,\nabla\tau)\right\|^{l}_{L^{1}_{T}\dot{B}^{s}}+\Big(\left\|u\right\|_{L^{\infty}_{T}\dot{B}^{s,\frac{d}{2}+1}}+\left\|\tau\right\|_{L^{\infty}_{T}\dot{B}^{\min{(s,\frac{d}{2})}+1,\frac{d}{2}+1}}\Big)\left\|v\right\|_{L^{1}_{T}\dot{B}^{\frac{d}{2}+1}}.
\end{align*}
From (\ref{case5 eq2}), we also bound $\tau$ as follows
\begin{align*}
 \left\|\tau\right\|^{l}_{L^{\infty}_{T}\dot{B}^{s+1}}+\left\|\tau\right\|^{l}_{L^{1}_{T}\dot{B}^{s+1}} \lesssim \left\|\tau_{0}\right\|^{l}_{\dot{B}^{s+1}}+\left\|u\right\|^{l}_{\dot{B}^{s+2}}+\left\|G\right\|^{l}_{L^{1}_{T}\dot{B}^{s+1}}+\left\|\tau\right\|_{L^{\infty}_{T}\dot{B}^{\min{(s,\frac{d}{2})}+1,\frac{d}{2}+1}}\left\|v\right\|_{L^{1}_{T}\dot{B}^{\frac{d}{2}+1}}.
\end{align*}
These two inequalities and (\ref{hybrid Besov embedding}) allow us to derive 
\eqn \label{case5 low estimate}
\begin{split}
\left\|u\right\|^{l}_{L^{\infty}_{T}\dot{B}^{s}}&+\left\|\tau\right\|^{l}_{L^{\infty}_{T}\dot{B}^{s+1}}+\left\|u\right\|^{l}_{L^{1}_{T}\dot{B}^{s+2}}+\left\|\tau\right\|^{l}_{L^{1}_{T}\dot{B}^{s+1}} \lesssim \left\|u_{0}\right\|^{l}_{\dot{B}^{s}}+ \left\|\tau_{0}\right\|^{l}_{\dot{B}^{s+1}}+\left\|F\right\|^{l}_{L^{1}_{T}\dot{B}^{s}} \\
&+\left\|G\right\|^{l}_{L^{1}_{T}\dot{B}^{s+1}}+\left\|\mathbb{B}(\nabla v,\nabla\tau)\right\|^{l}_{L^{1}_{T}\dot{B}^{s}}+\Big(\left\|u\right\|_{L^{\infty}_{T}\dot{B}^{s,\frac{d}{2}+1}}+\left\|\tau\right\|_{L^{\infty}_{T}\dot{B}^{\min{(s,\frac{d}{2})}+1,\frac{d}{2}+1}}\Big)\left\|v\right\|_{L^{1}_{T}\dot{B}^{\frac{d}{2}+1}}.
\end{split}
\een 
\subsubsection{\bf High frequency part}
From (\ref{case5 linear eq}), we obtain
\eqn \label{case5 eq4}
\frac{1}{2}\frac{d}{dt}\left\|\Delta_{j}\Lambda u\right\|^{2}_{L^{2}}-(\Delta_{j}\Lambda u,\Delta_{j}\Lambda\dv \tau)=(\Delta_{j}\Lambda F,\Delta_{j}\Lambda u)-(\Delta_{j}\Lambda(v\cdot\nabla u),\Delta_{j}\Lambda u).
\een
By using $(\ref{case5 eq4})+(\ref{case5 eq2})+K_{2}(\ref{case5 eq3})$, we have
\[
\frac{1}{2}\frac{d}{dt}f^{2}_{j}+h^{2}_{j}=\widetilde{G}^{1}_{j}-\widetilde{G}^{2}_{j}-\widetilde{G}^{3}_{j}, \qquad\text{for}\ j\geq j_{0}+1,
\]
where
\begin{align*}
f^{2}_{j} &=\left\|\Delta_{j}\Lambda u\right\|^{2}_{L^{2}}+\left\|\Delta_{j}\Lambda\tau\right\|^{2}_{L^{2}}+2K_{2}(\Delta_{j}u,\Delta_{j}\mathbb{P}\dv \tau), \\
h^{2}_{j} &=\alpha \left\|\Delta_{j}\Lambda\tau\right\|^{2}_{L^{2}}-K_{2}\left\|\Delta_{j}\mathbb{P}\dv \tau\right\|^{2}_{L^{2}}+\frac{K_{2}}{2}\left\|\Delta_{j}\Lambda u\right\|^{2}_{L^{2}}+\alpha K_{2}(\Delta_{j}u,\Delta_{j}\dv \tau),
\\
\widetilde{G}^{1}_{j} &=(\Delta_{j}\Lambda F,\Delta_{j}\Lambda u)+(\Delta_{j}\Lambda G,\Delta_{j}\Lambda\tau)+K_{2}(\Delta_{j}F,\Delta_{j}\mathbb{P}\dv \tau)+K_{2}(\Delta_{j}\mathbb{P}\dv G,\Delta_{j}u), \\
\widetilde{G}^{2}_{j}&=(\Delta_{j}\Lambda(v\cdot\nabla u),\Delta_{j}\Lambda u)+(\Delta_{j}\Lambda(v\cdot\nabla\tau),\Delta_{j}\Lambda\tau) +K_{2}(\Delta_{j}(v\cdot\nabla u),\Delta_{j}\mathbb{P}\dv \tau)\\
& +K_{2}(\Delta_{j}(v\cdot\nabla\mathbb{P}\dv \tau),\Delta_{j}u)\\
\widetilde{G}^{3}_{j} &=K_{2}(\Delta_{j}\mathbb{B}(\nabla v,\nabla\tau),\Delta_{j}u).
\end{align*}
For a sufficiently small $K_{2}>0$ and $j_{0}$ fixed above, we deduce that  
\[
f^{2}_{j}\simeq\left\|\Delta_{j}\Lambda u\right\|^{2}_{L^{2}}+\left\|\Delta_{j}\Lambda\tau\right\|^{2}_{L^{2}}, \quad h^{2}_{j}\simeq \left\|\Delta_{j}\Lambda u\right\|^{2}_{L^{2}}+\left\|\Delta_{j}\Lambda\tau\right\|^{2}_{L^{2}} \simeq f^{2}_{j}.
\]
Also using (\ref{hybrid Besov embedding}), (\ref{convection term estimate}) and (\ref{convection term estimate 2}), we estimate $\widetilde{G}_{j}$ as 
\begin{align*}
\left|\widetilde{G}^{1}_{j}\right|+\left|\widetilde{G}^{3}_{j}\right| &\lesssim c_{j}2^{-\frac{d}{2}j}\Big(\left\|(F,G)\right\|^{h}_{\dot{B}^{\frac{d}{2}+1}}+\left\|\mathbb{B}(\nabla v,\nabla \tau)\right\|^{h}_{\dot{B}^{\frac{d}{2}-1}}\Big)f_{j}, \\
\left|\widetilde{G}^{2}_{j}\right| &\lesssim c_{j}2^{-\frac{d}{2}j}\left\|v\right\|_{\dot{B}^{\frac{d}{2}+1}}\Big(\left\|u\right\|_{\dot{B}^{s,\frac{d}{2}+1}}+\left\|\tau\right\|_{\dot{B}^{\frac{d}{2}+1}}\Big)f_{j}
\end{align*}
whenever $j\geq j_{0}+1$ and $-\frac{d}{2}-1<s\leq \frac{d}{2}+1$. From these bounds, we get 
\eqn
\begin{split} \label{case5 high estimate}
\left\|u\right\|^{h}_{L^{\infty}_{T}\dot{B}^{\frac{d}{2}+1}}&+\left\|\tau\right\|^{h}_{L^{\infty}_{T}\dot{B}^{\frac{d}{2}+1}}+\left\|u\right\|^{h}_{L^{1}_{T}\dot{B}^{\frac{d}{2}+1}}+\left\|\tau\right\|^{h}_{L^{1}_{T}\dot{B}^{\frac{d}{2}+1}} \lesssim \left\|(u_{0},\tau_{0})\right\|^{h}_{\dot{B}^{\frac{d}{2}+1}}+\left\|(F,G)\right\|^{h}_{L^{1}_{T}\dot{B}^{\frac{d}{2}+1}} \\
& +\left\|\mathbb{B}(\nabla v,\nabla \tau)\right\|^{h}_{L^{1}_{T}\dot{B}^{\frac{d}{2}-1}}+\Big(\left\|u\right\|_{L^{\infty}_{T}\dot{B}^{s,\frac{d}{2}+1}}+\left\|\tau\right\|_{L^{\infty}_{T}\dot{B}^{\frac{d}{2}+1}}\Big)\left\|v\right\|_{L^{1}_{T}\dot{B}^{\frac{d}{2}+1}}.
\end{split}
\een
By (\ref{case5 low estimate}) and (\ref{case5 high estimate}), we finish  the proof of Lemma \ref{case5 lemma}.

\subsection{Uniqueness}
We assume that $(u_{1},\tau_{1})$ and $(u_{2},\tau_{2})$ are two solutions of (\ref{our model Case V}) with the same initial data. Let $(u,\tau)=(u_{1}-u_{2},\tau_{1}-\tau_{2})$. Then $(u,\tau)$ satisfies $\dv u_{2}=0$ and 
\eqn
\begin{split} \label{case5 difference eq}
& u_{t}+\mathbb{P}(u_{2}\cdot\nabla u)-\mathbb{P}\dv\tau=-\mathbb{P}(u\cdot\nabla u_{1})=\delta F, \\
& \tau_{t}+u_{2}\cdot\nabla\tau+\alpha\tau-D(u)=-u\cdot\nabla\tau_{1}-\left[\mathcal{Q}(\tau_{1},\nabla u_{1})-\mathcal{Q}(\tau_{2},\nabla u_{2})\right]=\delta G,
\end{split}
\een
where $\mathcal{Q}(\tau_{1},\nabla u_{1})-\mathcal{Q}(\tau_{2},\nabla u_{2})\simeq(\tau_{1}+\tau_{2})\tau$. We now show the uniqueness  in  $
\mathcal{F}_{T}=L^{\infty}_{T}\dot{B}^{\frac{d}{2}}\times L^{\infty}_{T}\dot{B}^{\frac{d}{2}}$. From (\ref{case5 difference eq}), we get
\begin{align*}
& \frac{1}{2}\frac{d}{dt}\left\|\Delta_{j}u\right\|^{2}_{L^{2}}-(\Delta_{j}\dv\tau,\Delta_{j}u)=(\Delta_{j}\delta F,\Delta_{j} u)-(\Delta_{j}(u_{2}\cdot\nabla u),\Delta_{j}u), \\
& \frac{1}{2}\frac{d}{dt}\left\|\Delta_{j}\tau\right\|^{2}_{L^{2}}+\alpha\left\|\Delta_{j}\tau\right\|^{2}_{L^{2}}-(\Delta_{j}D(u),\Delta_{j}\tau)=(\Delta_{j}\delta G,\Delta_{j}\tau)-(\Delta_{j}(u_{2}\cdot\nabla\tau),\Delta_{j}\tau).
\end{align*}
By bounding the last terms using (\ref{convection term estimate}), we obtain 
\begin{align*}
&\frac{1}{2}\frac{d}{dt}\left(\left\|\Delta_{j}u\right\|^{2}_{L^{2}}+\left\|\Delta_{j}\tau\right\|^{2}_{L^{2}}\right)+\alpha\left\|\Delta_{j}\tau\right\|^{2}_{L^{2}} \\
&\lesssim \Big[\left\|\Delta_{j}\delta F\right\|_{L^{2}}+\left\|\Delta_{j}\delta G\right\|_{L^{2}} +c_{j}2^{-\frac{d}{2}j}\left\|u_{2}\right\|_{\dot{B}^{\frac{d}{2}+1}}\Big(\left\|u\right\|_{\dot{B}^{\frac{d}{2}}}+\left\|\tau\right\|_{\dot{B}^{\frac{d}{2}}}\Big)\Big]\left(\left\|\Delta_{j}u\right\|_{L^{2}}+\left\|\Delta_{j}\tau\right\|_{L^{2}}\right).
\end{align*}
From this, we derive the following bound 
\begin{align*}
\left\|(u,\tau)\right\|_{\mathcal{F}_{T}}&\lesssim \left\|\delta F\right\|_{L^{1}_{T}\dot{B}^{\frac{d}{2}}}+\left\|\delta G\right\|_{L^{1}_{T}\dot{B}^{\frac{d}{2}}}+\left\|u_{2}\right\|_{L^{1}_{T}\dot{B}^{\frac{d}{2}+1}} \left\|(u,\tau)\right\|_{\mathcal{F}_{T}}\\
& \lesssim \left\|u\cdot\nabla u_{1}\right\|_{L^{1}_{T}\dot{B}^{\frac{d}{2}}}+\left\|u\cdot\nabla\tau_{1}\right\|_{L^{1}_{T}\dot{B}^{\frac{d}{2}}}+\left\|(\tau_{1}+\tau_{2})\tau\right\|_{L^{1}_{T}\dot{B}^{\frac{d}{2}}}+\left\|u_{2}\right\|_{L^{1}_{T}\dot{B}^{\frac{d}{2}+1}} \left\|(u,\tau)\right\|_{\mathcal{F}_{T}}\\
&\lesssim \Big(\left\|u_{1}\right\|_{L^{1}_{T}\dot{B}^{\frac{d}{2}+1}}+\left\|u_{2}\right\|_{L^{1}_{T}\dot{B}^{\frac{d}{2}+1}}+\left\|\tau_{1}\right\|_{L^{1}_{T}\dot{B}^{\frac{d}{2}+1}}+\left\|\tau_{1}\right\|_{L^{1}_{T}\dot{B}^{\frac{d}{2}}}+\left\|\tau_{2}\right\|_{L^{1}_{T}\dot{B}^{\frac{d}{2}}}\Big)\left\|(u,\tau)\right\|_{\mathcal{F}_{T}}\\
&\lesssim \left(\left\|(u_{1},\tau_{1})\right\|_{\mathcal{E}_{T}}+\left\|(u_{2},\tau_{2})\right\|_{\mathcal{E}_{T}}\right)\left\|(u,\tau)\right\|_{\mathcal{F}_{T}}.
\end{align*}
Since $\left\|(u_{1},\tau_{1})\right\|_{\mathcal{E}_{T}}+\left\|(u_{2},\tau_{2})\right\|_{\mathcal{E}_{T}}\lesssim\epsilon$, $(u,\tau)=0$ in $\mathcal{F}_{T}$ which completes the uniqueness part.

\subsection{Proof of Theorem \ref{case5 decay}}
If we do not take time integrations in the proof of Lemma \ref{case5 lemma} with $v=u$, $F=0$, $-G=\mathcal{Q}(\tau,\nabla u)\simeq\tau^{2}$, we can obtain the following
\[
\begin{split}
& \frac{d}{dt}\left(\left\|u\right\|_{\dot{B}^{s,\frac{d}{2}+1}}+\left\|\tau\right\|_{\dot{B}^{s+1,\frac{d}{2}+1}}\right)+\left\|u\right\|_{\dot{B}^{s+2,\frac{d}{2}+1}}+\left\|\tau\right\|_{\dot{B}^{s+1,\frac{d}{2}+1}} \\
& \lesssim \left\|\mathcal{Q}(\tau,\nabla u)\right\|_{\dot{B}^{s+1,\frac{d}{2}+1}}+\left\|\mathbb{B}(\nabla u,\nabla\tau)\right\|_{\dot{B}^{s,\frac{d}{2}-1}}+\left\|u\right\|_{\dot{B}^{\frac{d}{2}+1}}\left(\left\|u\right\|_{\dot{B}^{s,\frac{d}{2}+1}}+\left\|\tau\right\|_{\dot{B}^{\min{(s,\frac{d}{2})}+1,\frac{d}{2}+1}}\right)\\
&\lesssim \left\|\tau\right\|_{\dot{B}^{\frac{d}{2}}}\left\|\tau\right\|_{\dot{B}^{s+1,\frac{d}{2}+1}}+\left\|u\right\|_{\dot{B}^{\frac{d}{2}+1}}\left(\left\|u\right\|_{\dot{B}^{s,\frac{d}{2}+1}}+\left\|\tau\right\|_{\dot{B}^{\min{(s,\frac{d}{2})}+1,\frac{d}{2}+1}}\right) \\
&\lesssim \left(\left\|u\right\|_{\dot{B}^{\frac{d}{2}+1}}+\left\|\tau\right\|_{\dot{B}^{\frac{d}{2}}}\right)\left(\left\|u\right\|^{h}_{\dot{B}^{\frac{d}{2}+1}}+\left\|\tau\right\|_{\dot{B}^{s+1,\frac{d}{2}+1}}\right)+\left\|u\right\|_{\dot{B}^{\frac{d}{2}+1}}\left\|u\right\|^{l}_{\dot{B}^{s}}+\left\|u\right\|_{\dot{B}^{\frac{d}{2}+1}}\left\|\tau\right\|^{l}_{\dot{B}^{\frac{d}{2}+1}}.
\end{split}
\]
Using (\ref{bound for decay 1}) and (\ref{bound for decay 2}), we have for $\frac{d}{2}-1\leq s\leq\frac{d}{2}+1$
\begin{align*}
& \frac{d}{dt}\left(\left\|u\right\|_{\dot{B}^{s,\frac{d}{2}+1}}+\left\|\tau\right\|_{\dot{B}^{s+1,\frac{d}{2}+1}}\right)+\left\|u\right\|_{\dot{B}^{s+2,\frac{d}{2}+1}}+\left\|\tau\right\|_{\dot{B}^{s+1,\frac{d}{2}+1}} \\
& \lesssim \left(\left\|u\right\|_{\dot{B}^{\frac{d}{2}-1,\frac{d}{2}+1}}+\left\|\tau\right\|_{\dot{B}^{\frac{d}{2}}}\right)\left(\left\|u\right\|_{\dot{B}^{s+2,\frac{d}{2}+1}}+\left\|\tau\right\|_{\dot{B}^{s+1,\frac{d}{2}+1}}\right)  \lesssim \left\|(u,\tau)\right\|_{\mathcal{E}}\left(\left\|u\right\|_{\dot{B}^{s+2,\frac{d}{2}+1}}+\left\|\tau\right\|_{\dot{B}^{s+1,\frac{d}{2}+1}}\right).
\end{align*}
Since $\left\|(u,\tau)\right\|_{\mathcal{E}}\lesssim\epsilon$, we obtain
\[
 \frac{d}{dt}\left(\left\|u\right\|_{\dot{B}^{s,\frac{d}{2}+1}}+\left\|\tau\right\|_{\dot{B}^{s+1,\frac{d}{2}+1}}\right)+\left\|u\right\|_{\dot{B}^{s+2,\frac{d}{2}+1}}+\left\|\tau\right\|_{\dot{B}^{s+1,\frac{d}{2}+1}}\leq0.
\]
As in the proof of Theorem \ref{case1 decay}, using (\ref{interpolation}) and letting $s=\frac{d}{2}-1+s_{0}$, $s_{0}\in (0,2]$, we get
\[
 \frac{d}{dt}\left(\left\|u\right\|_{\dot{B}^{\frac{d}{2}-1+s_{0},\frac{d}{2}+1}}+\left\|\tau\right\|_{\dot{B}^{\frac{d}{2}+s_{0},\frac{d}{2}+1}}\right)+\left(\left\|u\right\|_{\dot{B}^{\frac{d}{2}-1+s_{0},\frac{d}{2}+1}}+\left\|\tau\right\|_{\dot{B}^{\frac{d}{2}+s_{0},\frac{d}{2}+1}}\right)^{1+\frac{2}{s_{0}}} \leq0,
\]
which implies Theorem \ref{case5 decay}.


\section*{Acknowledgments}
H. Bae and J. Shin were supported by the National Research Foundation of Korea(NRF) grant funded by the Korea government(MSIT) (grant No. 2022R1A4A1032094).

\end{document}